\documentclass[12pt]{article}
\usepackage{amsmath,amsthm,verbatim,amssymb,amsfonts,amscd,diagrams, graphics, mathrsfs, graphicx,hyperref}
\theoremstyle{plain}
\newtheorem{theorem}{Theorem}[section]

\newtheorem{lemma}[theorem]{Lemma}
\newtheorem{proposition}[theorem]{Proposition}

\theoremstyle{definition}
\newtheorem{definition}[theorem]{Definition}

\theoremstyle{remark}
\newtheorem{remark}[theorem]{Remark}
\newtheorem{example}[theorem]{Example}

\newarrow{ul}---->
\newarrow{Backwards}<----

\newcommand{\td}[1]{\tilde{#1}}
\newcommand{\into}{\hookrightarrow}

\newcommand{\R}{\mathbb{R}}

\newcommand{\bd}{\partial}

\newcommand{\mbf}[1]{\mathbf{#1}}
\newcommand{\mc}[1]{\mathcal{#1}}
\newcommand{\ms}[1]{\mathscr{#1}}

\newcommand{\Hom}{\text{Hom}}

\newcommand{\Homs}{\textbf{Hom}}

\newarrow{Onto}----{>>}
\newarrow{Equals}=====
\newarrow{Into}C--->

\newcommand{\Set}{\text{\textbf{Set}}}
\newcommand{\Top}{\text{\textbf{Top}}}
\newcommand{\Cat}{\mathcal{C}}

\makeatletter
\def\foo#1\endgraf\unskip#2\foo{\def\row@to@buffer{#1\endgraf\unskip\unskip#2}}
\expandafter\foo\row@to@buffer\foo
\makeatother

\begin{document}

\title{An elementary illustrated introduction to simplicial sets}
\author{Greg Friedman\\Texas Christian University}
\date{December 6, 2011 (minor corrections August 13, 2015; October 3, 2016; December 21, 2020; May 25, 2021; June 10, 2022; and June 8, 2023 --- see errata at end of paper)}

\maketitle
 


\begin{abstract}
This is an expository introduction to simplicial sets and simplicial homotopy theory with particular focus on relating the combinatorial aspects of the theory to their geometric/topological origins. It is intended to be accessible to students familiar with just the fundamentals of algebraic topology.  
\end{abstract}

\tableofcontents

\section{Introduction}

The following notes grew out of my own difficulties in attempting to learn the basics of simplicial sets and simplicial homotopy theory, and thus they are aimed at someone with roughly the same starting knowledge I had, specifically some amount of comfort with simplicial homology  and the basic fundamentals of  topological homotopy theory, including  homotopy groups. Equipped with this background, I wanted to understand a little of what simplicial sets and their generalizations to other categories are all about, as they seem ubiquitous in the literature of certain schools of topology. To name just a few important instances of which I am aware, 
simplicial objects occur in May's work on recognition principles for iterated loop spaces \cite{May72}, Quillen's approach to rational homotopy theory (see \cite{Quil69, FHT}), Bousfield and Kan's work on completions, localization, and limits in homotopy theory \cite{BouKan}, Quillen's abstract treatment of homotopy theory \cite{Quil67}, and various aspects of homological algebra, including group cohomology, Hochschild homology, and cyclic homology (see \cite{WEIB}).

 However, in attempting to learn the rudiments of simplicial  theory, I encountered immediate and discouraging difficulties, which led to serious frustration on several occasions. It was only after several different attempts from different angles that I finally began to ``see the picture,'' and my intended goal here is to aid future students (of all ages) to ease into the subject.

My initial difficulty with the classic expository sources such as May \cite{MAY67} and Curtis \cite{Cu71} was the extent to which the theory is presented purely combinatorially. And the combinatorial definitions are not often pretty; they tend to consist of long strings of axiomatic conditions (see, for example, the combinatorial definition of simplicial homotopy,   Definition \ref{D: comb hom}, below). Despite simplicial objects originating in very topological settings, these classic expositions often sweep this fact too far under the rug for my taste, as someone who likes to comprehend even algebraic and combinatorial constructions as visually as possible. There is a little bit more geometry in Moore's lecture notes \cite{MOORE}, though still not much, and these are also more difficult to obtain (at least not without some good help from a solid Interlibrary Loan Department). On the other hand, there is a much more modern point of view that sweeps both topology and combinatorics away in favor of axiomatic category theory! Goerss and Jardine \cite{GoeJar} is an excellent modern text based upon this approach, which, ironically, helped me tremendously to understand what the combinatorics were getting at!

So what are we getting at here? My goal, still as someone very far from an expert in either combinatorial or axiomatic simplicial theory, is to revisit the material covered in, roughly, the first chapters (in some cases the first few pages) of the texts cited above and to provide some concrete geometric signposts. Here, for the most part, you won't find many complete proofs of theorems, and so these notes will not be completely self-contained. Rather, I try primarily to show by example how the very basic combinatorics, including the definitions, arise out of geometric ideas and to show the geometric ideas underlying the most elementary proofs and properties. Think of this as an appendix or a set of footnotes  to the first chapters of the classic expositions, or perhaps as a Chapter 0. This may not sound like much, but during my earliest learning stages with this material, I would have been very grateful for something of the sort. Theoretically my reader will acquire  enough of ``the idea'' to go forth and read the more thorough (and more technical) sources equipped with enough intuition to \emph{see} what's going on. 

In Section \ref{S: simplicial complexes}, we lay the groundwork with a look at the more familiar topics of simplicial sets and, their slight generalizations, Delta sets. Simplicial sets are then introduced in Section \ref{S: simplicial sets}, followed by their geometric realizations in Section \ref{S: realization} and a detailed look at products of simplicial sets in Section \ref{S: product}. In Section \ref{S: simplicial cats}, we provide a brief look at how the notion of simplicial sets is generalized to other kinds of simplicial objects based in different categories. In Section \ref{S: Kan}, we introduce Kan complexes; these are the simplicial sets that lend themselves to simplicial analogues of homotopy theory, which we study in Section \ref{S: simp homotopy}. This section gets a bit more technical as we head toward more serious applications and theorems in simplicial theory, including the definition and properties  of the simplicial homotopy groups $\pi_n(X,*)$ in Section \ref{S: pin}. Finally, in Section \ref{S: concluding}, we make some concluding remarks and steer the reader toward more comprehensive expository sources.

\paragraph{Acknowledments.} I thank Jim McClure for his useful suggestions and Efton Park for his careful reading of and comments on the preliminary manuscript. Later corrections and improvements were suggested by Henry Adams, Daniel M\"{u}llner, Peter Landweber, and an anonymous referee. I am very grateful for the amount of attention this exposition has received since its initial posting at \texttt{arxiv.org}. 

One text diagram in this paper was typeset using the \TeX\, commutative
diagrams package by Paul Taylor.

\section{A build-up to simplicial sets}\label{S: simplicial complexes}

We begin at the beginning with the relevant geometric notions and their immediate combinatorial counterparts.

\subsection{Simplicial complexes}
Simplicial sets are, essentially, generalizations of the geometric simplicial complexes of elementary algebraic topology (in some cases quite extreme generalizations). So let's recall simplicial complexes, referring the absolute beginner to \cite{MK} for a complete course in the essentials.

Recall that a \emph{(geometric) $n$-simplex} is the convex set spanned by $n+1$ geometrically independent points $\{v_0,\ldots, v_n\}$ in some euclidean space. Here ``geometrically independent'' means that the collection of $n$ vectors $v_1-v_0,\ldots, v_n-v_0$ is linearly independent, and this implies that an $n$-simplex is homeomorphic to a closed $n$-dimensional ball. The points $v_i$ are called \emph{vertices}. A \emph{face} of the (geometric) $n$-simplex determined by $\{v_0,\ldots, v_n\}$ is the convex set spanned by some subset of these vertices. 

A \emph{(geometric) simplicial complex $X$ in $\R^N$} consists of a collection of simplices, possibly of various dimensions, in $\R^N$ such that \begin{enumerate}
\item every face of a simplex of $X$ is in $X$, and
\item the intersection of any two simplices of $X$, if non-empty, is a face of each them.
\end{enumerate}
This definition can be extended easily to handle geometric simplicial complexes containing collections of simplices of arbitrary cardinality and $n$-simplices for arbitrary non-negative integer $n$. Since we will head directly toward abstractions that will obviate this issue by other means, we refer the interested reader to \cite[Section 2]{MK}. We also observe that one is often interested in a geometric simplicial complex only for its homeomorphism type and its combinatorial information, in which case one tends to ignore the precise embedding into euclidean space. This will be the sense in which we shall generally think of simplicial complexes.

\begin{figure}[!htp]
\begin{center}
\scalebox{.6}{\includegraphics{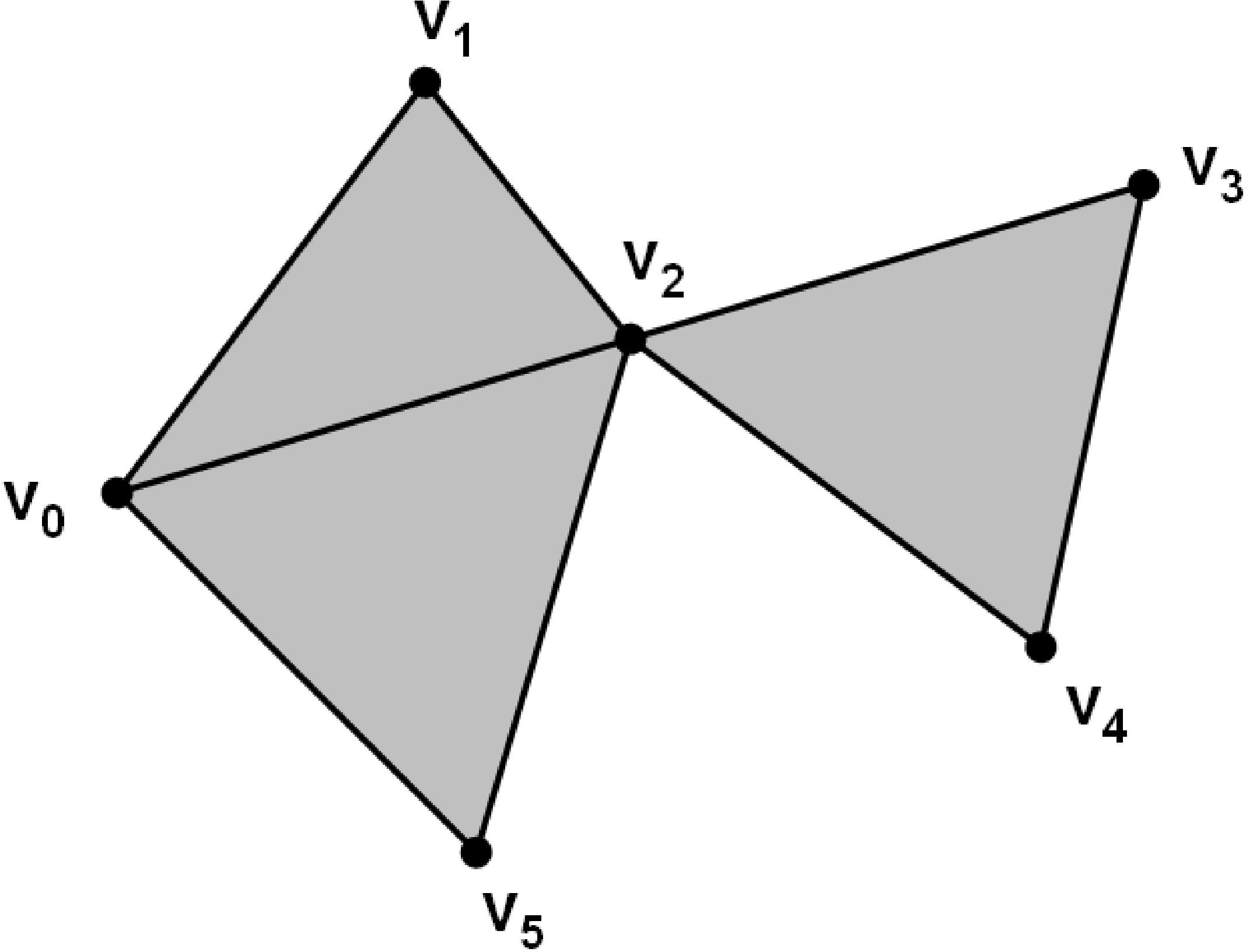}}
\end{center}
\caption{A simplicial complex. Note that  $[v_0,v_1,v_2]$ is a simplex, but $[v_1,v_2,v_4]$ is not.}\label{F: fig1}
\end{figure}

So, less formally, we think of a simplicial complex $X$ as made up of simplices (generalized tetrahedra) of various dimensions, glued together along common faces (see Figure \ref{F: fig1}). The most efficient description, containing all of the relevant information, comes from labeling the vertices (the $0$-simplices) and then specifying which collections of vertices together constitute the vertices of simplices of higher dimension. 
If  the collection of vertices is countable,  we can label them $v_0, v_1,v_2,\ldots$, though this assumption is not strictly necessary - we could label by $\{v_i\}_{i\in I}$ for any indexing set $I$. Then if some collection of vertices $\{v_{i_0},\ldots, v_{i_n}\}$ constitutes the vertices of a simplex, we can label that simplex as $[v_{i_0},\ldots, v_{i_n}]$. 

\begin{example}\label{E: simplex/complex}
If $X$ is a complex and $[v_{i_0},\ldots, v_{i_k}]$ is a simplex of $X$, then any subset of $\{v_{i_0},\ldots, v_{i_k}\}$ is a face of that simplex and thus itself a simplex of $X$. In particular, we can think of the $k$-simplex $[v_{i_0},\ldots, v_{i_k}]$  as a  \emph{geometric simplicial complex} consisting of itself and its faces. 
\end{example}

A nice way to organize the combinatorial information involved is to define the skeleta $X^k$, $k=0,1,\ldots$, of a simplicial complex so that $X^k$ is the set of all $k$-simplices of $X$. Notice that, having labeled our vertices so that  $X^0=\{v_i\}_{i\in I}$, we can think of each element of $X^k$ as a certain subset of $X^0$ of cardinality $k+1$. A subset $\{v_{i_0},\ldots, v_{i_k}\}\subset X^0$ is an \emph{element} of $X^k$ precisely if $[v_{i_0},\ldots, v_{i_k}]$ is a $k$-dimensional simplex of $X$.

To describe a geometric simplicial complex given its set of vertices, it is enough to know which collections of vertices $\{v_{i_0},\ldots, v_{i_k}\}$ correspond to simplices $[v_{i_0},\ldots, v_{i_k}]$ of the simplicial complex. Paring down to this information (which is purely combinatorial) 
leads us  to the notion of an \emph{abstract simplicial complex}. 

\begin{definition}
An \emph{abstract simplicial complex} consists of a set of ``vertices'' $X^0$ together with, for each integer $k$, a set $X^k$ consisting of subsets\footnote{Not necessarily all of them!} of $X^0$ of cardinality $k+1$. These must satisfy the condition that any $(j+1)$-element subset of an element of $X^k$ is an element of $X^j$. 
\end{definition}

Each element of $X^k$ is an abstract $k$-simplex, and 
the last requirement of the definition just guarantees that every face of an abstract simplex in an abstract simplicial complex is also a simplex of the simplicial complex. 

So, an abstract simplicial complex has exactly the same combinatorial information as a geometric simplicial complex. We have lost geometric information about how big a simplex is, how it is embedded in euclidean space, etc., but we have retained all of the information necessary to reconstruct the complex up to homeomorphism.
 It is straightforward that a geometric simplicial complex yields an abstract simplicial complex, but conversely, we can obtain a geometric simplicial complex (up to homeomorphism) from an abstract one by assigning to each element of $X^0$ a point  and to each  abstract simplex $[v_{i_0},\ldots, v_{i_k}]$ a geometric $k$-simplex spanned by the appropriate vertices and gluing these simplices together via the quotient topology. This process can be carried out either concretely geometrically by choosing specific (and sufficiently geometrically independent) points within some generalized euclidean space, or, as we shall prefer to think of it,  more purely topologically by choosing standard representative simplices of the homeomorphism type of euclidean simplices and then gluing abstractly. 
 
It is worth noting separately the important point that, just like for a geometric simplicial complex, a simplex in an abstract simplicial complex is completely determined by its collection of vertices.

\subsection{Simplicial maps}

The appropriate notion of a morphism between two geometric simplicial complexes  is the simplicial map. Such maps will play an important role as we transition  from simplicial complexes to simplicial sets. 

Recall (see \cite[Section 2]{MK}) that if $K$ and $L$ are geometric simplicial complexes, then a simplicial map $f\colon K\to L$ is determined by taking the vertices $\{v_i\}$ of $K$ to vertices $\{f(v_i)\}$ of $L$ such that if $[v_{i_0},\ldots, v_{i_k}]$ is a simplex of $K$ then $f(v_{i_0}),\ldots, f(v_{i_k})$ are all vertices (not necessarily unique) of some simplex in $L$. Given such a function $K^0\to L^0$, the rest of $f\colon K\to L$ is determined by linear interpolation on each simplex (if $x\in K$ can be represented by $x=\sum_{j=1}^n t_jv_{i_j}$ in barycentric coordinates of the simplex spanned by the $v_{i_j}$, then $f(x)=\sum_{j=1}^n t_jf(v_{i_j})$).  The resulting function $f\colon K\to L$ is continuous (see \cite{MK}). 

\begin{example}
A simple, yet interesting and important example, is the inclusion of an $n$-simplex into a simplicial complex (Figure \ref{F: fig5-0}). If $X$ is a simplicial complex and ${v_{i_0},\ldots, v_{i_n}}$ is a collection of vertices of $X$ that spans an $n$-simplex of $X$, then $K=[v_{i_0},\ldots, v_{i_n}]$ is itself a simplicial complex. We then have a simplicial map $K\to X$ that takes each $v_{i_j}$ to the corresponding vertex in $X$ and hence takes $K$ identically to itself inside $X$.

\begin{figure}[!htp]
\begin{center}
\scalebox{.6}{\includegraphics{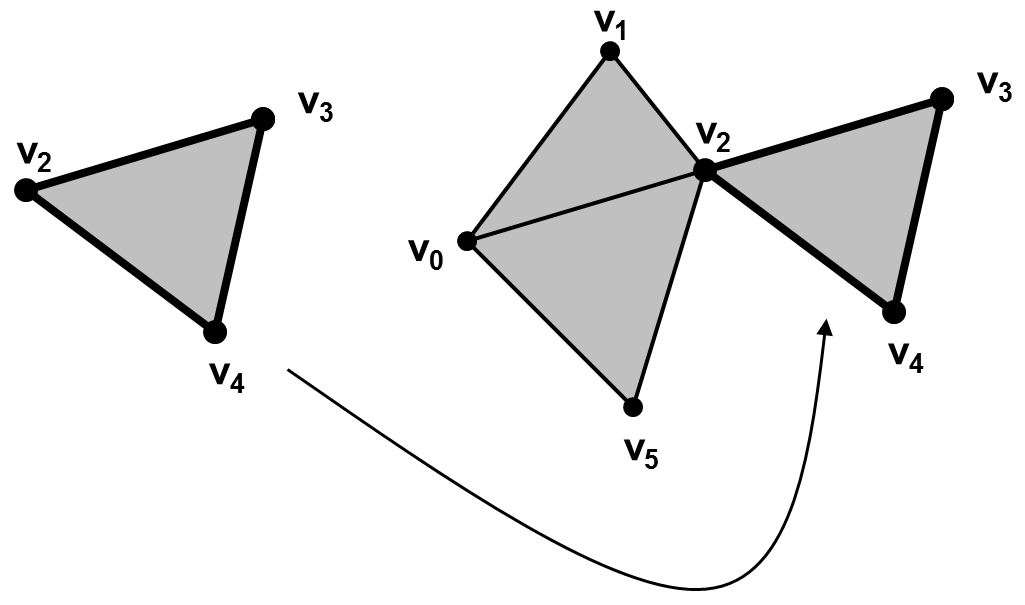}}
\end{center}
\caption{Including the simplex $[v_2,v_3,v_4]$ into a larger simplicial complex}\label{F: fig5-0}
\end{figure}
\end{example}

\begin{figure}[!htp]
\begin{center}
\scalebox{.6}{\includegraphics{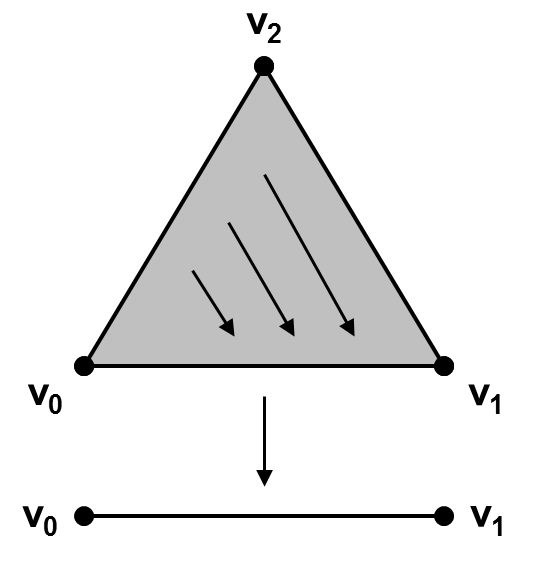}}
\end{center}
\caption{A collapse of a $2$-simplex to a $1$-simplex}\label{F: fig4}
\end{figure}

\begin{example}\label{E: collapse}
Some other  very interesting examples of simplicial maps, which will be critical for our development of  simplicial sets, are the simplicial maps that collapse simplices. For example, let $[v_0,v_1,v_2]$ be a $2$-simplex, one of whose $1$-faces is $[v_0,v_1]$.
Consider the simplicial map $f\colon [v_0,v_1,v_2]\to [v_0,v_1]$ determined by $f(v_0)=v_0$, $f(v_1)=v_1$, $f(v_2)=v_1$ that collapses the $2$-simplex down to the $1$-simplex (see Figure \ref{F: fig4}). The great benefit of the theory of simplicial sets is a way to generalize these kinds of maps in order to preserve information so that we can still see the image of the $2$-simplex hiding in the $1$-simplex as a \emph{degenerate} simplex (see Section \ref{S: simplicial sets}).
\end{example}

Of course simplicial maps of geometric simplicial complexes determine simplicial maps of abstract simplicial complexes by simply recording where each vertex of the domain goes. Conversely, observe that a simplicial map is described entirely in terms of abstract simplicial complex information; it is determined completely by specifying an image vertex for each vertex in the domain complex. Furthermore, once we have simplicial maps, we have a notion of simplicial homeomorphism, and this allows us once and for all to identify, up to simplicial homeomorphism, an abstract simplicial complex with all the geometric simplicial complexes that possess the same combinatorial data, all of which will be simplicially homeomorphic to each other. This will justify our use below of the phrase ``simplicial complex'', from which we may drop the word ``geometric'' or ``abstract''.

\subsection{Ordered simplicial complexes and face maps}

A slightly more specific way to do all this is to let the set of vertices $X^0$ of a simplicial complex $X$ be totally ordered, in which case we obtain an \emph{ordered simplicial complex}. When we do this, the symbol $[v_{i_0},\ldots, v_{i_k}]$ may stand for a simplex if and only if $v_{i_j}<v_{i_{l}}$ whenever $j<l$. This poses no undue complications as each collection $\{v_{i_0},\ldots, v_{i_k}\}$ of cardinality $k$ still corresponds to at most one simplex. We're just being picky and removing some redundancy in how many ways we can label a given simplex of a simplicial complex.

\begin{example}\label{E: n-simplex}
The prototypical example of an ordered simplicial complex is the (ordered) $n$-simplex itself\footnotemark. The \emph{ordered $n$-simplex} is simply an $n$-simplex with ordered vertices. It is an ordered simplicial complex when considered together with its faces as in Example \ref{E: simplex/complex}. We denote the ordered $n$-simplex $|\Delta^n|$; it will become clear later why we want to employ the notation $|\Delta^n|$ instead of just $\Delta^n$.  
The $n$-simplex is so fundamental that one often labels the vertices simply with the numbers $0,1,\ldots, n$, so that $|\Delta^n|=[0,\ldots,n]$ (see Figure \ref{F: fig2}). Each $k$-face of $|\Delta^n|$ then has the form $[i_0,\ldots, i_k]$, where $0\leq i_0<i_1<\ldots <i_k\leq n$. 

\footnotetext{Notice that we have already begun employing the abstraction promised at the end of the last section by referring to \emph{the} $n$-simplex. Of course, to be technical, \emph{the} $n$-simplex refers to the (abstract or geometric) simplicial homeomorphism class, as there are many different ways to realize the $n$-simplex in euclidean space as a specific geometric $n$-simplex (though of course, up to relabeling, there is only one way to describe it as an abstract simplicial complex - which is sort of the point of introducing abstract simplicial complexes in the first place). }

\begin{figure}[!htp]
\begin{center}
\scalebox{.6}{\includegraphics{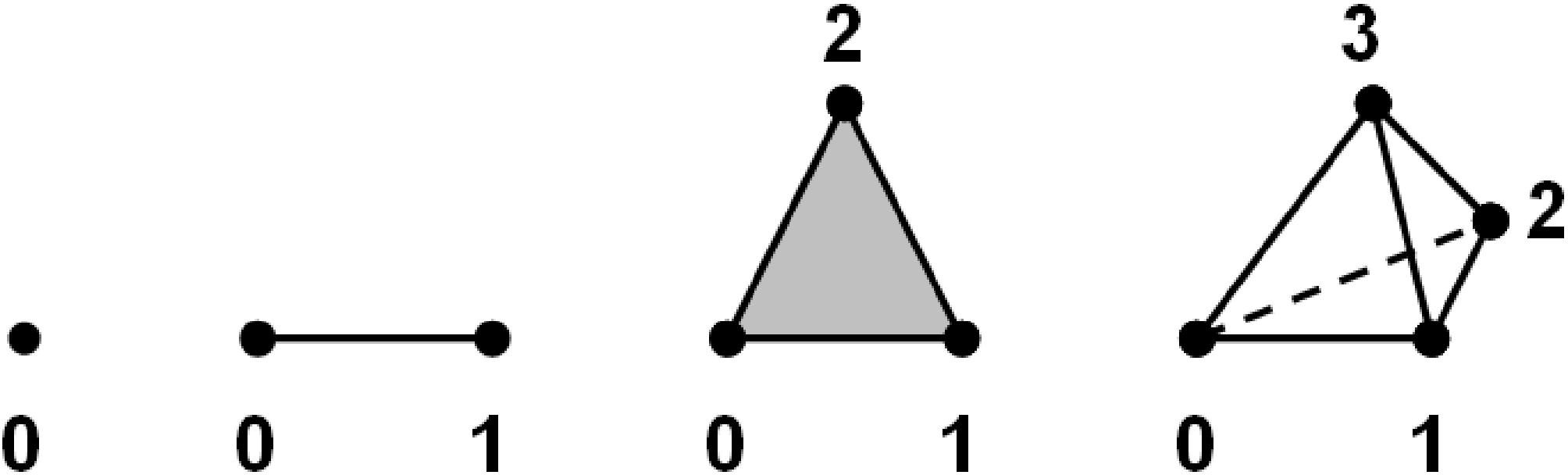}}
\end{center}
\caption{The standard ordered $0$-, $1$-, $2$-, and $3$-simplices}\label{F: fig2}
\end{figure}
\end{example}

\begin{figure}[!htp]
\begin{center}
\scalebox{.6}{\includegraphics{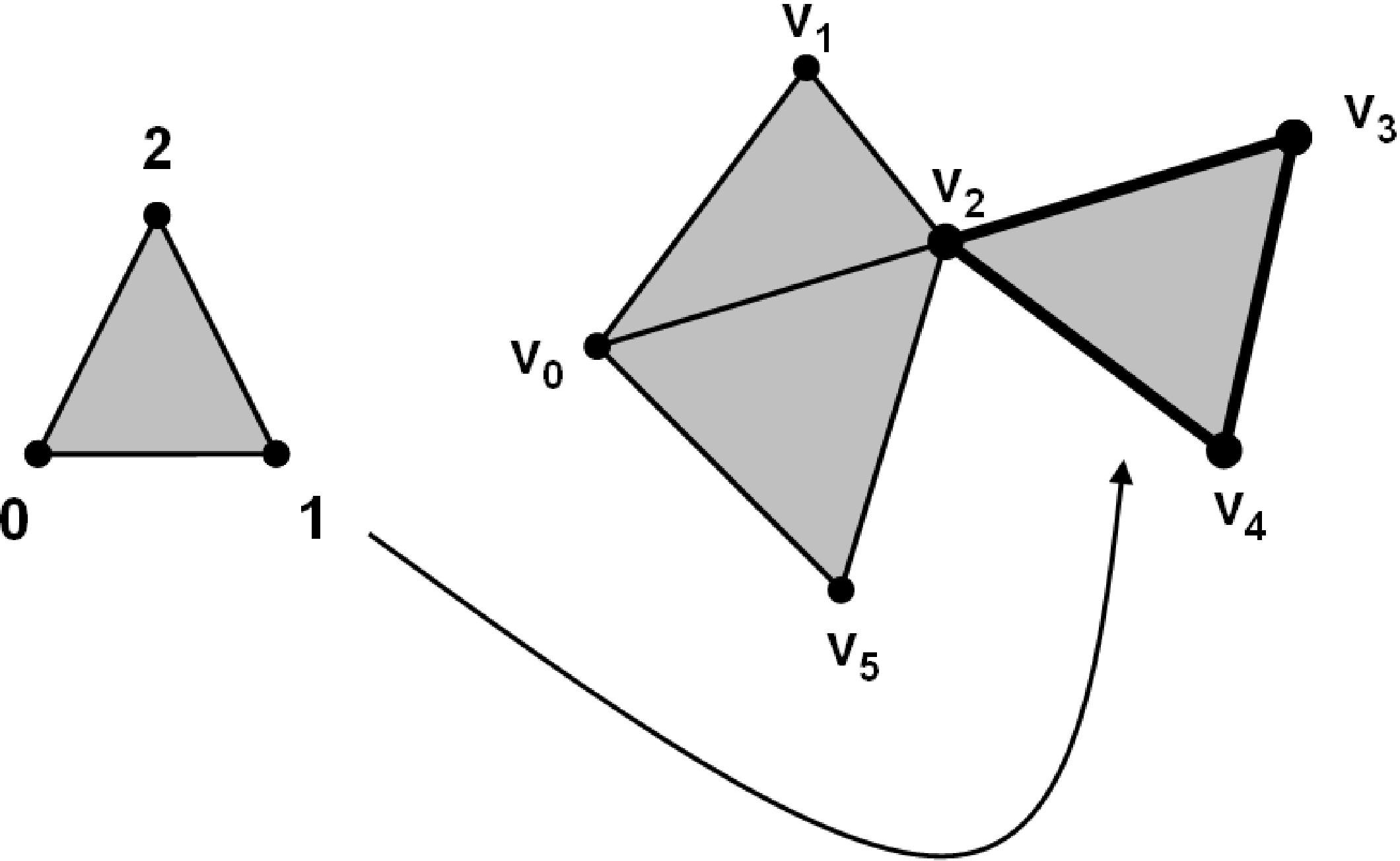}}
\end{center}
\caption{$[v_2,v_3,v_4]$ as the image of $|\Delta^2|$}\label{F: fig5}
\end{figure}

The notation $[0,\ldots,n]$ for the standard ordered $n$-simpex should be suggestive when compared with the simplices $[v_{i_0},\ldots, v_{i_n}]$ appearing within more general ordered simplicial complexes, and it is worth pointing out at this early stage that one can think of any such simplex in a complex $X$ as the image of $|\Delta^n|$ under a  simplicial map (order-preserving) taking $0$ to $v_{i_0}$, and so on. Since $X$ is an ordered simplicial complex, then there is precisely one way to do this for each $n$-simplex of $X$.
Thus another point of view on ordered simplicial complexes is that they are made up out of images of the  standard ordered simplices (Figure \ref{F: fig5}). This will turn out to be a very useful point of view as we progress.

\begin{figure}[!htp]
\begin{center}
\scalebox{.6}{\includegraphics{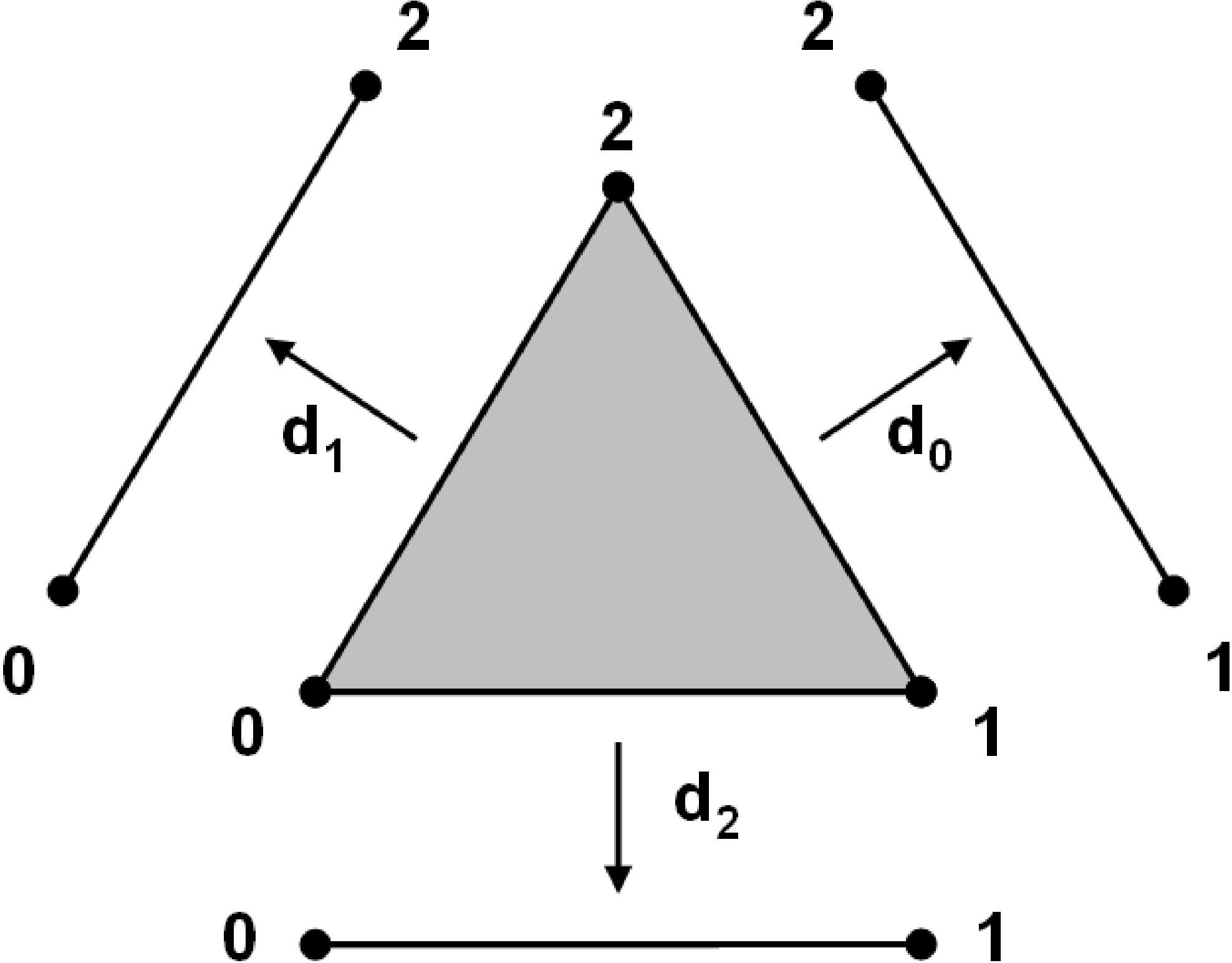}}
\end{center}
\caption{The face maps of $|\Delta^2|$. Note well: the arrows denote \emph{assignments}, not continuous maps of spaces. }\label{F: fig3}
\end{figure}

\paragraph{Face maps.} Another aspect of ordered simplicial complexes familiar to the student of basic algebraic topology is that, given an $n$-simplex, we would like a handy way of referring to its $(n-1)$-dimensional faces (its \emph{$(n-1)$-faces}). This is handled by the \emph{face maps}.  On the standard $n$-simplex, we have $n+1$ face maps $d_0,\ldots, d_n$, defined so that $d_j[0,\ldots, n]=[0,\ldots, \hat \jmath,\ldots, n]$, where, as usual, the $\hat{}$ denotes a term that is being omitted. Thus applying $d_j$ to $[0,\ldots, n]$ yields  the $(n-1)$-face missing the vertex $j$ (see Figure \ref{F: fig3}). It is important to note that each $d_j$ simply assigns to the $n$-simplex one of its faces;  there is no underlying point-set topological or simplicial map meant.

Within more general ordered simplicial complexes, we make the obvious extension: if $[v_{i_0},\ldots, v_{i_n}]\in X^n$ is a simplex of the complex $X$, then $d_j[v_{i_0},\ldots, v_{i_n}]=[v_{i_0},\ldots,\hat v_{i_j},\dots, v_{i_n}]$. Assembled all together, we get, for each fixed $n$, a collection of functions $d_0,\ldots, d_n\colon X^n\to X^{n-1}$. Note that here is where the ordering of the vertices of the simplices becomes important.

If one wanted to be a serious stickler, we might be careful to label the face maps from $X^n$ to $X^{n-1}$ as $d_0^n, \ldots, d_n^n$, but this is rarely done in practice, for which we should probably be grateful. Thus $d_j$ is used to represent the face map leaving out the $j$th vertex in any dimension where this makes sense (i.e. dimensions $\geq j$).

Furthermore, one readily sees by playing with $|\Delta^n|$ that there are certain relations satisfied by the face maps. In particular, if $i<j$, then 
\begin{equation}\label{E: d com}
d_id_j=d_{j-1}d_i.
\end{equation}
 Indeed, $d_id_j[0,\ldots, n]=[0,\ldots, \hat \imath,\ldots, \hat \jmath,\ldots, n]=d_{j-1}d_i[0,\ldots,n]$ (notice the reason that we have $d_{j-1}$ in the last expression is that removing the $i$ first shifts the $j$ into the $j-1$ slot). 

Clearly, the relation $d_id_j=d_{j-1}d_i$ must hold for any simplex in a complex $X$ (which is made up of copies of $|\Delta^n|$). This relation  will become one of the axioms in the definition of a simplicial set when we get there.

Another observation that will come up later is that there are more general face maps. We could, for example, assign to $[0,1,2,3,4,5,6]$ the face $[1,3,4]$, and we could define such general face maps systematically. However, any such face can be obtained as a composition of face maps that lower dimension by $1$. For example, we can decompose the map just described as $d_0d_2d_5d_6$. It may entertain the reader to use the ``face map relations'' and some basic reasoning to show that any generalized face map can be obtained as a composition $d_{i_1}\cdots d_{i_m}$ \emph{uniquely} if we require that $i_j<i_{j+1}$ for all $j$.

\subsection{Delta sets and Delta maps}

Delta sets (sometimes called $\Delta$-sets) constitute an intermediary between simplicial complexes and simplicial sets. These allow a degree of abstraction without yet introducing the degeneracy maps we have begun hinting at.

\begin{definition}
A \emph{Delta set}\footnote{It seems to be at least fairly usual to capitalize the word ``Delta'' in this context, probably because it is essentially a stand-in for the Greek capital letter $\Delta$. However, for reasons that will become clear, it is probably best to avoid the notation ``$\Delta$-set'' and to use instead the English stand-in.} consists of a sequence of  sets $X_0, X_1, \ldots$ and, for each $n\geq 0$, maps $d_i\colon X_{n+1}\to X_n$ for each $i$, $0\leq i\leq n+1$, such that $d_id_j=d_{j-1}d_i$ whenever $i<j$.
\end{definition}

Of course this is just an abstraction, and generalization, of the definition of an ordered simplicial complex, in which the $X_n$ are the sets of $n$-simplices and the $d_i$ are the face maps. However, there are Delta sets that are not simplicial complexes:

\begin{example}\label{E: cone}
Consider the cone $C$ obtained by starting with the standard ordered $2$-simplex $|\Delta^2|=[0,1,2]$ and gluing the edge $[0,2]$ to the edge $[1,2]$ (see Figure \ref{F: fig6}). This space is no longer a simplicial complex (at least not with the ``triangulation'' given), since in a simplicial complex, the faces of a given simplex must be unique. This is no longer the case here as, for example, the ``edge [0,1]'' now has both endpoint vertices equal to each other.

\begin{figure}[!htp]
\begin{center}
\scalebox{.6}{\includegraphics{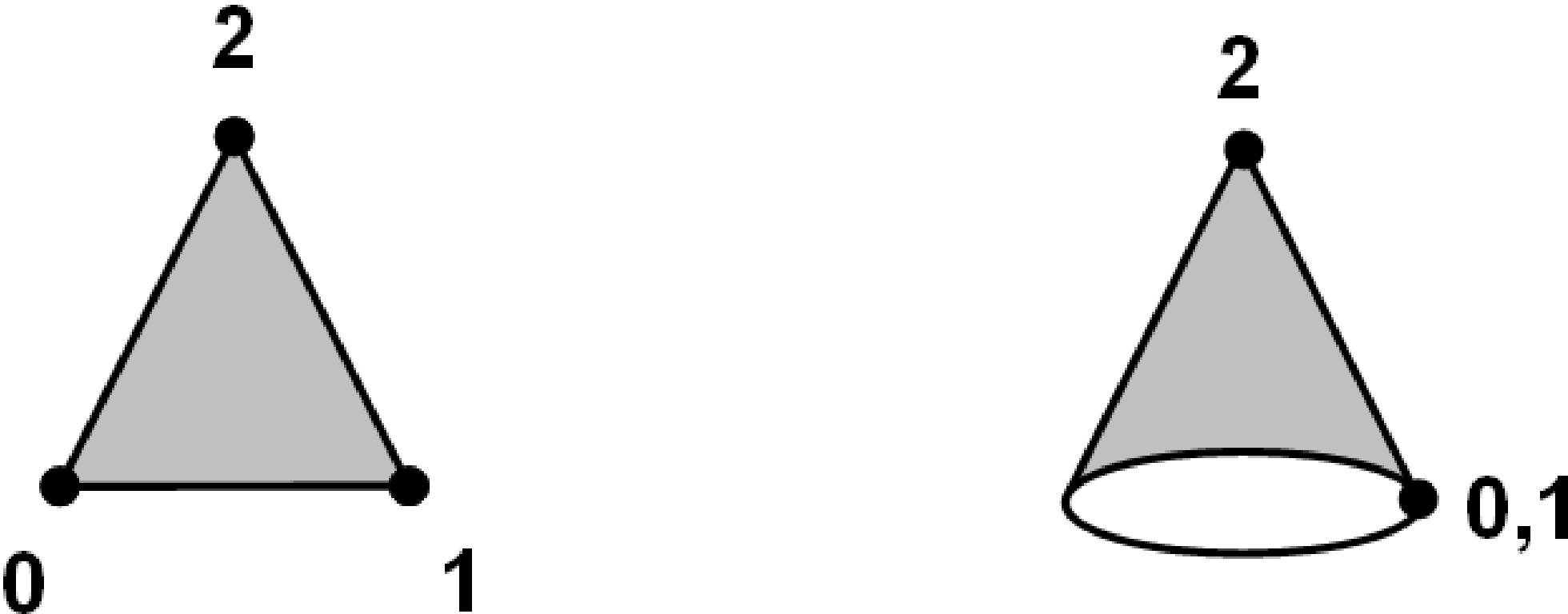}}
\end{center}
\caption{Gluing $|\Delta^2|$ into a cone}\label{F: fig6}
\end{figure}

However, this is a Delta set. Without (I hope!) too much risk of confusion, we use the notation for the simplices in the triangle to refer also to their images in the cone. So, for example $[0]$ and $[1]$ now both stand for the same vertex in the cone and $[0,1]$ stands for the circular base edge. Then $C_0=\{[0],[2]\}$, $C_1=\{[0,1], [0,2]\}$, $C_2=[0,1,2]$, and $C_n=\emptyset$ for all $n>2$. The face maps are the obvious ones, also induced from the triangle, so that, e.g. $d_2[0,1,2]=[0,1]$ and $d_0[0,1]=d_1[0,1]=[0]=[1]$. It is not hard to see that the face map relation \eqref{E: d com} is satisfied - it comes right from the fact that it holds for the standard $2$-simplex.
\end{example}  

\begin{example}\label{E: nonuniqueness}
One feature of Delta sets we need to be careful about is that, unlike for simplicial complexes, a collection of vertices does not necessarily specify a unique simplex. For example, consider the Delta set with $X_0=\{v_0,v_1\}$, $X_1=\{e_0,e_1\}$, $d_0(e_0)=d_0(e_1)=v_0$, and $d_1(e_0)=d_1(e_1)=v_1$. Both $1$-simplices have the same endpoints. See Figure \ref{F: fig15}.
\begin{figure}[!htp]
\begin{center}
\scalebox{.6}{\includegraphics{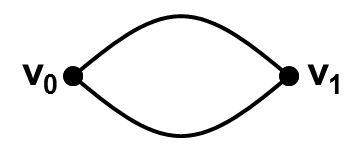}}
\end{center}
\caption{A Delta set containing two edges with the same vertices}\label{F: fig15}
\end{figure}
\end{example}

Thus Delta sets afford some greater flexibility beyond ordered simplicial complexes. One may continue to think of the sets $X_n$ as collections of simplices and interpret from the face maps how these are meant to be glued together (Exercise: Give each ``simplex'' of the cone $X$ of the preceding example an abstract label, write out the full set of face maps in these labels, then reverse engineer how to construct the cone from this information. One sees that everything is forced. For example, there is one $2$-simplex, two of whose faces are the same, so they must be glued together!). However, it is common in the fancier literature not to think of the $X_n$ as collections of simplices at all but simply as abstract sets with abstract collections of face maps. At least this is what authors would have us believe - I tend to picture simplices in my head anyway, while keeping in mind that this is more of a cognitive aid than it is ``what's really going on.''

\paragraph{The category-theoretic definition.}
While we're walking the tightrope of abstraction, let's take it a step further. Recall that 
we discussed in Example \ref{E: n-simplex} that we can think of an ordered simplicial complex as a collection of isomorphic  images of the standard $n$-simplices (for various $n$). Of course to describe the simplicial complex fully we need to know not just about these copies of the  standard simplices but also about how their  faces are attached together. This  information is contained in the face maps, which tell us when two simplices share a face. There's an alternative definition of Delta complexes that takes more of this point of view. It might be a little scary if you're not that comfortable with category theory, but don't worry, I'll walk you through it (though I do assume you know the basic language of categories and functors). 

First, we define a category $\widehat \Delta$:

\begin{definition} The category $\widehat \Delta$ has as  objects  the finite ordered sets $[n]=\{0,1,2\ldots, n\}$.  The morphisms of $\widehat \Delta$ are the strictly order-preserving functions $[m]\to [n]$ (recall that $f$ is strictly order-preserving if $i<j$ implies $f(i)<f(j)$).
\end{definition}

\begin{figure}[!htp]
\begin{center}
\scalebox{.6}{\includegraphics{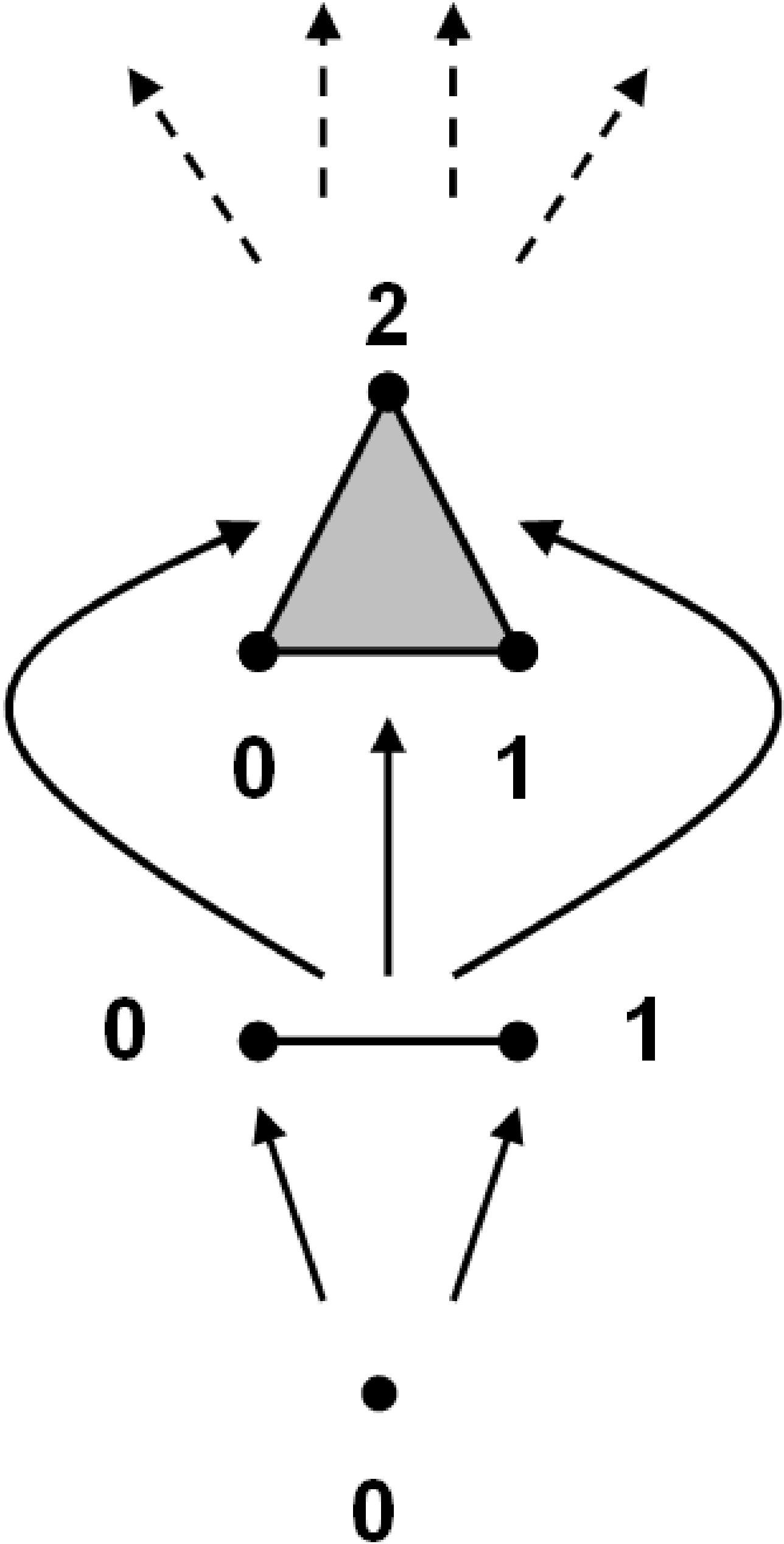}}
\end{center}
\caption{A partial illustration of the category $\widehat \Delta$}\label{F: fig7}
\end{figure}

The objects of $\widehat \Delta$ should be thought of as our prototype \emph{ordered} $n$-simplices. The morphisms are only defined when $m\leq n$, and you can think of these morphisms as taking an $m$-simplex and embedding it as a face of an $n$-simplex (see Figure \ref{F: fig7}). Note that, since order matters, there are exactly as many ways to do this as there are strictly order-preserving maps $[m]\to [n]$.

Next, we think about the opposite category $\widehat \Delta^{op}$. Recall that this means that we keep the same objects $[n]$ of $\widehat \Delta$, but for every morphism $[m]\to [n]$ in $\widehat \Delta$, we instead have a map $[n]\to [m]$ in $\widehat \Delta^{op}$. What should this mean? Well if a given morphism $[m]\to[n]$ was the inclusion of a face, then the new opposite map $[n]\to[m]$ should be thought of as taking  the $n$-simplex $[n]$ and prescribing a given face. This is just a generalization of what we have seen already: if we consider in $\widehat \Delta$ the morphism $D_i\colon[n]\to [n+1]$ defined by the strictly order-preserving map $\{0,\ldots, n\}\to \{0,\ldots, \hat \imath,\ldots, n+1\}$, then in $\widehat\Delta^{op}$ this corresponds precisely to the simplex face map $d_i$. Even better, it is easy to check once again that, with this definition, $d_id_j=d_{j-1}d_i$ when $i<j$, simply as an evident property of strictly order-preserving maps. This is really how we argued for this axiom in the first place!

So, in summary, the category $\widehat \Delta^{op}$ is just the collection of elementary $n$-simplices together with the face maps (satisfying the face map axiom) and the iterations of face maps. But this should be precisely the prototype for all Delta sets:

\begin{definition}[Alternative definition for Delta sets]
A Delta set is a covariant functor $X\colon\widehat \Delta^{op}\to \Set$, where $\Set$ is the category of sets and functions. Equivalently, a Delta set is a contravariant functor $\widehat \Delta \to \Set$. 
\end{definition}

Let's see why this makes sense. A functor takes objects to objects and morphisms to morphisms, and it obeys composition rules. So, unwinding the definition, a covariant functor $\widehat \Delta^{op}\to \Set$ assigns to $[n]\in \widehat \Delta^{op}$  a set $X_n$ (which we can think of, and which we refer to, as a set of simplices) and gives us, for each strictly order-preserving $[m]\to [n]$ in $\widehat \Delta$ (or its corresponding opposite in $\widehat \Delta^{op}$) a generalized face map $X_n\to X_m$ (which we think of as assigning an $m$-face to each simplex in $X_n$). As noted previously, these generalized face maps are all compositions of our standard face maps $d_i$, so the $d_i$ (and their axioms) are the only ones we usually bother focusing on.

\begin{figure}[!htp]
\begin{center}
\scalebox{.6}{\includegraphics{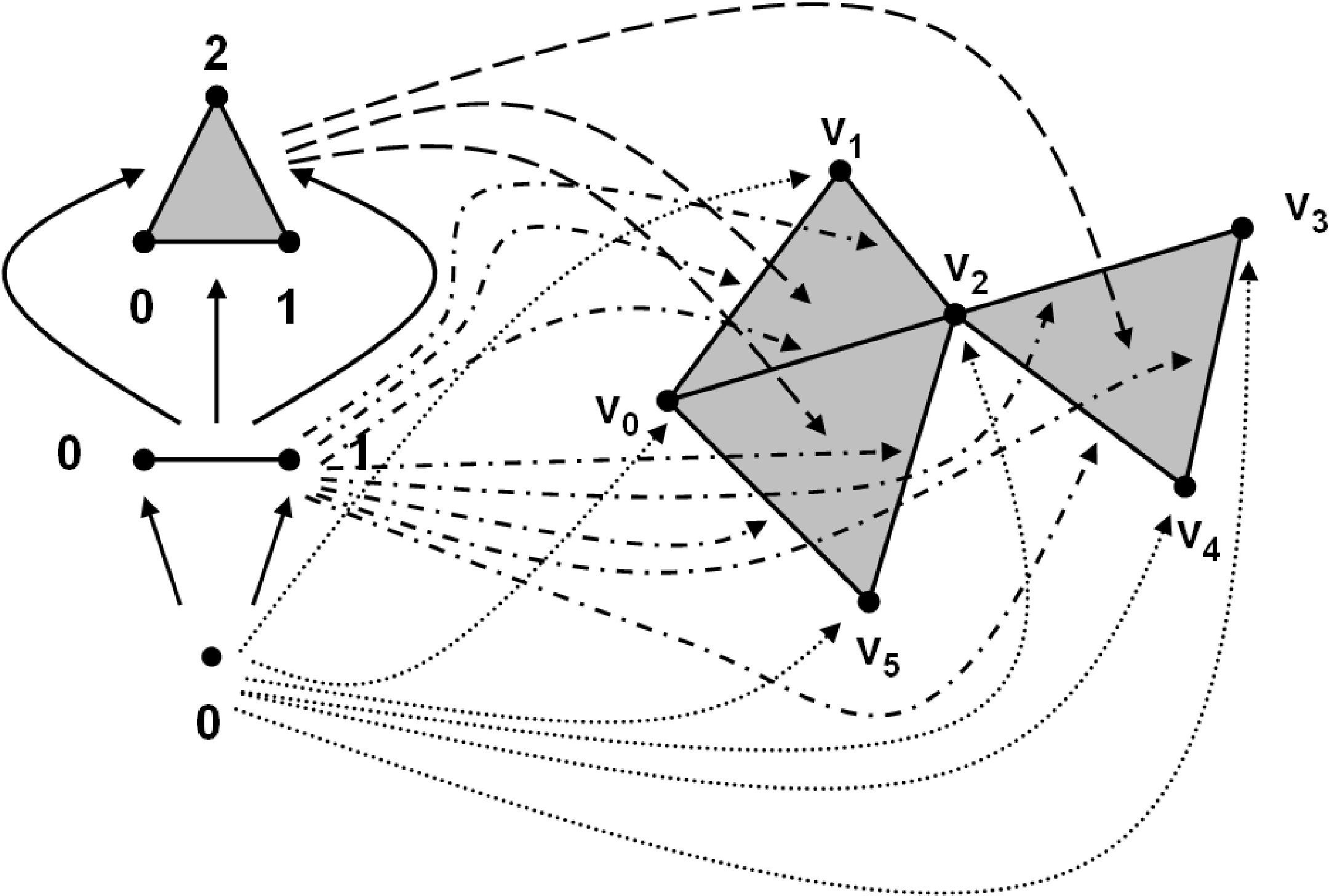}}
\end{center}
\caption{A Delta complex as the functorial image of $\widehat \Delta$}\label{F: fig8}
\end{figure}

So what just happened? The power of this definition is really in its point of view. Instead of thinking of a Delta set as being made up of a whole bunch of simplices one at a time, we can now think of the standard $n$-simplex as standing for all of the simplices in $X_n$, all at once - the functor $X$ assigns to $[n]$ the collection of all of the simplices of $X_n$ (see Figure \ref{F: fig8}). The face map $d_i$ applied to the standard simplex $[n]$ represents all of the $i$th faces of all the $n$-simplices simultaneously.

At the same time, we see how any argument in $X$ really comes from what happens back in $\widehat \Delta$. The axiom $d_id_j=d_{j-1}d_i$ in a Delta set $X$ is just a consequence of this being true in the prototype simplex $[n]$ and inherent properties of functors. We'll get a lot of mileage  out of this kind of thinking: things we'd like to prove in a Delta set $X$ can often be proved just by proving them in the prototype standard simplex and applying functoriality.

\paragraph{Delta maps.}

We won't dwell overly long on Delta maps, except to observe that they, too, point toward the need for simplicial sets (however, see \cite{RSDelta} where Delta complexes and Delta maps are treated in their own right).

Going directly to the category theoretic definition, given two Delta sets $X,Y$, thought of as contravariant functors $\widehat \Delta\to \Set$, a morphism $X\to Y$ is a natural transformation of functors from $X$ to $Y$. In other words, a morphism consists of a collection of set maps $X_n\to Y_n$ that commute with the face maps. 

\begin{figure}[!htp]
\begin{center}
\scalebox{.6}{\includegraphics{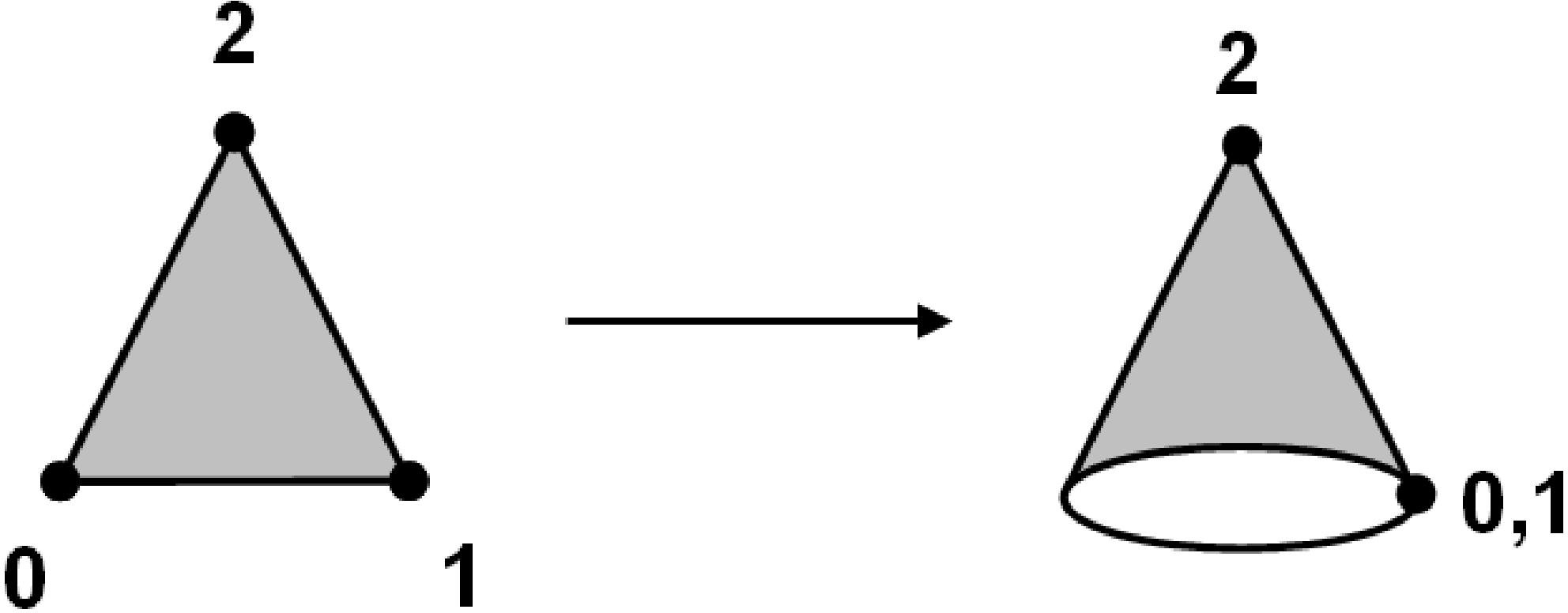}}
\end{center}
\caption{The Delta map from $|\Delta^2|$ to the cone}\label{F: fig9}
\end{figure}

\begin{example}
There is an evident Delta map from the standard $2$-simplex $[0,1,2]$ to the  cone $C$ of Example \ref{E: cone}. See Figure \ref{F: fig9}.
\end{example}

The astute reader will notice something fishy here. We would hope that simplicial maps of simplicial complexes would yield morphisms of Delta sets. However, consider the collapse $\pi\colon|\Delta^2|=[0,1,2]\to |\Delta^1|=[0,1]$ defined by $\pi(0)=0$ and $\pi(1)=\pi(2)=1$ (see Figure \ref{F: fig4}). To be a Delta set morphism, the simplex $[0,1,2]\in |\Delta^2|_2$ would have to be taken to an element of $|\Delta^1|_2$. But this set is empty! There are no $2$-simplices of $|\Delta^1|$. Something is amiss. We need simplicial sets.

\section{Simplicial sets and morphisms}\label{S: simplicial sets}

Simplicial sets generalize both simplicial complexes and Delta sets. 

When approaching the literature, the reader should be very careful about terminology. Originally (\cite{EZ50}), Delta sets were referred to as \emph{semi-simplicial complexes}, and, once the degeneracy operations we are about to discuss were discovered, the term \emph{complete semi-simplicial complex} (\emph{c.s.s. set}, for short) was introduced. Over time, with Delta sets becoming of less interest,  ``complete semi-simplicial'' was abbreviated back to ``semi-simplicial'' and eventually to ``simplicial,'' leaving us with the \emph{simplicial sets} of today. Meanwhile, some modern authors have returned to using ``semi-simplicial complexes'' to refer to what we are calling Delta sets, on the grounds that, as we will see, the category $\Delta$ (``Delta'') is the prototype for simplicial sets, not Delta sets, for which we have been using the prototype category $\widehat \Delta$. This all sounds very confusing because it is, and the reader is advised to be very careful when reading the literature.\footnote{I thank Jim McClure for explaining to me this historical progression.}

We try to be careful and use only the three terms ``simplicial complex,'' ``Delta set,'' and ``simplicial set.'' In particular, be sure to note the difference between  ``simplicial complex'' and ``simplicial set'' going forward.

\paragraph{Degenerate simplices.}

Recall from Example \ref{E: collapse} that a simplicial map can collapse a simplex. In that example, we had a simplicial map $\pi\colon|\Delta^2|\to|\Delta^1|$ defined on vertices so that $\pi(0)=0$ and $\pi(1)=\pi(2)=1$. Recall also that we have begun to think of simplicial complexes and Delta sets as collections of images of standard simplices under appropriate maps. Well, here is a map of the standard $2$-simplex $|\Delta^2|$. What image simplex does it give us in $|\Delta^1|$ under $\pi$? In the land of simplicial sets, the image $\pi(|\Delta^2|)$ is an example of a \emph{degenerate} simplex.

Roughly speaking, degenerate simplices are simplices that don't have the ``correct'' number of dimensions. A degenerate $3$-simplex might be realized geometrically as a $2$-dimensional, $1$-dimensional, or $0$-dimensional object. Geometrically, degenerate simplices are ``hidden''; thus the clearest approach to dealing with them lies in the combinatorial notation we have been developing all along. 

The key both to the idea and to the notation is in allowing vertices to repeat. The natural way to label $\pi(|\Delta^2|)=\pi([0,1,2])$ in our example is as $[0,1,1]$, reflecting where the vertices of $|\Delta^2|$ go under the map. This violates our earlier principle that simplices in complexes should be written $[v_0,\ldots,v_n]$ with the $v_i$  distinct  vertices written in order, but sometimes in mathematics we need a new, more general principle. For degenerate simplices, we'll keep the orderings but dispense with the uniqueness. Thus, officially, a \emph{degenerate} simplex is a $[v_{i_0},\ldots, v_{i_n}]$ for which  the $v_{i_j}$ are \emph{not} all distinct, though we do still require  $i_k\leq i_\ell$ if $k<\ell$.

\begin{figure}[!htp]
\begin{center}
\scalebox{.6}{\includegraphics{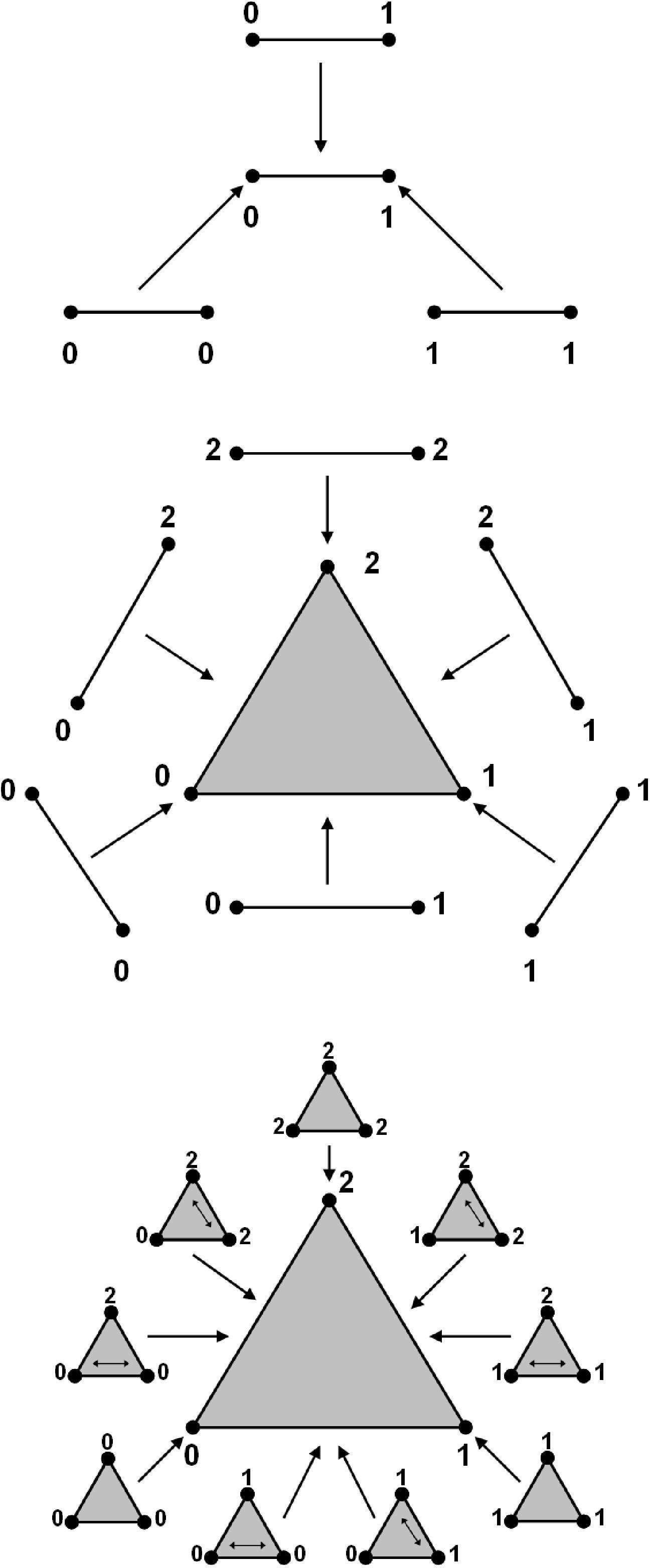}}
\end{center}
\caption{The first picture represents all of the $1$-simplices in $|\Delta^1|$, including the degenerate ones that are taken to individual vertices. The second picture represents all the $1$-simplices in $|\Delta^2|$, and the last picture represents all of the degenerate $2$-simplices in $|\Delta^2|$.}\label{F: fig10}
\end{figure}

\begin{example}
How many  $1$-simplices, including degenerate ones, are lurking within the elementary $2$-simplex $[0,1,2]$? A $1$-simplex is still written $[a,b]$, with $a\leq b$, but now repetition is allowed. The answer is six: $[0,1]$, $[0,2]$, $[1,2]$, $[0,0]$, $[1,1]$, and $[2,2]$. See the middle picture in Figure \ref{F: fig10}.

Similarly, within $|\Delta^2|=[0,1,2]$ there are now three kinds of $2$-simplices. We have the nondegenerate $[0,1,2]$, the $2$-simplices that degenerate to $1$-dimension such as $[0,1,1]$ and $[0,0,2]$, and we have the $2$-simplices that degenerate to $0$-dimensions such as $[0,0,0]$ and $[2,2,2]$. 
\end{example}

Working with degenerate simplices makes drawing diagrams much more difficult. We take  a crack at it in Figure \ref{F: fig10}.

As implied by the diagram, we can think of degenerate simplices as being the images of collapsing maps such as that in Example \ref{E: collapse}.

Of course any simplicial complex or Delta set can be expanded conceptually to include degenerate simplices. In the example of Figure \ref{F: fig1}, we might have the degenerate $5$-simplex $[v_2,v_2, v_2,v_3, v_3]$.

Notice also that our innocent little $n$-dimensional simplicial complexes suddenly contain degenerate simplices of arbitrarily large dimension. Even the $0$-simplex $|\Delta^0|=[0]$ becomes host to degenerate simplices such as $[0,0,0,0,0,0,0,0,0,0,0,0,0,0,0,0,0,0]$. 

The situation has degenerated indeed! To keep track of it all, we need \emph{degeneracy maps}.

\paragraph{Degeneracy maps.}

Degeneracy maps are, in some sense, the conceptual converse of face maps. Recall that the face map $d_j$ takes an $n$-simplex and give us back its $j$th $(n-1)$-face. On the other hand, the $j$th degeneracy map $s_j$ takes an $n$-simplex and gives us back the $j$th degenerate $(n+1)$-simplex living inside it. 

As usual, we illustrate with the standard $n$-simplex, which will be a model for what happens in all simplicial sets. Given the standard $n$-simplex $|\Delta^n|=[0,\ldots, n]$, there are $n+1$ degeneracy maps $s_0,\ldots, s_n$, defined by $s_j[0,\ldots, n]=[0,\ldots, j,j,\ldots, n]$. In other words, $s_j[0,\ldots,n]$ gives us the unique degenerate $n+1$ simplex in $|\Delta^n|$ with only the $j$th vertex repeated.

Again, the geometric concept is that $s_j|\Delta^n|$ can be thought of as the process of collapsing $\Delta^{n+1}$ down into $\Delta^n$ by the simplicial map $\pi_j$ defined by $\pi_j(i)=i$ for $i<j$, $\pi_j(j)=\pi_j(j+1)=j$ and $\pi_j(i)=i-1$ for $i>j+1$.

This idea extends naturally to simplicial complexes, to Delta sets, and to simplices that are already degenerate. If we have a  (possibly degenerate) $n$-simplex $[v_{i_0},\ldots, v_{i_n}]$ with $i_k\leq i_{k+1}$ for each $k$, $0\leq k< n$, then we set $s_j[v_{i_0},\ldots, v_{i_n}]=[v_{i_0},\ldots, v_{i_j},v_{i_{j}},\ldots, v_{i_n}]$, i.e. repeat $v_{i_j}$. This is a degenerate simplex in $[v_{i_0},\ldots, v_{i_n}]$.

It is not hard to see that any degenerate simplex can be obtained from an ordinary simplex by repeated application of degeneracy maps. Thus, just as any face of a simplex can be obtained by using compositions of the $d_i$,  any degenerate simplex can be obtained from compositions of the $s_i$. 

Also, as for the $d_i$, there are certain natural relations that the degeneracy maps possess. In particular, if $i\leq j$, then $s_is_j[0,\ldots,n]=[0,\ldots, i,i,\ldots, j,j,\ldots,n]=s_{j+1}s_i [0,\ldots,n]$. Note that we have $s_{j+1}$ in the last formula, not $s_j$, since the application of $s_i$ pushes $j$ one slot to the right. 

Furthermore, there are relations amongst the face and degeneracy operators. These are a little more awkward to write down since there are three possibilities: 
\begin{align*}
 d_i s_j &= s_{j-1} d_i\quad\text{ if $i < j$},\\
d_j s_j &= d_{j+1} s_j=\text{id},\\
 d_i s_j &= s_j d_{i-1}\quad  \text{ if $i > j + 1$}.
\end{align*}
These can all be seen rather directly. For example, applying either side of the first formula to $[0,\ldots,n]$ yields $[0,\ldots,\hat \imath,\ldots ,j,j,\ldots, n]$. Note also that the middle formula takes care of both $i=j$ and $i=j+1$.

\paragraph{Simplicial sets.}

We are finally ready for the definition of simplicial sets:
\begin{definition}\label{D: ss}
A \emph{simplicial set} consists of a sequence of sets $X_0, X_1, \ldots$ and, for each $n\geq 0$, functions $d_i\colon X_{n}\to X_{n-1}$ and $s_i\colon X_n\to X_{n+1}$ for each $i$ with $0\leq i\leq n$  such that 
\begin{align}
 d_id_j&=d_{j-1}d_i\quad \text{ if $i<j$}, \notag\\
 d_i s_j &= s_{j-1} d_i\quad\text{ if $i < j$},\notag \\
d_j s_j &= d_{j+1} s_j=\text{id},\label{E: axioms}\\
 d_i s_j &= s_j d_{i-1}\quad  \text{ if $i > j + 1$},\notag\\
  s_is_j&=s_{j+1}s_i\quad \text{ if $i\leq j$}.\notag
\end{align}
\end{definition}

\begin{example}\label{E: comp to set}
Our first example is the critical observation that every ordered simplicial complex  can be made into a simplicial set by adjoining all possible degenerate simplices. More precisely, suppose $X$ is an ordered simplicial complex. Then we obtain a simplicial set\footnote{The notation transition $X$ to $\bar X$ from an ordered simplicial complex to a simplicial set is not standard notation; we simply use it for expediency in this example.} $\bar X$ such that $\bar X_n$ consists of all the simplices $[v_{i_0},\ldots,v_{i_n}]$ where $v_{i_k}\leq v_{i_{k+1}}$ and the \emph{set} of vertices $\{v_{i_0},\ldots,v_{i_n}\}$  spans a simplex of $X$; note that the $v_{i_j}$ are not required to be unique. Another way to say this is that for every simplex $[v_{i_0},\ldots,v_{i_m}]$ of $X$, we have in $\bar X$ all simplices of the form $[v_{i_0},\ldots,v_{i_0}, v_{i_1},\ldots, v_{i_1},\ldots, v_{i_m}]$ for any number of repetitions of each of the vertices. The face and degeneracy maps are defined on these simplices in the evident ways. Similarly, every Delta set can be ``completed'' to a simplicial set by an analogous process, though some additional care is necessary as we know that an element of a Delta set is not necessarily determined by its vertices; we leave the precise construction as an exercise for the reader.

Conversely, each simplicial set yields a Delta set by neglect of structure (throw away the degeneracy maps). However, a simplicial set does not necessarily come from an ordered simplicial complex by the process described above as, for example, not every Delta set is an ordered simplicial complex.
\end{example}

\begin{example}
The standard $0$-simplex $X=[0]$, now thought of as a simplicial set, is the unique simplicial set with one element in each $X_n$, $n\geq 0$. The element in dimension $n$ is $\overset{n+1 \text{ times}}{\overbrace{[0,\ldots,0]}}$. 
\end{example}

\begin{example}
As a simplicial set, the standard ordered $1$-simplex $X=[0,1]$ already has $n+2$ elements in each $X_n$. For example, $X_2=\{[0,0,0], [0,0,1],[0,1,1],[1,1,1]\}$. 
\end{example}

\begin{remark}
We will   use $\Delta^n$ or $[0,\ldots, n]$ to refer to the standard ordered $n$-simplex, thought of  as a simplicial set. 
\end{remark}

\begin{example}\label{E: singular}
Now for an example familiar from algebraic topology. Given a topological space $X$, let  $\ms S(X)_n$ be the set of continuous functions from $|\Delta^n|$  to $X$. Together with face and degeneracy maps that we will describe in a moment, these constitute a simplicial set called the singular set of $X$. The singular chain complex $S_*(X)$  from algebraic topology has each $S_n(X)$ equal to the free abelian group generated by $\ms S(X)_n$. 

To define the face and degeneracy maps,  let $\sigma\colon|\Delta^n|\to X$ be a continuous map representing a singular simplex (Figure \ref{F: fig11a}). The singular simplex $d_i\sigma$ is defined as the  restriction of  $\sigma$ to the $i$th face of $|\Delta^n|$. Equivalently it is the composition of $\sigma$ and the simplicial inclusion map $[0,\ldots, n-1]\to [0,\ldots, \hat \imath,\ldots, n]$ (Figure \ref{F: fig11b}). These are precisely the same as the terms that show up in the boundary map of the singular chain  complex where $\bd=\sum_{i=0}^n(-1)^id_i$. 
\begin{figure}[!htp]
\begin{center}
\scalebox{.6}{\includegraphics{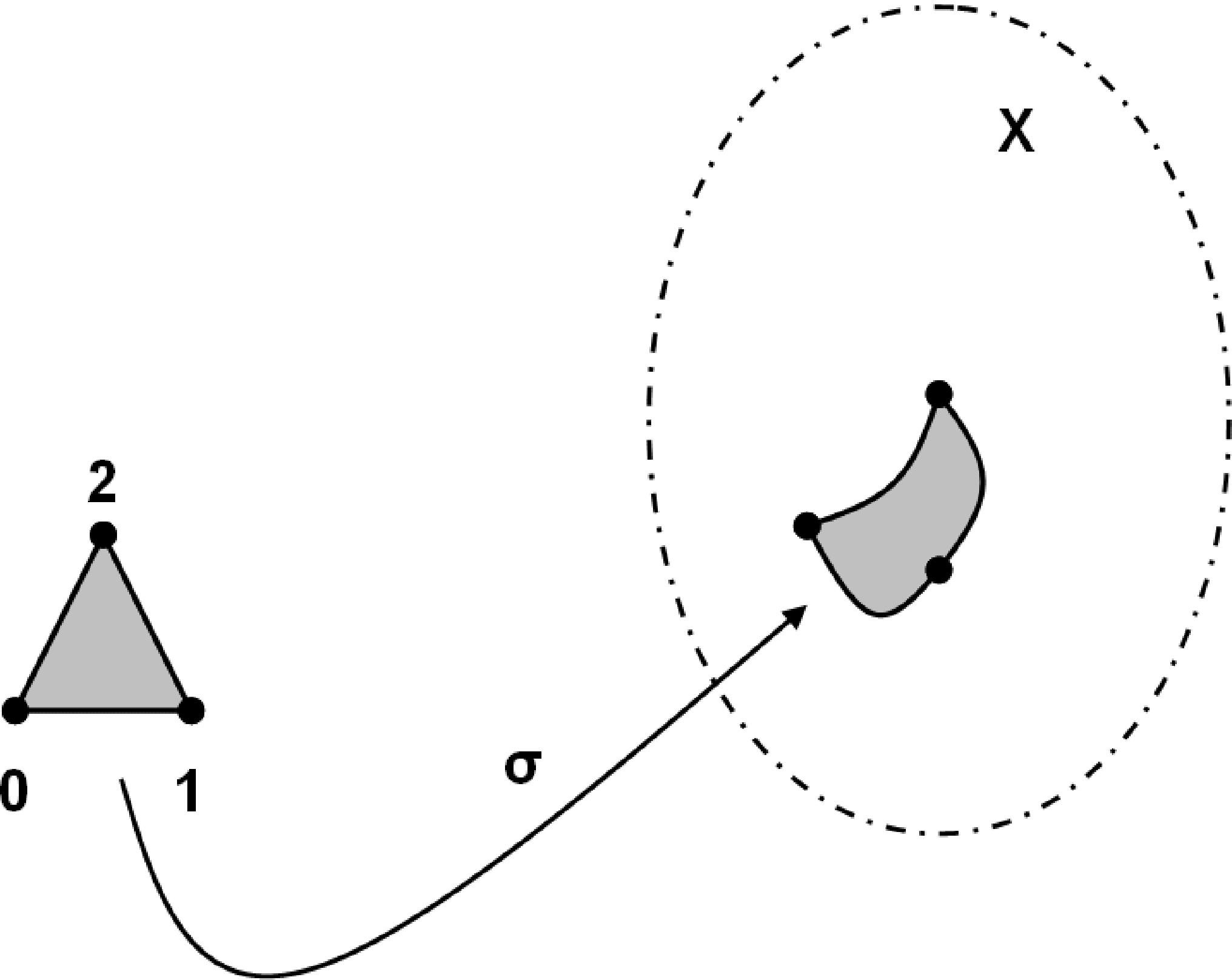}}
\end{center}
\caption{A singular simplex}\label{F: fig11a}
\end{figure}

\begin{figure}[!htp]
\begin{center}
\scalebox{.6}{\includegraphics{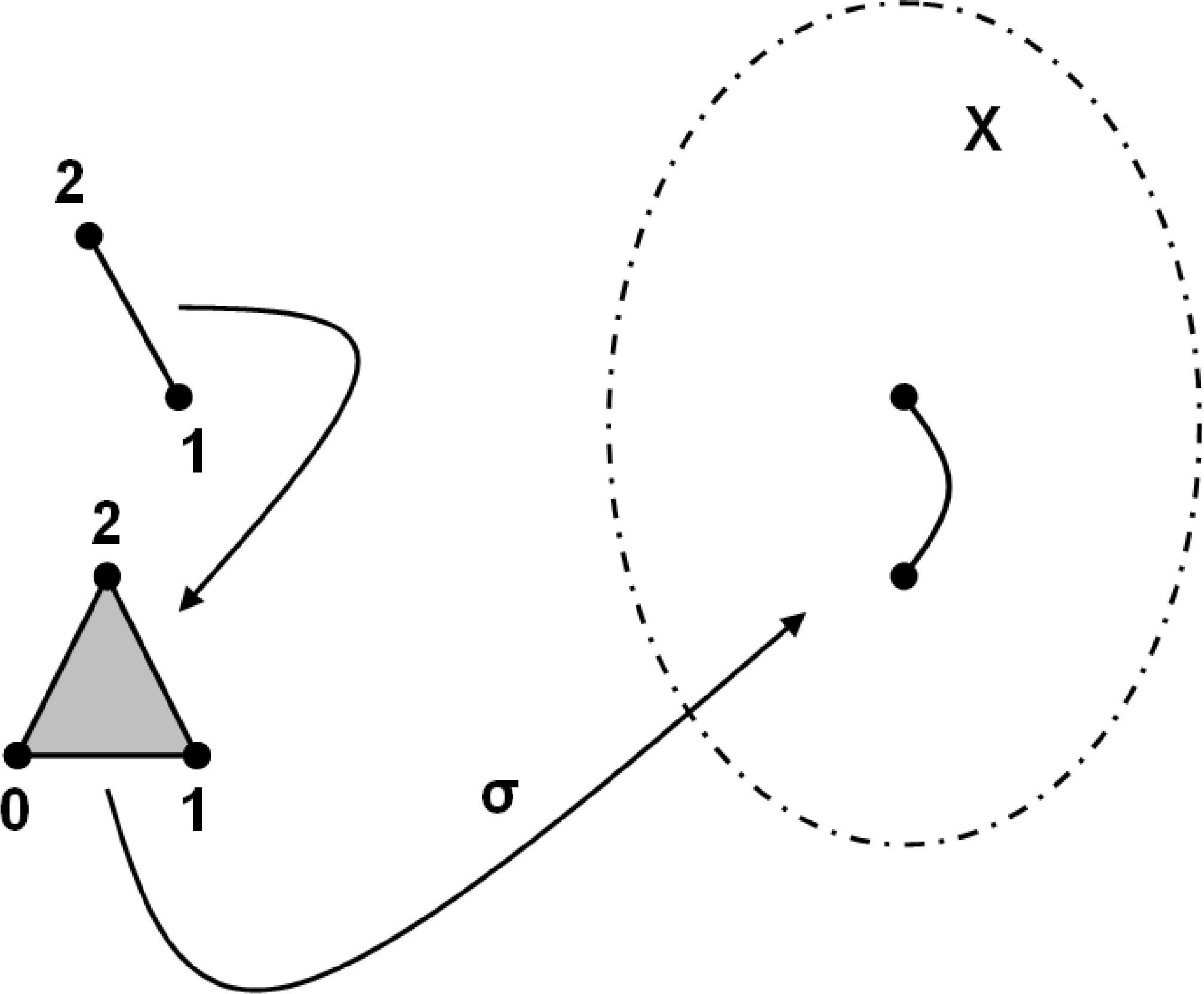}}
\end{center}
\caption{A face of a singular simplex}\label{F: fig11b}
\end{figure}

On the other hand, the degeneracy $s_i$ takes the singular simplex $\sigma$ to the composition of $\sigma\colon|\Delta^n|=[0,\ldots, n]\to X$ with the geometric collapse represented by the degeneracy $[0,\ldots, n+1]\to[0,\ldots, i,i,\ldots,n]$.  Once again, a degenerate simplex is a collapsed version of another simplex (Figure \ref{F: fig11c}).

\begin{figure}[!htp]
\begin{center}
\scalebox{.6}{\includegraphics{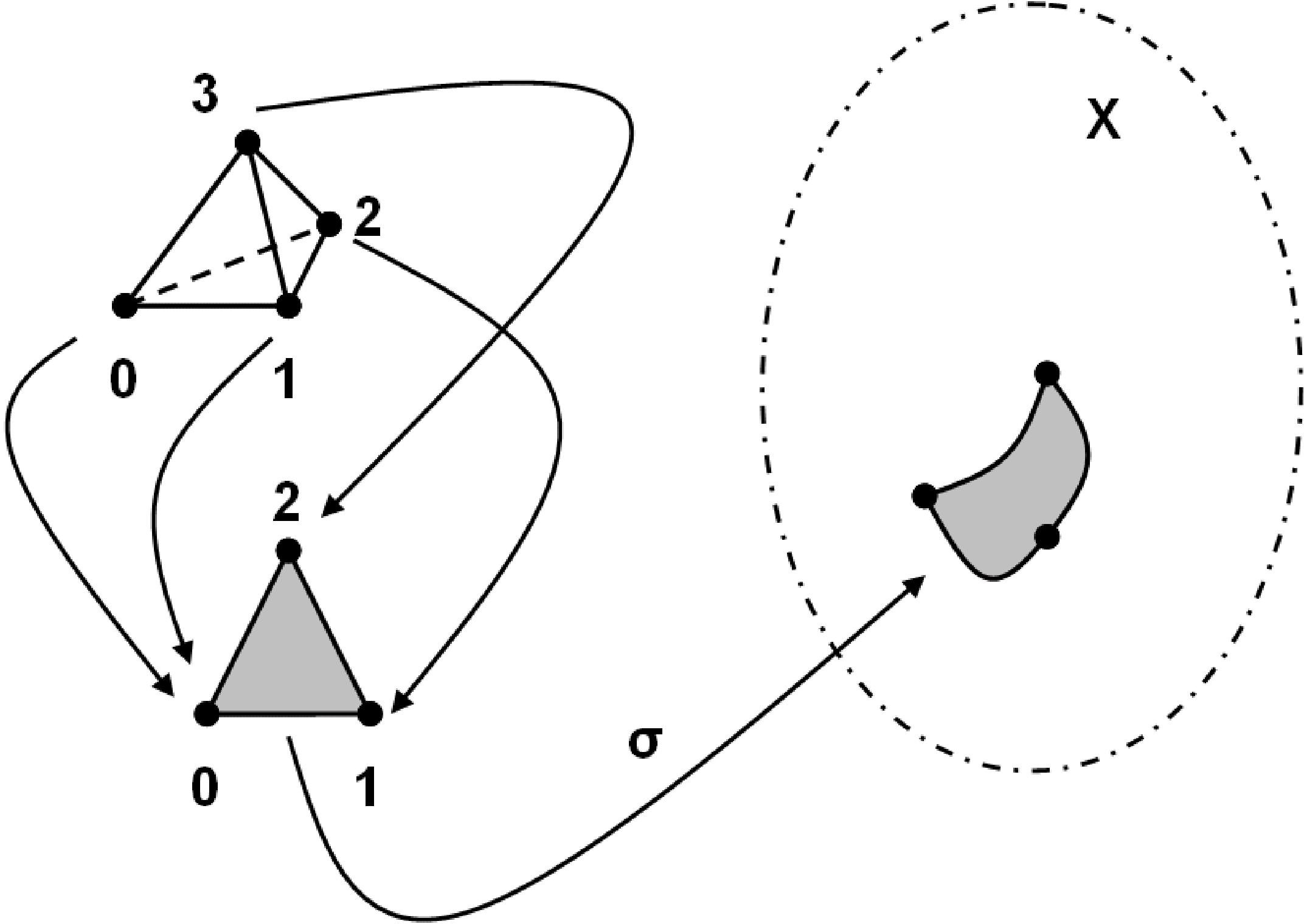}}
\end{center}
\caption{A degenerate singular simplex}\label{F: fig11c}
\end{figure}

$\ms S(X)$ turns out to be simplicial set, and we invite the reader to think through why the relations \eqref{E: axioms} hold as a consequence of their holding for  the standard ordered simplex. 
In some sense, this is our usual model, just redesigned within the context of the continuous map $\sigma$.  
\end{example}

Some more examples of simplicial sets are given below in Section \ref{S: realization}, where we can better study their geometric manifestations.

\paragraph{Nondegenerate simplices.}

\begin{definition}
A simplex $x\in X_n$ is called \emph{nondegenerate} if $x$ cannot be written as $s_iy$ for any $y\in X_{n-1}$ and any $i$.
\end{definition}

Every simplex in the sense of Section \ref{S: simplicial complexes} of a simplicial complex or Delta set  is a nondegenerate simplex
of the corresponding simplicial set.  If $Y$ is a topological space, an $n$-simplex of $\ms S(Y)$ is nondegenerate if it cannot be written as the composition $\Delta^n\overset{\pi}{\to} \Delta^k\overset{\sigma}{\to} Y$, where $\pi$ is a simplicial collapse with $k<n$ and $\sigma$ is a singular $k$-simplex.

Note that it is possible for a nondegenerate simplex to have a degenerate face (see Example \ref{E: sphere}, below, though it might be good practice to try to come up with your own example first).  It is also possible  for a degenerate simplex to have a nondegenerate face (for example, we know $d_js_jx=x$ for any $x$, degenerate or not).

\paragraph{The categorical definition.} As for Delta sets, the basic properties of simplicial sets derive from those of the standard ordered $n$-simplex. In fact, that is where the prototypes of both the face and degeneracy maps live and where we first developed the axioms relating them. Thus it is not surprising (at this point) that there is a categorical definition of simplicial sets, analogous to the one for Delta sets, in which each simplicial set is the functorial image of a category, $\Delta$, built from the standard simplices.

\begin{definition} The category $\Delta$ has as  objects  the finite ordered sets $[n]=\{0,1,2\ldots, n\}$.  The morphisms of $\Delta$ are order-preserving functions $[m]\to [n]$.
\end{definition}

Notice that the only difference between the definitions of $\widehat \Delta$ and $\Delta$ is that the morphisms in $\Delta$ only need to be order-preserving and not strictly order-preserving. Thus, equating the objects $[n]$ with the ordered simplices $\Delta^n$, the morphisms no longer need to represent only inclusions of simplices but may represent degeneracies as well. In more familiar notation, a typical morphism, say,  $f\colon[5]\to [3]$ might be described by $f[0,1,2,3,4,5]=[0,0,2,2,2,3]$, which can be thought of as a simplicial complex map taking the $5$-simplex degenerately to the $2$-face of the $3$-simplex spanned by $0$, $2$, and $3$. 

As in $\widehat \Delta$, the morphisms in $\Delta$ are generated by certain maps between neighboring cardinalities $D_i\colon[n]\to [n+1]$ and $S_i\colon[n+1]\to [n]$, $0\leq i\leq n$. The $D_i$ are just as for $\widehat \Delta$: $D_i[0,\ldots, n]= [0,\ldots, \hat \imath,\ldots, n+1]$. The new maps, which couldn't exist in $\widehat \Delta$, are defined by $S_i[0,\ldots, n+1]= [0,\ldots, i, i,\ldots, n]$. It is an easy exercise to verify that all morphisms in $\Delta$ are compositions of the $D_i$ and $S_i$ and that these satisfy axioms analogous to those in the definition of simplicial set. Later on, we will also use $D_i$ and $S_i$ to stand for the geometric maps they induce on the  standard geometric simplices. 

\begin{figure}[!htp]
\begin{center}
\scalebox{.6}{\includegraphics{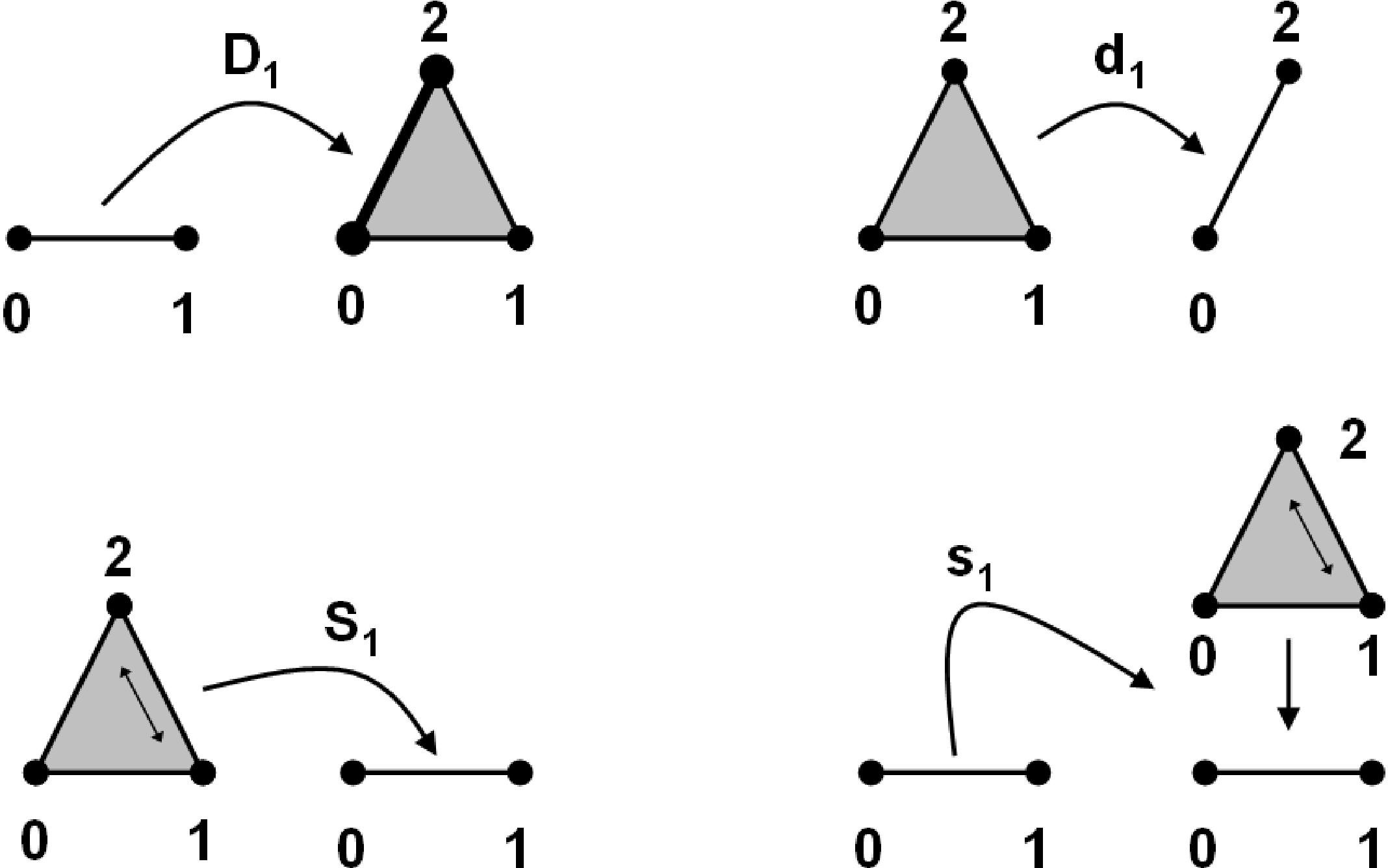}}
\end{center}
\caption{How to visualize $D_i$, $d_i$, $S_i$, and $s_i$. Our difficulty with drawing degeneracies extends here so that we represent the image of $s_i$ pictorially by the picture for $S_i$. In other words, the image of $s_1$ in the bottom right is the degenerate $2$-simplex arising from  the collapse map $S_1$.}\label{F: fig12}
\end{figure}

To get to our categorical definition of simplicial sets, we must, as for Delta sets, consider $\Delta^{op}$. The maps $D_i$ become their opposites, denoted $d_i$, and these correspond to the face maps as before: the opposite of the inclusion $D_i\colon[n]\to [n+1]$ of the $i$th face is the $i$th face map, $d_i$, which assigns to the $n$-simplex its $i$th face. The opposites of the $S_i$ become the degeneracies; the opposite of the collapse  $S_i\colon[n+1]\to [n]$  that pinches together the  $i$-th and $i+1$-th vertices of an $n+1$ simplex is the $i$th degeneracy map, $s_i$, which assigns to the $n$-simplex $\Delta^n$ the degenerate $n+1$-simplex within $\Delta^{n}$ that repeats the  $i$th vertex. See Figure \ref{F: fig12}.

Of course, one can check that the $d_i$ and $s_i$ satisfy the axioms in the definition of simplicial set given above.

\begin{definition}[Categorical definition of simplicial set]
A \emph{simplicial set} is a contravariant functor $X\colon \Delta\to \Set$ (equivalently, a covariant functor $X\colon  \Delta^{op}\to \Set$). 
\end{definition}

The reader should compare this with the categorical definition of Delta sets and reassure himself/herself that this definition is equivalent to Definition \ref{D: ss}. 
As for Delta sets, the power in this definition is that we can think of the standard ordered $n$-simplex as standing for all of the simplices in $X_n$, all at once - the functor $X$ assigns to $[n]$ all of the $n$-simplices in $X_n$ - and the standard face and degeneracy maps $d_i$ and $s_i$  pick out all of the faces and degeneracies  of $X_n$ by functoriality.

\begin{example}
Let's re-examine the singular set  $\ms S(Y)$ of the topological space $Y$ from this point of view. The singular set $\ms S(Y)$ is a functor $\Delta\to \Set$ that assigns to $[n]$ the set $\Hom_{\Top}(|\Delta^n|,Y)$, the set of all continuous maps from $|\Delta^n|$ to $Y$. It assigns to the face and degeneracy maps of $\Delta$  the face and degeneracy maps of Example \ref{E: singular}, i.e. we have the following correspondences:
\begin{diagram}
[n] & & \Hom_{\Top}(|\Delta^n|,Y)&&&  [n] & &  \Hom_{\Top}(|\Delta^{n}|,Y)\\
\dTo^{d_i} &\Rightarrow& \dTo_{d_i} &&&\dTo^{s_i} &\Rightarrow& \dTo_{s_i}\\
[n-1] && \Hom_{\Top}(|\Delta^{n-1}|,Y)&&&[n+1] & & \Hom_{\Top}(|\Delta^{n+1}|,Y).
\end{diagram}
The reader  should check that the definitions for the face and degeneracy maps of the singular set defined above are consistent with the claimed functoriality. (Notice that the maps on the right sides of these diagrams should more appropriately be labeled $\ms S(Y)(d_i)$ and $\ms S(Y)(s_i)$, but we stick with common practice and use $d_i$ and $s_i$ for face and degeneracy maps wherever we find them.) 
\end{example}

\paragraph{Simplicial morphisms.}
Simplicial sets themselves constitute a category $\mbf S$. The morphisms in this category are the \emph{simplicial morphisms}:

\begin{definition} If $X$ and $Y$ are simplicial sets (and thus functors $X,Y:\Delta\to \Set$), then a \emph{simplicial morphism} $f\colon X\to Y$ is a natural transformation of these functors. 
\end{definition}

Unwinding this to more concrete language, $f$ consists of set maps $f_n\colon X_n\to Y_n$ that commute with face operators and with degeneracy operators.

\begin{example}\label{E: simp mor v}
At last we have a context in which to explore properly the collapse map $\pi\colon|\Delta^2|\to |\Delta^1|$ of Example \ref{E: collapse}. We can extend $\pi$  to a morphism of simplicial sets $\pi\colon\Delta^2\to\Delta^1$ by prescribing $\pi(0)=0$ and $\pi(1)=\pi(2)=1$. Then as in Example \ref{E: collapse},  $\Delta^2=[0,1,2]$ is taken to the degenerate simplex $[0,1,1]=s_1([0,1])$. At the same time, the morphism $\pi$ is doing an infinite number of other things: it takes the vertex $[0]\in \Delta^2$ to $[0]\in \Delta^1$, it takes the vertices $[1],[2]\in \Delta^2$ to $[1]\in \Delta^1$, it takes the $1$-simplex $[0,1]\in \Delta^2$ to $[0,1]\in \Delta^1$, it takes the $1$-simplex $[1,2]\in \Delta^2$ to the degenerate $1$-simplex\footnote{Careful: $[1]$ is a $0$-simplex, so $s_0$ is the appropriate (indeed the only well-defined) degeneracy map. Remember that $s_0$ tells us to repeat what occurs in the $0$th place - it doesn't know what's in that place.} $[1,1]=s_0[1]\in \Delta^1$, and it even takes the degenerate simplex $[0,1,1,2,2,2]=s_4s_3s_1[0,1,2]\in \Delta^2$ to the degenerate simplex 
$s_4s_3s_1[0,1,1]=[0,1,1,1,1,1]\in \Delta^1$. And much much more. 
\end{example}

\begin{example}Notice that, unlike simplicial maps on simplicial complexes,  morphisms on simplicial sets are not completely determined by what happens on vertices. For example, consider the possible simplicial morphisms from $\Delta^1$ to the simplicial set corresponding to the Delta set of Example \ref{E: nonuniqueness}. If we have a simplicial morphism that takes $[0]$ to $[v_0]$ and $[1]$ to $[v_1]$, there are still two possibilities for where to send $[0,1]$. 
\end{example}

\begin{example}
On the other hand, given a map of ordered simplicial complexes $f\colon X\to Y$, this induces a map of the associated simplicial sets as constructed in Example \ref{E: comp to set}. In this case, a function on vertices \emph{does} determine a simplicial map because simplices of ordered simplicial complexes are determined uniquely by their vertices. 
This was the case for the simplicial morphism of Example \ref{E: simp mor v}.
\end{example}

\begin{remark}\label{R: map faces}
Notice that it is always enough to define a simplicial morphism by what it does to nondegenerate simplices. What happens to the degenerate simplices is forced by the definition since, e.g. $f(s_i(x))=s_i(f(x))$. Similarly, what happens on faces is forced by what happens on the simplices of which they are faces. Thus, altogether,  simplicial morphisms can be described  by specifying  what they do to a comparatively small collection of nondegenerate simplices.
\end{remark}

From here on, we'll abandon the distinction between ``simplicial map'' and ``simplicial morphism'' and use the terms interchangeably as applied to simplicial sets.

\section{Realization}\label{S: realization}

If the idea of simplicial objects is to abstract from geometry/topology to combinatorics, there should be a way to reverse that process and turn simplicial sets into geometric/topological objects. Indeed that is the case. The definition looks a bit off-putting at first (what concerning simplicial sets doesn't?), but, in fact, we'll see that simplicial realization is a very natural thing to do.

\begin{definition}
Let $X$ be a simplicial set. Give each set $X_n$ the discrete topology and let $|\Delta^n|$ be the $n$-simplex with its standard topology. The \emph{realization} $|X|$ is given by $$|X|=\coprod_{n=0}^\infty X_n\times |\Delta^n|/\sim,$$
where $\sim$ is the equivalence relation generated by the relations $(x,D_i(p))\sim(d_i(x),p)$ for $x\in X_{n+1}, p\in |\Delta^n|$ and the relations $(x,S_i(p))\sim(s_i(x),p)$ for $x\in X_{n-1},p\in |\Delta^n|$. Here $D_i$ and $S_i$ are the face inclusions and collapses induced on the standard geometric simplices as in our discussion above of the category $\Delta$. 
\end{definition}

\begin{figure}[!htp]
\begin{center}
\scalebox{.6}{\includegraphics{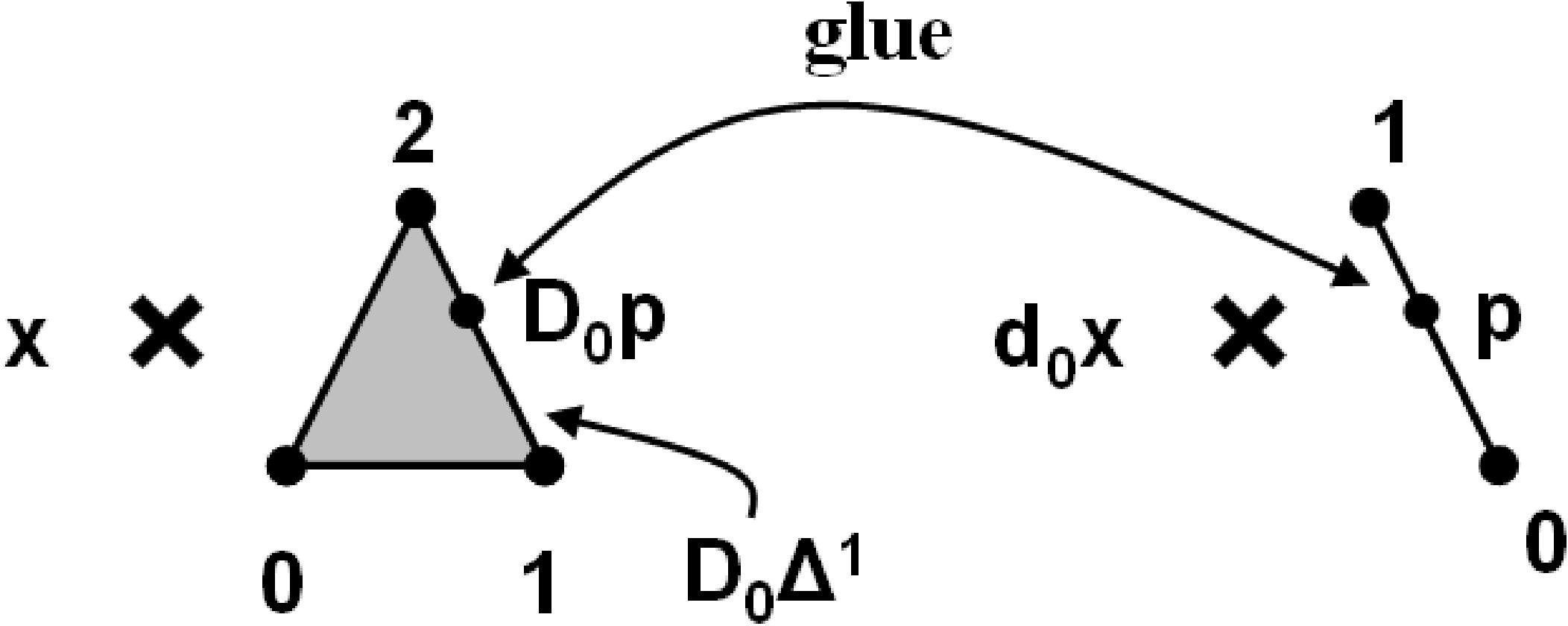}}
\end{center}
\caption{In the realization, the $1$-simplex representing $d_0x$, pictured on the right, is glued to the $2$-simplex representing $x$, pictured on the left, along the appropriate face.}\label{F: fig13}
\end{figure}

To see why this definition makes sense, let's think about how we would like to form a simplicial complex out of the data of a simplicial set. From the get-go, we have been thinking of the $X_n$ as collections of simplices. So this is just what $X_n\times |\Delta^n|$ is: a disjoint collection of simplices, one for each element of $X_n$. The next natural thing to do is to identify common faces. This is precisely what the relation $(x,D_i(p))\sim(d_i(x),p)$ encodes (see Figure \ref{F: fig13}): The first term of $(x, D_i(p))\subset (x,|\Delta^{n+1}|)$ is an $(n+1)$-simplex of $X$ and the second term $D_i(p)$ is a point on the $i$th face of a geometric $(n+1)$-simplex. On the other hand, $(d_i(x),p)$ is the $ith$ face of $x$ together with the same point, now in a stand-alone $n$-simplex. So the identification described just takes the $n$-simplex corresponding to $d_i(x)$ in $X_{n}\times |\Delta^n|$ and glues it as the $i$th face of the $(n+1)$-simplex assigned to $x$ in $X_{n+1}\times |\Delta^{n+1}|$. Since a similar gluing is done for any other $y$ and $j$ such that $d_j(y)=d_i(x)$, the effect is to glue faces of simplices together.

\begin{figure}[!htp]
\begin{center}
\scalebox{.6}{\includegraphics{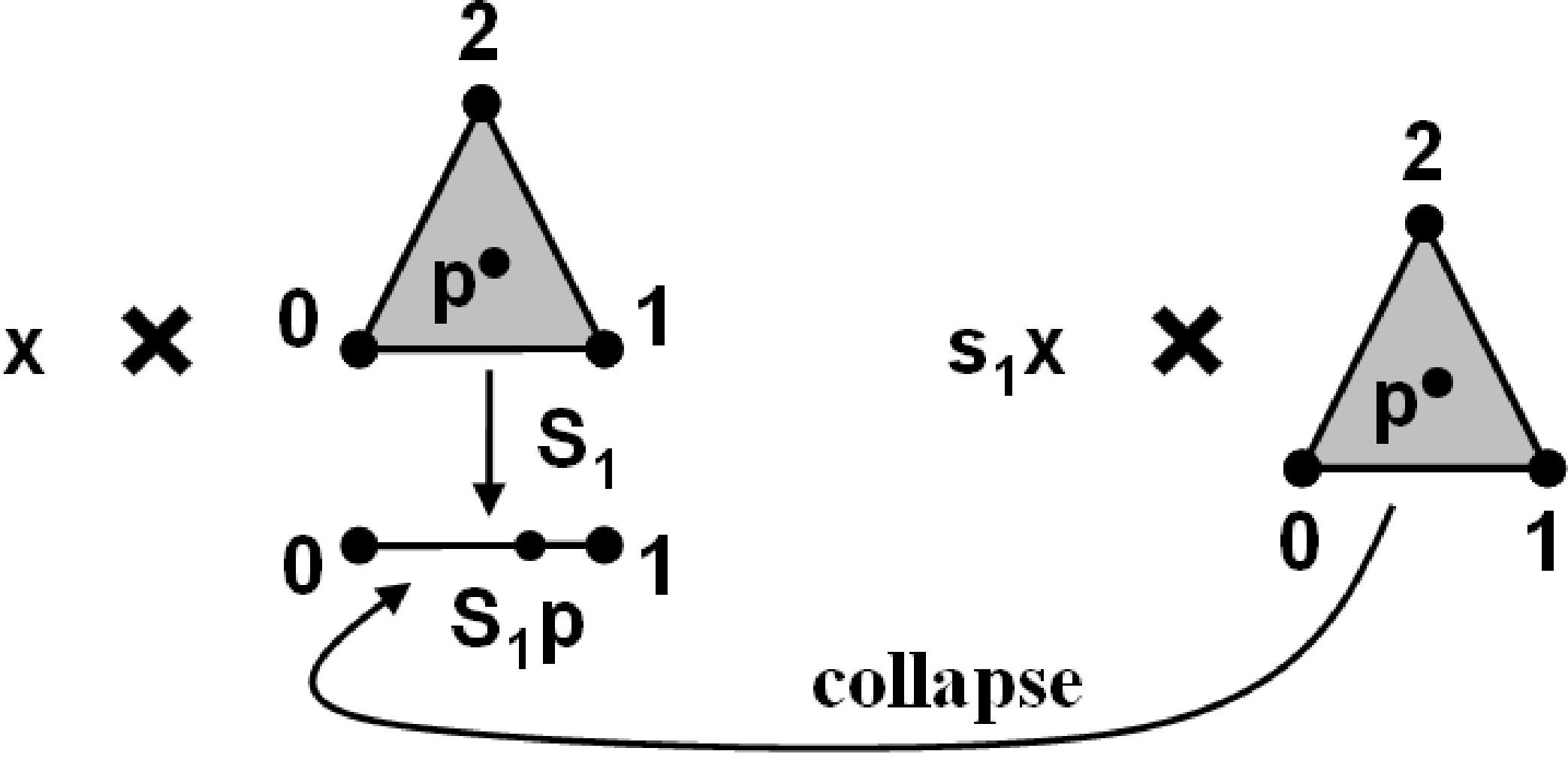}}
\end{center}
\caption{In the realization, the $2$-simplex representing $s_1x$, pictured on the right, is glued to the $1$-simplex representing $x$, pictured on the left, via the appropriate collapse, depicted by $S_1$.}\label{F: fig14}
\end{figure}

The next natural thing to do is suppress the degenerate simplices, since they're encoded within nondegenerate simplices anyway. This is what the relation $(x,S_i(p))\sim(s_i(x),p)$ for $x\in X_{n-1},p\in |\Delta^n|$ does, although more elegantly. This relation tells us that given a degenerate $n$-simplex $s_i(x)$ and a point $p$ in the pre-collapse $n$-simplex $|\Delta^n|$, we should glue $p$ to the $(n-1)$-simplex represented by $x$ at the point $S_i(p)$ in the image of the collapse map. That still sounds a little confusing, but the idea is straightforward: the $|\Delta^n|$ corresponding to degenerate $n$-simplices get collapsed in the natural way into the $(n-1)$-simplices they are degeneracies of. See Figure \ref{F: fig14}. We note also that there is no reason to believe that $x$ itself is nondegenerate. It might be, in which case the simplex corresponding to $x$ is itself collapsed. This provides no difficulty.

\begin{example}\label{E: 0}
Recall that the $0$-simplex $[0]$, thought of as a simplicial set, has one simplex $[0,\ldots,0]$ in each dimension $\geq 0$. Thus $|[0]|=\coprod_{i=0}^\infty |\Delta^i|/\sim$. So in dimension $0$ we have a single vertex $v$. In dimension $1$, we have a single simplex $[0,0]=s_0[0]$. The gluing instructions tell us to identify each $(s_0[0],p)=([0,0],p)\in ([0,0],|\Delta^1|)$
with $([0],S_0(p))=([0],v)$. Thus the $|\Delta^1|$ in dimension $1$ gets collapsed to the vertex. Similarly, since each point of the $2$-simplex $([0,0,0],|\Delta^2|)$ gets identified to a point of $([0,0],|\Delta^1|)$, and so on, we see that the whole situation collapses down to a single vertex. Thus $|[0]|$ is a point.
\end{example}

\begin{example}
Generalizing the preceding example, $|[0,\ldots, n]|=|\Delta^n|$ is just the standard geometric $n$-simplex, justifying our earlier use of notation.  We encourage the reader to explore this example on his or her own, noting that all of the degenerate simplices wind up tucked away within actual faces of $|\Delta^n|$, just where we expect them.
\end{example}

\begin{example}\label{E: return}
More generally, given any ordered simplicial complex, the realization of the simplicial set associated to it by adjoining all degenerate simplices (see Example \ref{E: comp to set}) returns the original simplicial complex.
\end{example}

\begin{example}
There is an analogous realization procedure for Delta sets. Given a Delta set $X$, we can define
the realization $|X|_{\Delta}$ by $$|X|_{\Delta}=\displaystyle\coprod_{n=0}^\infty X_n\times |\Delta^n|/\sim,$$
where $\sim$ is the equivalence relation generated by $(x,D_i(p))\sim(d_i(x),p)$ for $x\in X_{n+1}, p\in |\Delta^n|$. These realizations yield the types of spaces we have been drawing already to represent Delta sets. These are sometimes called Delta complexes; see, e.g., \cite{Ha}.

However, given a simplicial set $X$, the simplicial set realization of  $X$ is not generally going to be the same as the Delta set realization of the associated Delta set, say $X_{\Delta}$, that we obtain by neglect of structure. 

For example, consider the simplicial set $\Delta^0$. As seen in Example \ref{E: 0}, its simplicial realization, $|\Delta^0|$ is the topological space consisting of a single point. But recall that the simplicial set $\Delta^0$ has exactly one simplex in each dimension, and the neglect of structure that turns this into a Delta set $\Delta^0_\Delta$ drops the degeneracy relation but still leaves a Delta set with one simplex in each dimension and all face maps the unique possible ones. Thus the Delta set realization  $|\Delta^0_{\Delta}|_{\Delta}$ is an infinite dimensional CW complex with one cell in each dimension whose $n$-dimensional cell is attached by gluing each face of an $n$-simplex, in an order-preserving manner,   to the image of the unique $(n-1)$-simplex in the $(n-1)$-skeleton. Thus the $1$-skeleton of $|\Delta^0_{\Delta}|_{\Delta}$ is a circle, the $2$-skeleton is the ``dunce cap'' (see, e.g., \cite[Section 14]{BRTG}), and so on. This is evidently not homeomorphic to $|\Delta^0|$. However, it turns out that $|\Delta^0|$ and $|\Delta^0_{\Delta}|_{\Delta}$ \emph{are} homotopy equivalent; in fact $|\Delta^0_{\Delta}|_{\Delta}$ is contractible. In general, it is true that the realization of a simplicial set $|X|$ and the Delta set realization of its corresponding Delta set $|X_\Delta|_\Delta$  will be homotopy equivalent; see \cite{RSDelta}.

In what follows, discussion of ``realization'' and the notation $|X|$ will refer exclusively to simplicial set realization unless noted otherwise.
\end{example}

\begin{example}
Let $Y$ be a topological space, and let $\ms S(Y)$ be its singular set. $|\ms S(Y)|$ will be huge, with uncountably many simplices in each dimension (unless $Y$ is discrete - what will it be then?). While this looks discouraging, it turns out that  the natural map $|\ms S(Y)|\to Y$ (which acts on the realization  of each singular simplex by the map defining that singular simplex) induces isomorphisms on all homotopy  groups; see \cite[Theorem 4]{Mi57}. In particular, if $Y$ is a CW complex, this is enough to assure $|\ms S(Y)|$ and $Y$ are homotopy equivalent as a consequence of the Whitehead Theorem (see \cite[Corollary VII.11.14]{BRTG}), as we will see below in Theorem \ref{T: CW} that the realization of a simplicial set is always a CW complex. Thus, for many of the purposes of algebraic topology, $Y$ and $|\ms S(Y)|$ are virtually indistinguishable. So perhaps, wearing the appropriate glasses, $Y$ and $\ms S(Y)$ can be treated as the same thing, especially if $Y$ is a CW complex? We'll return to this idea later.
\end{example}

\begin{example}\label{E: sphere}
As noted in Example \ref{E: return}, the realization of a simplicial set that we obtained from an ordered simplicial complex is the original simplicial complex. So, for example, we can obtain a topological $(n-1)$-sphere as the realization of the boundary of the $n$-simplex, $\bd \Delta^n$. Here $\bd \Delta^n$ denotes the simplicial set obtained from the boundary $\bd |\Delta^n|$ of the ordered simplicial complex $|\Delta^n|$ by adjoining all degeneracies as in Example \ref{E: comp to set}. 
Let's find a good description of $\bd \Delta^n$  as a simplicial set. Since every $m$-simplex of $\bd \Delta^n$ should also be a simplex of $\Delta^n$, each can be written $[i_0,\ldots, i_m]$, where $0\leq i_0\leq \cdots\leq i_m\leq n$. The only caveat is that we do not allow any $m$-simplex that contains all of the vertices $0,\ldots, n$, since any such simplex would either be the ``top face'' $[0,\ldots, n]$, itself, or a degeneration of it, and these should not be faces of $\bd \Delta^n$. In summary, then,  $\bd \Delta^n$ is the simplicial set consisting of all nondecreasing sequences of the numbers $0,\ldots,n$ that do not contain all of the numbers $0,\ldots ,n$, and since this is the simplicial set arising from the ordered simplicial complex $\bd |\Delta^n|$, we have $|\bd \Delta^n|\cong S^{n-1}$.

Is this the most efficient way to obtain $S^{n-1}$ as the realization of a simplicial set? After all, $\bd \Delta^n$ contains quite a number of simplices, many of which are nondegenerate (the interested reader might go and count them). Here is another way to do it, at least for $n\geq 2$, suggested by CW complexes. Let $X$ be a simplicial set whose only  nondegenerate simplices are denoted by $[0]\in X_0$ and $[0,\ldots, n-1]\in X_{n-1}$. All simplices in $X_i$, $0<i<n-1$, are the degenerate simplices $[0,\ldots, 0]$. This, of course, forces all of the faces of $[0,\ldots, n-1]$ to be $[0,\ldots, 0]$, and we see that the realization $|X|$  is equivalent to the standard construction of $S^{n-1}$ as a CW complex by collapsing the boundary of an $(n-1)$-cell to a point. See Figure \ref{F: fig16}. 

\begin{figure}[!htp]
\begin{center}
\scalebox{.6}{\includegraphics{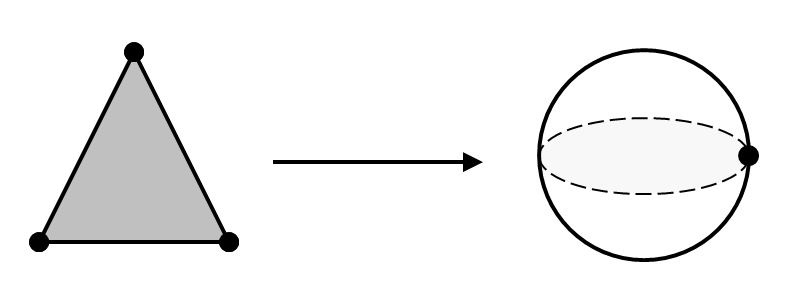}}
\end{center}
\caption{The realization of the simplicial set consisting of only two nondegenerate simplices, one in dimension $0$ and the other in dimension $2$, is the sphere $S^2$; this picture represents the image of the nondegenerate simplex of dimension $2$ in the realization. The entire boundary of the $2$-simplex is collapsed to the unique $0$-simplex. }\label{F: fig16}
\end{figure}

\end{example}

The preceding example is instructive on several different points:

\begin{enumerate}
\item The second part of Example \ref{E: sphere} relies strongly on the existence of degenerate simplices. For $n>2$, we cannot construct $S^{n-1}$ this way as the realization of a Delta set. A Delta set with an $(n-1)$-simplex would require actual (nondegenerate) $(n-2)$-simplices as its faces. Of course we can still get $S^{n-1}$ as the realization of the Delta set corresponding to $\bd \Delta^n$.  

\item Notice that the realization of a simplicial set does not necessarily inherit the structure of a simplicial complex, at least not in any obvious way from the data of the simplicial set. 

\item Realizations are non-unique, in the sense that very different looking simplicial sets can have the same geometric realization up to homeomorphism. This is not surprising, since there are many ways to triangulate a piecewise-linear space.

\end{enumerate}

Example \ref{E: sphere} is also disconcerting in that the reader may be getting worried that realizations of simplicial sets might be very complicated to understand with all of the gluing and collapsing that can occur. To mitigate these concerns somewhat, we first observe that all degenerate simplices do get collapsed down into the simplices of which they are degeneracies, and so constructing a realization depends only on understanding what happens to the nondegenerate simplices. A second concern would be that two nondegenerate simplices might be glued together. This would happen if it were possible for two nondegenerate simplices to have a common degeneracy (why?). Luckily, this does not happen, as we demonstrate in the following proposition. As a corollary, we can conclude that the realization of a simplicial set is made up of the disjoint union of the interiors of the nondegenerate simplices. We must limit this statement to the interiors as the faces of a nondegenerate simplex may be degenerate, as in the second part of Example \ref{E: sphere} - meanwhile, nondegenerate faces will look out for themselves!

\begin{proposition}\label{P: unique deg}
A degenerate simplex is a degeneracy of a unique nondegenerate simplex. In other words, if $z$ is a degenerate simplex, then there is a unique nondegenerate simplex $x$ such that $z=s_{i_1}\cdots s_{i_k}x$, for some collection of degeneracy maps $s_{i_1},\ldots, s_{i_k}$. 
\end{proposition}
\begin{proof}
Suppose  $z$ is a degenerate simplex. Then $z=s_{i_1}x_1$ for some $x_1$ and some degeneracy map $s_{i_1}$. If $x_1$ is degenerate, we can make a similar replacement and  continue inductively until eventually we have  $z=s_{i_1}\cdots s_{i_k}x_k$ for some nondegenerate $x_k$. The process stops because each successive $x_j$ has lower dimension than the preceding, and there are no simplices of dimension less than zero. Thus $z$ can be written in the desired form.

Next, suppose $x$ and $y$ are nondegenerate simplices, possibly of different dimensions, and that $Sx=Ty$, where $S$ and $T$ are compositions of degeneracy operators.  Suppose $S=s_{i_1}\cdots s_{i_k}$. Let $D=d_{i_k}\cdots d_{i_1}$. Then $x=DSx=DTy$, using the simplicial set axioms for the first equality. By using the simplicial set axioms to trade face maps to the right, we obtain $x=\td T\td D y$ for some composition of face operators $\td D$ and some composition of degeneracies $\td T$. But, by hypothesis, $x$ is nondegenerate, so $\td T$ must be vacuous, and we must have $x=\td Dy$. That is $x$ is a face of $y$. But we could repeat the argument reversing $x$ and $y$ to obtain that $y$ is also face of $x$. But this is impossible unless $x=y$.  
\end{proof}

Another comforting fact is the following theorem:

\begin{theorem}\label{T: CW}
If $X$ is a simplicial set, then $|X|$ is a CW complex with one $n$-cell for each nondegenerate $n$-simplex of $X$. 
\end{theorem}
\begin{proof}
 We refer to Milnor's paper on geometric realization \cite{Mi57} (or, alternatively, to \cite[Theorem 14.1]{MAY67}) for the proof, which is not difficult and which formalizes our discussion preceding Proposition \ref{P: unique deg}.
\end{proof}

\paragraph{The adjunction relation.}

The realization functor $|\cdot|$ turns out to be adjoint to the singular set functor $\ms S(\cdot )$. 

\begin{theorem}\label{T: adjoint}
If $X$ is a simplicial set and $Y$ is a topological space, then $$\text{\emph{Hom}}_{\textbf{\emph{Top}}}(|X|,Y)\cong \text{\emph{Hom}}_{\mbf S}(X,\ms S(Y)),$$
where $\text{\emph{Hom}}_{\mbf S}$ denotes morphisms of simplicial sets and $\text{\emph{Hom}}_{\textbf{\emph{Top}}}$ denotes continuous maps of topological spaces.  
\end{theorem}
\begin{proof}[Sketch of proof]
We identify the two maps $\Psi\colon\Hom_{\Top}(|X|,Y)\to \Hom_{\mbf S}(X,\ms S(Y)) $ and $\Phi\colon\Hom_{\mbf S}(X,\ms S(Y))\to \Hom_{\Top}(|X|,Y)$ and leave it to the reader both to check carefully that these are well-defined and to show that they are mutual inverses.

A map $f\in \Hom_{\mbf S}(X,\ms S(Y))$ assigns to each $n$-simplex  $x\in X$ a continuous function $\sigma_x\colon |\Delta^n|\to Y$. Let $\Phi(f)$ be the continuous function that acts on the simplex $(x,|\Delta^n|)\in |X|$ by applying $\sigma_x$ to $|\Delta^n|$. 

Conversely, given a function $g\in \Hom_{\Top}(|X|,Y)$, then the restriction of $g$ to a nondegenerate simplex $(x, |\Delta^n|)$  yields a continuous function $|\Delta^n|\to Y$ and thus an element of $\ms S(Y)_n$. If $(x,| \Delta^n|)$ represents a degenerate simplex, then we precompose with the appropriate collapse map of $\Delta^n$ into $|X|$  before applying $g$. 
\end{proof}

One can say much more on the relation between simplicial sets and categories of topological spaces. For example, see Theorem \ref{T: equiv} below, according to which the homotopy category of CW complexes is equivalent to the homotopy category of simplicial sets satisfying a condition called the \emph{Kan condition}. The Kan condition is defined in Section \ref{S: Kan}.

\section{Products}\label{S: product}

Before we move on to  look at simplicial homotopy, we will need to know about products of simplicial sets. For those accustomed to products of simplicial complexes or products of chain complexes, the definition of the product of simplicial sets looks surprisingly benign by comparison.

\begin{definition}\label{D: product}
Let $X$ and $Y$ be simplicial sets. Their product $X\times Y$ is defined by 
\begin{enumerate}
\item $(X\times Y)_n=X_n\times Y_n=\{ (x,y)\mid x\in X_n, y\in Y_n\}$,

\item if $(x,y)\in (X\times Y)_n$, then $d_i(x,y)=(d_ix,d_iy)$,

\item if $(x,y)\in (X\times Y)_n$, then $s_i(x,y)=(s_ix,s_iy)$.
\end{enumerate}
\end{definition}

Notice that there are evident projection maps $\pi_1\colon X\times Y\to X$ and $\pi_2\colon X\times Y\to Y$ given by $\pi_1(x,y)=x$ and $\pi_2(x,y)=y$. These maps are clearly simplicial morphisms.

Definition \ref{D: product} looks disturbingly simple-minded, but it is vindicated by the following important theorem.

\begin{theorem}\label{T: product}
If $X$ and $Y$ are simplicial sets, then $|X\times Y|\cong |X|\times |Y|$ (in the category of compactly generated Hausdorff spaces). In particular, if $X$ and $Y$ are countable or if one of $|X|,|Y|$ is locally finite as a CW complex, then $|X\times Y|\cong |X|\times |Y|$ as topological spaces.
\end{theorem}
We refer the reader to \cite[Theorem 14.3]{MAY67} or \cite{Mi57} for a proof in the latter situations and to \cite[Chapter III]{GabZis} for a proof of the general case. However, since an example is perhaps worth a thousand proofs, we will take a detailed look at some special cases.

\begin{example}
Let $X$ be any simplicial set, and let $Y=\Delta^0=[0]$. Since $\Delta^0$ has a unique element in each dimension, $X\times \Delta^0\cong X$. So indeed, $|X\times \Delta^0|\cong |X|\times |\Delta^0|\cong  |X|$. 
\end{example}

\begin{example}
The first interesting example is $\Delta^1\times \Delta^1$. We would like to see that $|\Delta^1\times \Delta^1|\cong |\Delta^1|\times |\Delta^1|$, the square.  As discussed in Section \ref{S: realization}, we need to focus on the nondegenerate simplices of $\Delta^1\times \Delta^1$. The reader can refer to Figure \ref{F: fig21} for the following discussion.

\begin{figure}[!htp]
\begin{center}
\scalebox{.6}{\includegraphics{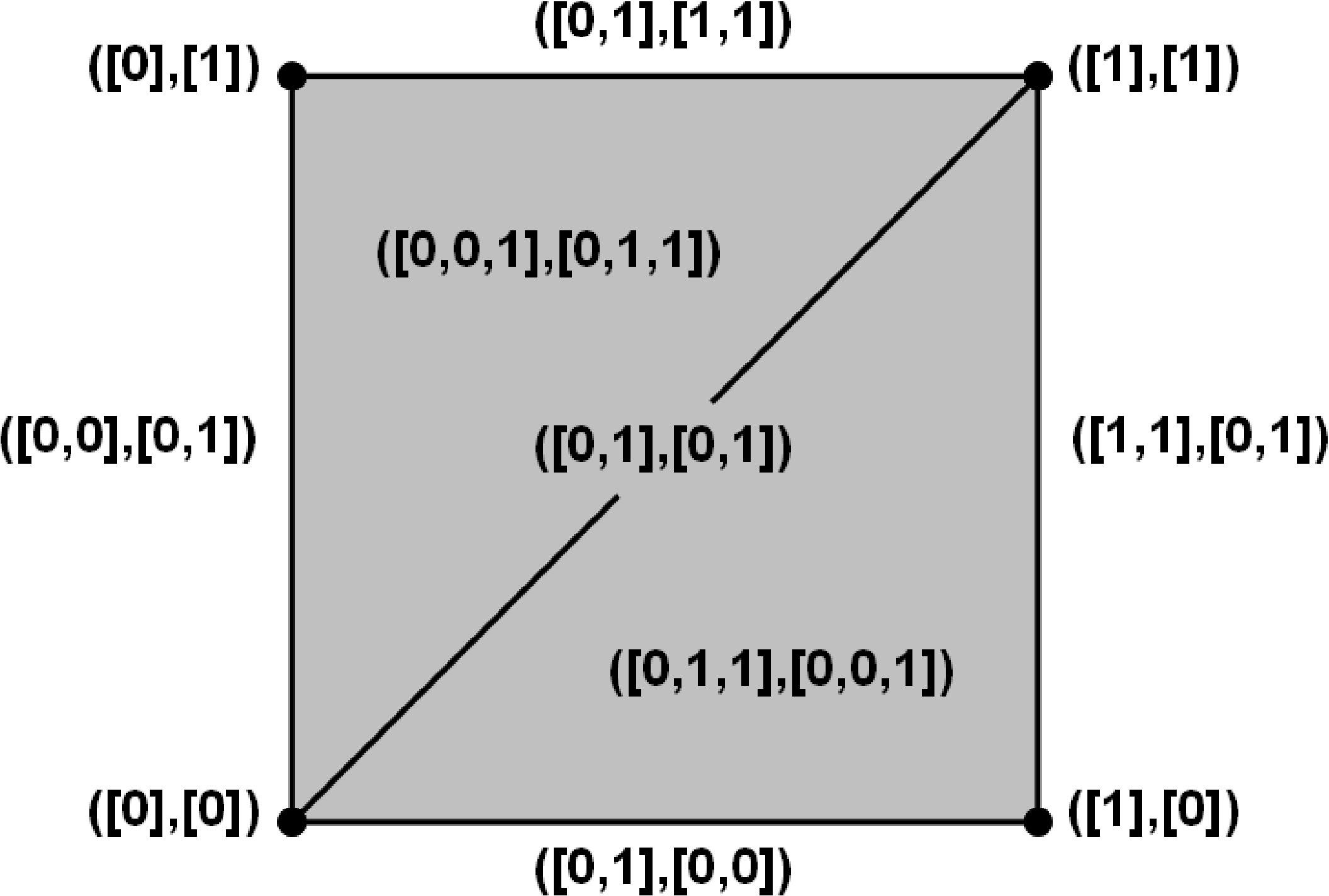}}
\end{center}
\caption{The realization of  $\Delta^1\times \Delta^1$}\label{F: fig21}
\end{figure}

First, in dimension $0$, we have the product $0$-simplices $$X_0=\{([0],[0]),([1],[0]),([0],[1]),([1],[1])\},$$ the four vertices of the square. 

In dimension $1$, we have the pairs $(e,f)$, where $e$ and $f$ are $1$-simplices of $\Delta^1$. There are three possibilities for each of $e$ and $f$ - $[0,0]$, $[0,1]$, and $[1,1]$. So there are nine $1$-simplices of $\Delta^1\times \Delta^1$.

There is only one $1$-simplex that is made up completely of nondegenerate simplices: $([0,1],[0,1])$. Since $d_0([0,1],[0,1])=(1,1)$ and $d_1([0,1],[0,1])=(0,0)$, the simplex $([0,1],[0,1])$ must be the diagonal. Those with one nondegenerate and one degenerate $1$-simplex are  $([0,0],[0,1])$,  $([0,1],[0,0])$, $([1,1],[0,1])$ and $([0,1],[1,1])$, which, as we see by checking the endpoints, are respectively the left, bottom, right, and top of the square. The other four $1$-simplices are the degeneracies of the vertices. For example, $([0,0],[1,1])=(s_0[0],s_0[1])=s_0([0],[1])$.

Now for the $2$-simplices - here's where things get a little tricky. There are four $2$-simplices of $\Delta^1$: $[0,0,0]$, $[0,0,1]$, $[0,1,1]$, and $[1,1,1]$.  So there are sixteen $2$-simplices of $\Delta^1\times \Delta^1$. There are two possible degeneracy maps,  $s_0$ and $s_1 $, from $(\Delta^1\times \Delta^1)_1$ to $(\Delta^1\times \Delta^1)_2$. These act on the nine $1$-simplices, but there are not eighteen degenerate $2$-simplices since $s_0s_0=s_1s_0$, and we know there are four degenerate $1$-simplices $s_0v_i$ of $\Delta^1\times \Delta^1$ corresponding to the degeneracies of the four vertices. Removing these redundancies leaves fourteen degenerate $2$-simplices. There are no other redundancies since $s_0s_0=s_1s_0$ is the only relation on $s_1$ and $s_0$. The remaining two $2$-simplices are nondegenerate. These turn out to be $([0,0,1],[0,1,1])$ and $([0,1,1],[0,0,1])$, which are the two triangles, as one can check by computing face maps. 

Next, we need to  see that all $3$-simplices and above of $\Delta^1\times \Delta^1$  are degenerate. We first observe that each $3$-simplex of $\Delta^1$  must be a double degeneracy of a $1$-simplex (since there are no nondegenerate simplices of $\Delta^1$ of dimension greater than $1$). But there are only six such options, of the forms $s_0s_0e$, $s_0s_1e$, $s_1s_0e$, $s_1s_1e$, $s_2s_0e$, and $s_2s_1e$ for a (possibly degenerate) $1$-simplex $e$. However, the simplicial set axioms reduce this to the possibilities $s_1s_0e$, $s_2s_0e$, and $s_2s_1e$. But then, again by the axioms, 
\begin{align*}
(s_1s_0e  ,s_1s_0f)&= s_1(s_0e,s_0f) \\
(s_1s_0e  ,s_2s_0f)&=(s_0s_0e  ,s_0s_1f)=s_0(s_0e  ,s_1f) \\
(s_1s_0e  ,s_2s_1f)&= (s_1s_0e  ,s_1s_1f)= s_1(s_0e  ,s_1f) \\
(s_2s_0e  ,s_1s_0f)&= (s_0s_1e    , s_0s_0f  )=s_0(s_1e    , s_0f  )\\
(s_2s_0e  ,s_2s_0f)&=s_2(s_0e  ,s_0f)   \\
(s_2s_0e  ,s_2s_1f)&= s_2(s_0e  ,s_1f)  \\
(s_2s_1e  ,s_1s_0f)&=(s_1s_1e  ,s_1s_0f)= s_1(s_1e  ,s_0f) \\
(s_2s_1e  ,s_2s_0f)&= s_2(s_1e  ,s_0f)  \\
(s_2s_1e  ,s_2s_1f)&= s_2(s_1e  ,s_1f).  
\end{align*}

So all $3$-simplices of $\Delta^1\times \Delta^1$ are degenerate. It also follows that all higher dimension simplices are degenerate:  the terms in any such product  must be further degeneracies of these particular doubly degenerate $1$-simplices, and using the simplicial set axioms, we can move $s_0$ and $s_1$ to the left in all expressions. Then we can proceed as in the above list of computations.
\end{example}

\smallskip

That last bit isn't very intuitive, but the low-dimensional part makes some sense. If we take the product of two CW complexes, the cells of the product will be product  cells of the form $C_1\times C_2$, where $C_1 $ and $C_2$ are not necessarily of the same dimension. These mixed dimensional cells occur here as products of nondegenerate simplices with degenerate simplices. What makes matters difficult is that we must preserve a simplicial structure. This forced ``triangulation'' is what makes matters somewhat complicated.

It will be useful for us to look even more closely at the products $\Delta^p\times \Delta^q$. After all, all products will be made up of these building blocks. The main point of interest for us is that the simplicial product construction yields the same triangulation structure that may be familiar from homotopy arguments in courses in beginning algebraic topology. 

\begin{example}\label{E: prism}
Suppose $p,q>0$. 
Since we know that $|\Delta^p\times \Delta^q|=|\Delta^p|\times |\Delta^q|$, let us focus on the nondegenerate $(p+q)$-simplices  of $\Delta^p\times \Delta^q$. We let $E_j$ stand for the unique nondegenerate $j$-simplex of $\Delta^j$. We note immediately that any nondegenerate $(p+q)$-simplex $s$ of $\Delta^p\times \Delta^q$ (and hence the only ones that appear nondegenerately in the realization) must have the form $s=(SE_p, S'E_q)$, where $S$ and $S'$ are sequences of degeneracy maps. Why? Otherwise $s$ would have to be of the form  $s=(\bar St, \bar S't')$, where $\bar S$ and $\bar S'$ are again sequences of degeneracy maps and $t$ and $t'$ are faces of $E_p$ and $E_q$, respectively, at least one of which is a proper face. But in this case, we would have $s\in F\times F'$, where $F$ and $F'$ are the simplicial subsets corresponding to faces of $\Delta^p$ and $\Delta^q$, at least one of which is a proper face. Consequently the image of $s\times |\Delta^{p+q}|$ in the realization of $\Delta^p\times \Delta^q$ will in fact lie within the realization $|F|\times |F'|$. In other words, $s$ is a simplex of some $\Delta^r\times \Delta^s$ with $r+s<p+q$, and this will imply that $s$ must actually be a degenerate simplex.  We invite the reader to think through why by generalizing the above argument that all $m$-simplices, $m\geq 3$, of $\Delta^1\times \Delta^1$ are degenerate (alternatively, $|F|\times |F'|$ has geometric dimension less than $p+q$ and so can contain no $(p+q)$-dimensional subspace).

So now we see that  $s=(SE_p, S'E_q)$, and for dimensional reasons, we can write this as $s=(s_{i_q}\cdots s_{i_{1}}E_p, s_{j_p}\cdots s_{j_{1}}E_q)$. Furthermore, using the simplicial set axioms, we can assume that $0\leq i_1< \cdots < i_q< p+q$ and $0\leq j_1< \cdots <j_q< p+q$. Now notice that the collection $\{i_1,\ldots, i_q,j_1,\ldots, j_p\}$ consists of $p+q$ numbers from $0$ to $p+q-1$. Furthermore, there can be no redundancy, since if $i_k=j_{k'}$ for some $k$ and $k'$, then again by the axioms, we can pull these indices to the front to get $s=(s_{i}\td S E_p,s_{i}\td S'E_q)=s_{i}(\td S E_p,\td S'E_q)$ for some $i, \td S,\td S'$, making $s$ degenerate. 

Thus we conclude that the nondegenerate $(p+q)$-simplices of $\Delta^p\times \Delta^q$ are precisely those of the form $s=(s_{i_q}\cdots s_{i_{1}}E_p, s_{j_p}\cdots s_{j_{1}}E_q)$, where the $i_k$ and $j_k$ are increasing series of integers from $0$ to $p+q-1$, all completely distinct.

In the special case $\Delta^p\times \Delta^1=\Delta^p\times I$, this rule for nondegenerate $(p+1)$-dimensional simplices reduces to the form $s=(s_i E_p, s_{j_p}\cdots s_{j_1} e)$, where $e$ is the edge $[0,1]$ of $I$, and the sequence $j_1,\ldots, j_p$ is increasing from $0$ to $p$, omitting only $i$. Thus there are precisely $p+1$ nondegenerate $(p+1)$-simplices. Since $e=[0,1]$, notice that all of the degeneracy maps before the ``gap'' at $i$ must adjoin another $0$ and all of those after the ``gap'' adjoin more $1$s. Thus we can also label these nondegenerate $(p+1)$-simplices exactly by the $p+1$ sequences of length $p+2$ of the form $[0,\ldots, 0,1,\ldots,1]$ that must start with a $0$ and end with a $1$.

If this looks familiar, it's because the standard way to triangulate the product prism $\Delta^p\times I$ when studying simplicial homology theory is by the $(p+1)$-simplices $[0,\ldots, k,k',\ldots, p']$, where  the unprimed numbers represent vertices in $\Delta^p\times 0$ and the primed numbers represent vertices in $\Delta^p\times 1$. The simplex $[0,\ldots,k,k',\ldots, p']$ corresponds to $k+1$ zeros and $p-k+1$ ones. See Figure \ref{F: fig22}.

\begin{figure}[!htp]
\begin{center}
\scalebox{.6}{\includegraphics{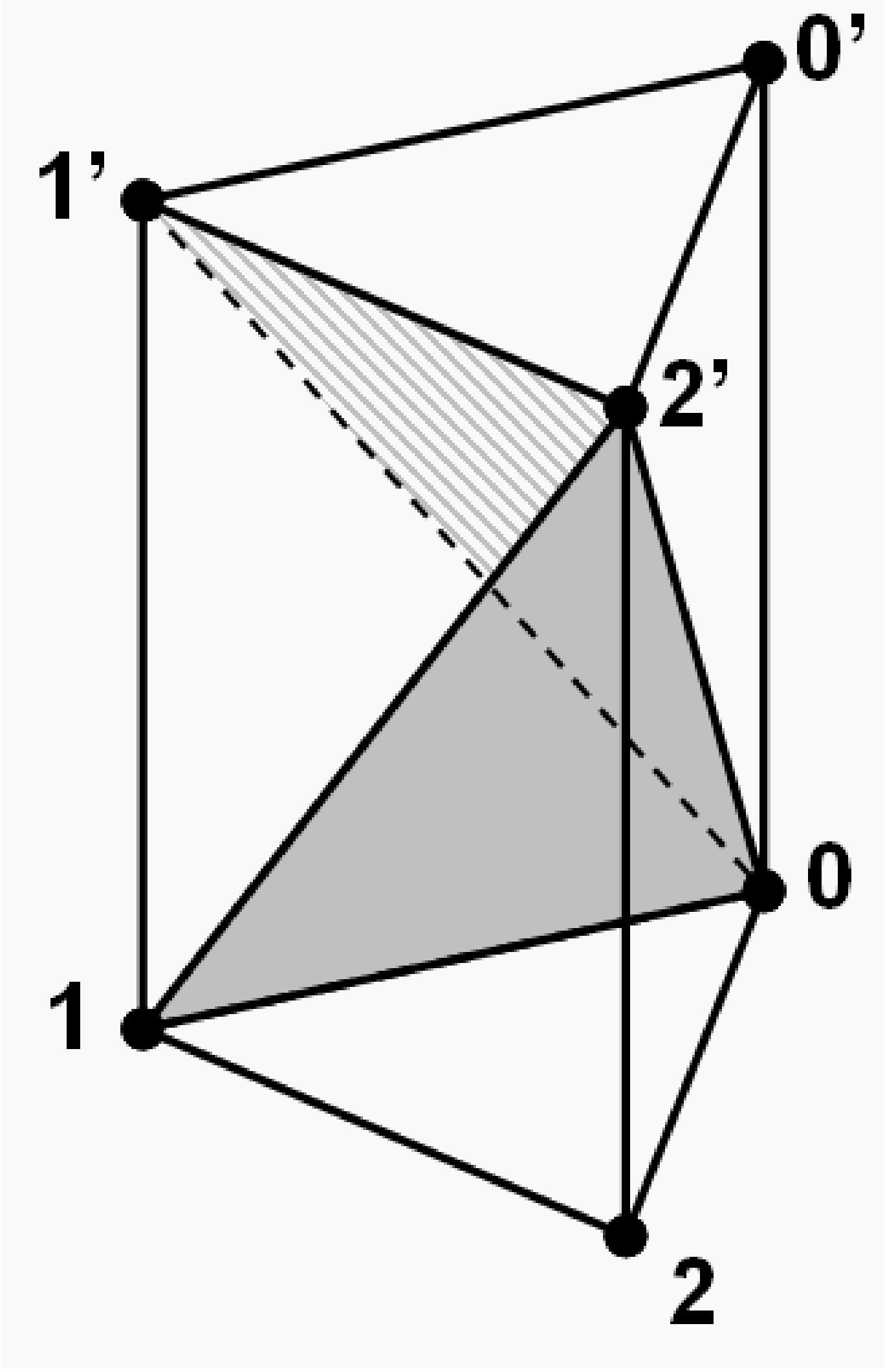}}
\end{center}
\caption{The realization of $|\Delta^2\times \Delta^1|$ with nondegenerate $3$-simplices $[0,0',1',2']$, $[0,1,1',2']$, and $[0,1,2,2']$}\label{F: fig22}
\end{figure}

For our upcoming discussion of simplicial homotopy, it's also worth looking at how these simplices are joined together along their boundaries. Let's first look from the point of view of writing the $(p+1)$-simplices of $\Delta^p\times I$ in the form $P_k=[0,\ldots, k,k',\ldots, p']$, where $0\leq k\leq p$. If $i<k$, then $d_iP_k=[0,\ldots,i-1,i+1,\ldots,k,k',\ldots, p']$. But this can be thought of as a $p$-simplex of $[0,\ldots,i-1,i+1,\ldots, p]\times I$ and so is part of the boundary $(\bd \Delta^p)\times I$. Similar considerations hold if $i>k+1$. The interesting ``interior cases'' are
\begin{align*}
d_kP_k&=[0,\ldots,k-1,k',\ldots,p']\\
d_{k+1}P_k&=[0,\ldots,k,(k+1)',\ldots,p'].
\end{align*}
To understand the assembly of the prism $\Delta^p\times I$ from the $P_k$,  notice that $d_kP_k=d_kP_{k-1}$ for $k>0$ and $d_{k+1}P_k=d_{k+1} P_{k+1}$ for $k<p$. This tells us how to glue the $(p+1)$-simplices together to form $|\Delta^p\times I|$. 

In our other notation, if we have 
$P_k=(s_k E_p, s_{p}\cdots s_{k+1}s_{k-1} \cdots s_{0} e)$, then for $i<k$ we have, using the axioms, $$d_iP_k=(s_{k-1}d_i E_p, s_{p-1}\cdots  s_{k}s_{k-2}\cdots s_i (d_is_i)s_{i-1}\cdots s_{0} e)=(s_{k-1}d_i E_p, s_{p-1}\cdots s_{k}s_{k-2} \cdots s_{0} e).$$ Notice that we use the axioms to   ``pass $d_i$ through,'' converting each $s_j$ to $s_{j-1}$ along the way, until it ``annihilates'' with the original $s_i$ (leaving the previous $s_{i+1}$ converted to the new $s_i$). We wind up with a $p$-simplex that is recognizable as a $p$-simplex in $d_iE_p\times I$. Similarly, for $i>k+1$, we get $d_iP_k=(s_{k}d_{i-1} E_p, s_{p-1}\cdots s_{k+1}s_{k-1} \cdots s_{0} e)$. The two ``interior'' cases correspond to $d_kP_k$ and $d_{k+1}P_k$:
\begin{align*}
d_kP_k&=(d_ks_k E_p, s_{p-1}\cdots s_{k}s_{k-2} \cdots s_{0} e)=(E_p, s_{p-1}\cdots s_{k}s_{k-2} \cdots s_{0} e)\\
d_{k+1}P_k&=(d_{k+1}s_k E_p, s_{p-1}\cdots s_{k+1}s_{k-1} \cdots s_{0} e)=(E_p, s_{p-1}\cdots s_{k+1}s_{k-1} \cdots s_{0} e).
\end{align*}
These are not in $\bd \Delta^p\times I$. However, we do again see that $d_kP_k=d_kP_{k-1}$ for $k>0$ and $d_{k+1}P_k=d_{k+1} P_{k+1}$ for $k<p$.
\end{example}

\subsection{Simplicial $\Hom$}\label{S: Hom}

We have just seen that there is a product functor internal to the category of simplicial sets; in other words the product of two simplicial sets is again a simplicial set. Many other important categories in algebraic topology possess an analogous internal product functor for which the product of two objects in the category is again an object of that category. Examples include the category of sets, the category of topological spaces, and the category of bimodules over a commutative ring $R$ (for which the appropriate product is the tensor product). In these categories there are important interplays (via adjunction) between the product functor and an internal $\Hom$ functor, which also takes a pair of objects of the category to an object of the category (in our examples, the set of set maps, the space of maps of spaces (with an appropriate topology), or the $R$-module of $R$-module homomorphisms). This is a feature shared by the category of simplicial sets, in which it is possible to define an internal $\Hom$ functor that takes two simplicial sets $X,Y$ and outputs a simplicial set $\Homs(X,Y)$. This simplicial set extends the set of morphisms between two simplicial sets, which occurs as the set of vertices $\Homs(X,Y)_0=\Hom_{\mbf S}(X,Y)$. We will briefly describe the construction in this section; however, a detailed study of the internal $\Hom$ functor would take us too far afield, so we simply provide the basic definitions  and leave a more detailed treatment to other sources\footnote{Note that notation for $\Homs(X,Y)$ varies widely across sources. Another common notation is $Y^X$. }. 

As a motivation for the definition of $\Homs(X,Y)$, recall the adjunction relation between products and $\Hom$ functors in our other familiar categories: $$\Hom(A\times B,C)\cong \Hom(A,\Hom(B,C)).$$ We would like to end up with something similar for simplicial sets. Furthermore, notice that for any simplicial set $Z$, the set of $n$-simplices $Z_n$ can be identified with the set of simplicial morphisms $\Hom_{\mbf S}(\Delta^n, Z)$. So whatever the simplicial set $\Homs(X,Y)$ is, it must satisfy $\Homs(X,Y)_n=\Hom_{\mbf S}(\Delta^n, \Homs(X,Y))$. But in a category for which an adjunction relation holds, we would then hope to be able to identify this expression with something like $\Hom_{\mbf S}(\Delta^n\times X,Y)$. Note that this is not yet a property we can \emph{check} because we do not yet have a definition of $\Homs(X,Y)$. However, we can turn it around into a definition by \emph{defining} $$\Homs(X,Y)_n= \Hom_{\mbf S}(\Delta^n\times X,Y).$$

Setting $\Homs(X,Y)_n= \Hom_{\mbf S}(\Delta^n\times X,Y)$ gives us the simplices of $\Homs(X,Y)$. If $f\in \Homs(X,Y)_n$, we obtain its $i$th face $d_if\in \Homs(X,Y)_{n-1}= \Hom_{\mbf S}(\Delta^{n-1}\times X,Y)$ as the composite $$\Delta^{n-1}\times X\xrightarrow{D_i\times \text{id}}\Delta^{n}\times X\xrightarrow{f}Y.$$ The degeneracy maps are defined analogously. 

With this definition, one can check that the following adjunction relationship holds in the category of simplicial sets:
$$\Hom_{\mbf S}(Z,\Homs(X,Y))\cong \Hom_{\mbf S}(Z\times X,Y).$$

Furthermore, it follows that $$\Homs(Z,\Homs(X,Y))\cong \Homs(Z\times X,Y).$$

For an excellent discussion of these various $\Hom$ objects and adjunctions from the more general point of view of simplicial model categories, see \cite[Section II.2]{GoeJar}.

\section{Simplicial objects in other categories}\label{S: simplicial cats}

Before moving on to discuss simplicial homotopy, we pause to note that the categorical definition of simplicial sets suggests a sweeping generalization.

\begin{definition}
Let $\Cat$ be a category. A \emph{simplicial object in $\Cat$} 
is a contravariant functor $X\colon \Delta\to \Cat$ (equivalently, a covariant functor $X\colon  \Delta^{op}\to \Cat$). A morphism of simplicial objects in $\Cat$ is a natural transformation of such functors.
\end{definition}

Another common notation, when $\Cat$ is a familiar category with objects of a given type, is to refer to a simplicial object in $\Cat$ as a simplicial [insert type of object]. In other words, when $\Cat$ is the category of groups and group homomorphisms, we speak of simplicial groups. This is consistent with referring to a simplicial object in the category $\Set$ as a simplicial set. One also commonly encounters simplicial $R$-modules, simplicial spaces, and even simplicial categories!

\begin{figure}[!htp]
\begin{center}
\scalebox{.6}{\includegraphics{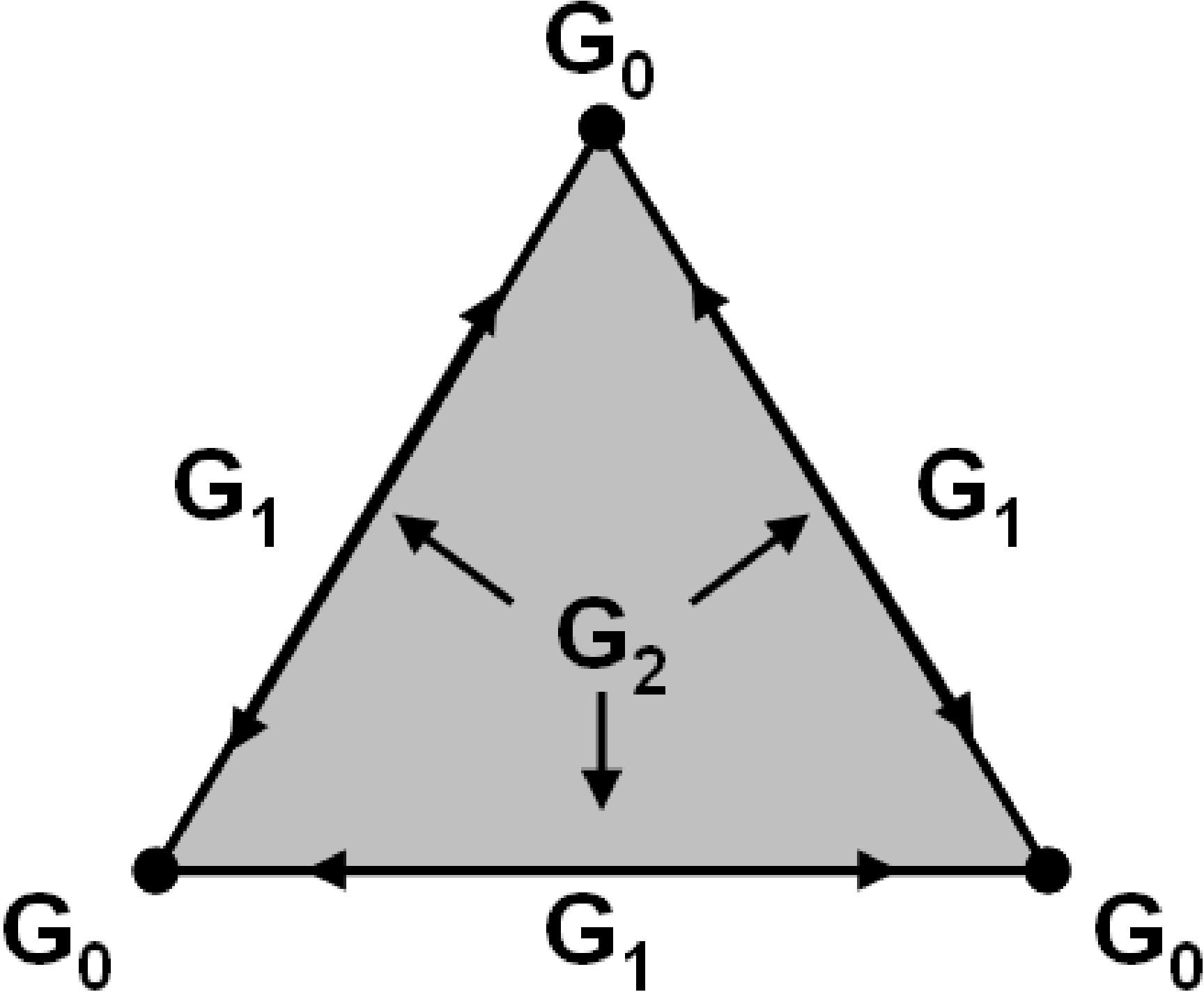}}
\end{center}
\caption{A pictorial representation of a $2$-simplex of a simplicial group with arrows representing the face morphisms from dimension $2$ to dimension $1$ and from dimension $1$ to dimension $0$}\label{F: fig17}
\end{figure}

\begin{example}
Let's unwind the definition in the case of simplicial groups. By definition, a simplicial group $\mc G$ consists of a sequence of groups $\mc G_n$ and collections of group homomorphisms $d_i\colon\mc G_n\to \mc G_{n-1}$ and $s_i\colon\mc G_n\to\mc G_{n+1}$, $0\leq i\leq n$, that satisfy the axioms \eqref{E: axioms}.

At this point, unfortunately, trying to picture group elements as simplices breaks down a little bit since there is so much extra structure around (what does it mean geometrically to multiply two simplices?). Nonetheless, it is still helpful to refer mentally to the category $\Delta$, in which we can visualize each simplex $[n]$ as representing  a group and picture movement toward each $n-1$ face as representing a different group homomorphism to the group represented by $[n-1]$. See Figure \ref{F: fig17}.
\end{example}

\begin{example}
Suppose $X$ is a simplicial set. Then we can form the simplicial group $C_*(X)$ with $(C_*X)_n=C_n(X)$ defined to be the free abelian group generated by the elements of $X_n$ with $d_i$ and $s_i$ in $C_*(X)$ taken to be the linear extensions of the face maps $d_i$ and $s_i$ of $X$. We can also form the total face map $d =\displaystyle \sum_{i=0}^n (-1)^id_i\colon C_n(X)\to C_{n-1}(X)$ and  then define the homology $H_*(X)$ as the homology of the chain complex $(C_*(X),d)$.

If $X=\ms S(Y)$, the singular set as defined in Example \ref{E: singular}, then we have $H_*(X)=H_*(Y)$, the singular homology of the space $Y$.
\end{example}

\begin{example}
Here's an example of a simplicial group that is important in the theory of homology of groups. Let $G$ be a group, and let $BG$ be the simplicial set defined as follows. Let $BG_n=G^{\times n}$, the product of $G$ with itself $n$ times. $G^{\times 0}$ is just the trivial group $\{e\}$. For an element $(g_1,\ldots,g_n)\in BG_n$, let 
\begin{align*}
d_0(g_1,\ldots, g_n)&= (g_2,\ldots, g_n)\\
d_i(g_1,\ldots, g_n)&= (g_1,\ldots, g_ig_{i+1},\ldots g_n)\quad\text{if $0<i<n$}\\
d_n(g_1,\ldots, g_n)&= (g_1,\ldots, g_{n-1})\\
s_i(g_1,\ldots, g_n)&= (g_1,\ldots,g_i,e,g_{i+1},\ldots, g_n).
\end{align*}
\end{example}

The reader can check that this defines a simplicial set. Unfortunately, it is not in general a simplicial group as $d_i$ will not necessarily be a homomorphism for $0<i<n$. But if $G$ is abelian, we will have a simplicial group. 
The realization of this simplicial set turns out to be  the classifying space of the group $G$, and so the homology $H_*(BG)$ coincides with group homology of the group $G$. For more on this simplicial set and its uses, the reader may consult \cite[Chapter 8]{WEIB}.

\section{Kan complexes}\label{S: Kan}

One of the goals of the development of simplicial sets (and other simplicial objects) was to find a  combinatorial way to study homotopy theory, just as simplicial homology theory allows us to derive invariants of simplicial complexes in a purely combinatorial manner (at least in principle).
Unfortunately, it turns out that not all simplicial sets are created equal as regards their usefulness toward this goal. The underlying reason turns out to be (once again, at least in principle) related to the reason that homotopy theorists prefer to work with CW complexes and not arbitrary topological spaces. Pairs of CW complexes satisfy the homotopy extension property, i.e. inclusions of subcomplexes are cofibrations (see, e.g., \cite{DK}). The condition we need to impose on simplicial sets to make them appropriate for the study of homotopy is similarly an extension condition. When seen through sufficiently advanced lenses, such as from the model category viewpoint presented in  \cite{GoeJar}, the extension condition on simplicial sets and the homotopy extension property in topology are essentially equivalent. 

As with much else in the theory of simplicial sets, the extension condition comes from a fairly straightforward idea that is often completely obfuscated in the formal definition. 

To explain the idea, we first need the following definition. 

\begin{figure}[!htp]
\begin{center}
\scalebox{.6}{\includegraphics{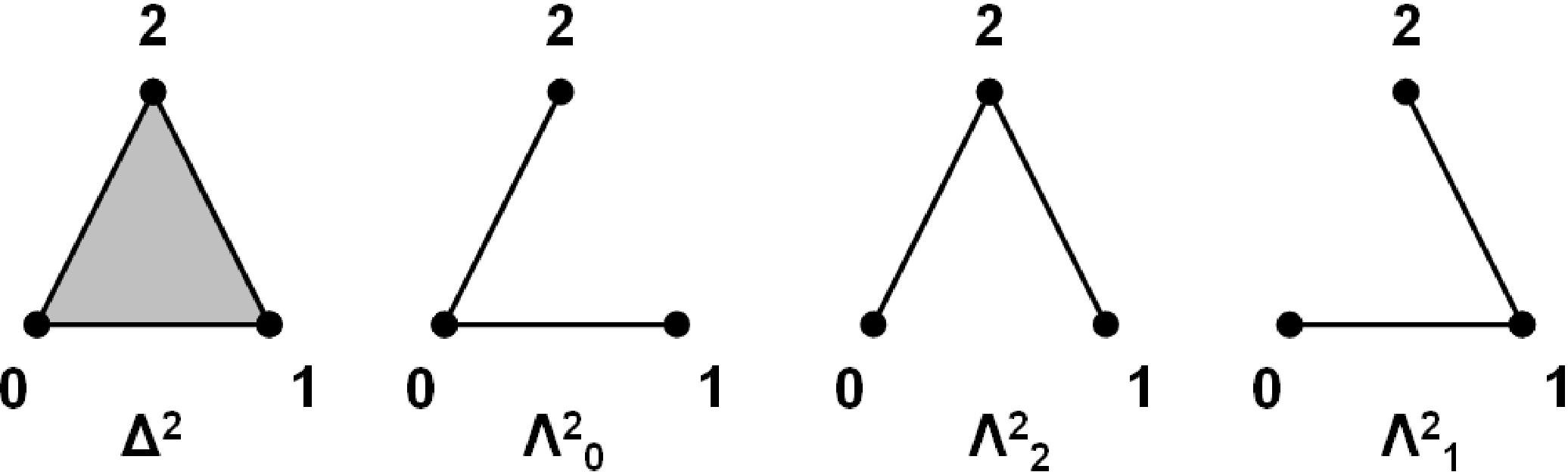}}
\end{center}
\caption{The three horns on $|\Delta^2|$}\label{F: fig18}
\end{figure}

\begin{definition}
As a simplicial complex, the $k$th horn $|\Lambda^n_k|$ on the $n$-simplex $|\Delta^n|$ is the subcomplex of $|\Delta^n|$ obtained by removing the interior of $|\Delta^n|$ and the interior of the face $d_k\Delta^n$. See Figure \ref{F: fig18}. We let $\Lambda^n_k$ refer to the associated simplicial set. This simplicial set consists of simplices $[i_0,\ldots,i_m]$ with $0\leq i_0\leq \cdots\leq i_m\leq n$ such that 1) not all numbers $0,\ldots, n$ are represented (this would be the top face or a degeneracy thereof) and 2) we never have all numbers except $k$ represented (this would be the missing $(n-1)$-face or a degeneracy thereof).  
\end{definition}

The extension condition, also known as the Kan condition (after Daniel Kan), says that whenever we see a horn on an $n$-simplex within a simplicial set, the rest of the simplex is there, too. Here's an elegant way to say this:

\begin{definition}
The simplicial object $X$  satisfies the \emph{extension condition} or \emph{Kan condition} if any morphism of simplicial sets $\Lambda^n_k\to X$ can be extended to a simplicial morphism $\Delta^n\to X$.

Such an $X$ is often called a \emph{Kan complex}\footnote{Note the risk of confusion here between simplicial sets and simplicial complexes. ``Kan complexes'' are simplicial \emph{sets}.} or, in more modern language, is referred to as being \emph{fibrant}.  
\end{definition}

We next present an equivalent formulation that is often used. This version has its advantages from the point of view of conciseness of combinatorial information, but it is much less conceptual.

\begin{definition}[Alternate version of the Kan condition]
The simplicial set $X$ satisfies the  \emph{Kan condition} if for any collection of $(n-1)$-simplices $x_0,\ldots, x_{k-1}, x_{k+1}, \ldots, x_{n}$ in $X$ such that $d_ix_j=d_{j-1}x_i$ for any $i<j$ with $i\neq k$ and $j\neq k$, there is an $n$-simplex $x$ in $X$ such that $d_ix=x_i$ for all $i\neq k$. 
\end{definition}

The condition on the simplices $x_i$ of the alternative definition glues them together to form the horn $\Lambda_k^n$, possibly with degenerate faces, within $X$, and the definition says that we can extend this horn to a (possibly degenerate) $n$-simplex in $X$.

\begin{example}\label{E: not Kan}
Not even the standard simplices $\Delta^n$, $n>0$, satisfy the Kan condition! Let $\Delta^1=[0,1]$ be the standard $1$-simplex, and consider the horn $\Lambda^2_0$, which consists of the edges $[0,2]$ and $[0,1]$ of $\Delta^2$, along with their degeneracies. Now consider the simplicial morphism that takes $[0,2]\in\Lambda^2_0$ to $[0,0]\in \Delta^1$ and $[0,1]\in\Lambda^2_0$ to $[0,1]\in\Delta^1$. There is a unique such simplicial map since we've specified what happens on all the nondegenerate simplices of $\Lambda^2_0$.  Notice that this is perfectly well-defined as a simplicial map since all functions on all simplices are order-preserving. However, this cannot be extended to a map $\Delta^2\to \Delta^1$ since we have already prescribed that $0\to 0$, $1\to 1$, and $2\to 0$, which is clearly not order-preserving on $\Delta^2$. 

For the same reason, no ordered simplicial complex $X$ (augmented to be a simplicial set) can ever satisfy the Kan condition unless $X$ is a discrete set of points! 
\end{example}

\begin{example}
It is easy to check that $\Delta^0$ does satisfy the Kan condition.
\end{example}

The following example is critical. 

\begin{figure}[!htp]
\begin{center}
\scalebox{.5}{\includegraphics{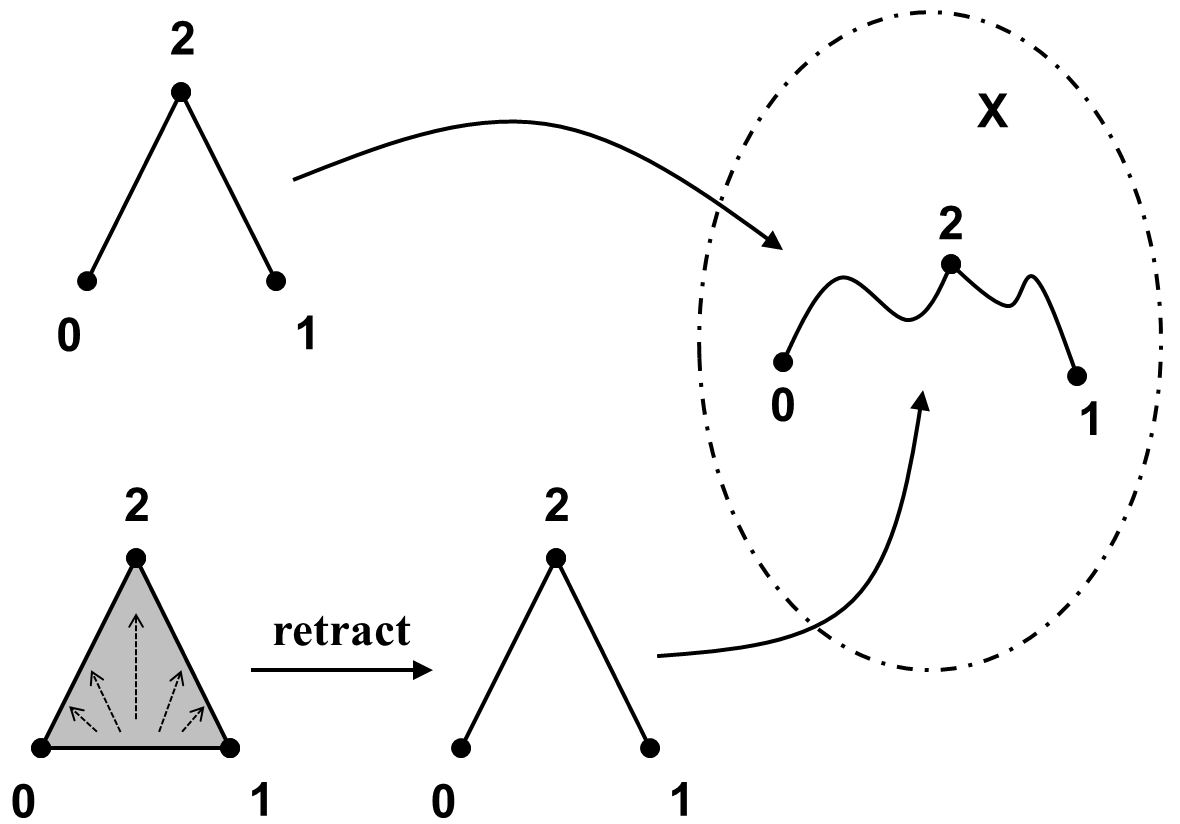}}
\end{center}
\caption{A demonstration that the singular set satisfies the Kan condition}\label{F: fig19}
\end{figure}

\begin{example}\label{E: simp Kan}
Given a topological space $Y$, the simplicial set $\ms S(Y)$ \emph{does} satisfy the Kan extension condition. It is actually fairly straightforward to see this. Consider any morphism of simplicial sets $f\colon\Lambda^n_k\to \ms S(Y)$. This is the same as specifying for each $n-1$ face, $d_i\Delta^n$, $i\neq k$, of $\Delta^n$ a singular simplex $\sigma_i\colon|\Delta^{n-1}|\to Y$. Every other simplex of $\Lambda^n_k$ is a face or a degeneracy of a face of one of these $(n-1)$-simplices, and so the rest of the map $f$ is determined by this data. Furthermore, the compatibility conditions coming from the simplicial set axioms ensure that the topological maps $\sigma_i$ piece together to yield, collectively, a continuous function $f\colon|\Lambda^n_k|\to Y$. It is now a simple matter to extend this function to all of $|\Delta^n|$: let $\pi\colon|\Delta^n|\to |\Lambda^n_k|$ be any continuous retraction (which certainly exists: $(|\Delta^n|,|\Lambda^k_n|)$ is homeomorphic to $(I^{n-1}\times I,I^{n-1}\times 0)$), and define $\sigma=f \pi\colon|\Delta^n|\to Y$. This is a singular $n$-simplex whose faces $d_if$, $i\neq k$, are precisely the singular simplices $\sigma_i$. Thus this determines the desired extension. See Figure \ref{F: fig19}.
\end{example}

\begin{example}
Any simplicial group is also, by neglect of structure, a simplicial set. All such simplicial sets arising from simplicial  groups satisfy the Kan condition. The proof is not difficult, but I don't know of a version that is particularly illuminating. Since we will not have much further use for this fact in these notes, we refer the reader to \cite[Theorem 2.2]{MOORE} for a proof. 
\end{example}

\section{Simplicial homotopy}\label{S: simp homotopy}

In this section we begin to look at the homotopy properties of simplicial sets. This is one of the key reasons that the theory of simplicial sets exists - to allow us to turn homotopy theoretic problems, at least in principle, into combinatorial problems by studying the homotopy groups of simplicial sets instead of those of topological spaces. In order to get started with simplicial homotopy, it is necessary to restrict attention to simplicial sets satisfying the Kan condition. This is not as large a handicap as it first appears, however, since we have already seen that, given a topological space $Y$, the singular set $\ms S(Y)$ satisfies the Kan condition, and eventually, we will see that $\ms S(Y)$  constitutes an appropriate combinatorial stand-in for $Y$. 

As we proceed, the reader should bear in mind the extent to which many of the ideas and definitions mirror those in topological homotopy theory. This may prove a comfort (or cause serious worry!) at those junctures where the mirror appears somewhat warped by the combinatorial complexity of the simplicial version.

We begin, naturally enough, with $\pi_0$, corresponding to the homotopy relationship between maps of points. This is a quite tractable warm-up for what is to follow.

\subsection{Paths and path components}

As in topology, when talking about homotopy, we will let $I$ stand for the simplicial set $\Delta^1=[0,1]$. As a simplicial set, $I$ has the nondegenerate $1$-simplex $[0,1]$, the nondegenerate $0$-simplices $[0]$ and $[1]$, and all other simplices are degenerate. Each simplex has the form $[0,\ldots,0,1,\ldots, 1]$ (possibly with no $0$s or no $1$s). 

\begin{definition}
A \emph{path} in a simplicial set $X$ is a simplicial morphism $p\colon I\to X$. Equivalently, a path in $X$ is a $1$-simplex $p\in X_1$. If $p$ is a path in $X$, $d_1p=p[0]$ is called the \emph{initial point} of the path and $d_0p=p[1]$ is called the \emph{final point} or \emph{terminal point}. 
\end{definition}

\begin{definition}
Two $0$-simplices $a$ and $b$ of the simplicial set $X$ are said to be \emph{in the same path component} of $X$ if there is a path $p$ with initial point $a$ and final point $b$.
\end{definition}

Already this definition appears slightly odd if you're used to working with simplicial complexes. In a connected simplicial complex, one might have to traverse several edges to link two vertices. Here we require it to be done all with one edge. Furthermore, we would expect ``being in the same path component'' to be an equivalence relation. This will not be the case in, say, an ordered simplicial complex in which we can have $a<b$ or $b<a$ but not both. What rescues this definition is the Kan condition.

\begin{theorem}
If $X$ is a Kan complex, then ``being in the same path component'' is an equivalence relation.
\end{theorem}
\begin{proof}
We will go through the proof in detail as it is very illuminating of how to think geometrically about simplicial sets.

\subparagraph{Reflexivity.} This one is easy: for any vertex $[a]$, $s_0[a]$ is a path from $a$ to $a$.

\subparagraph{Transitivity.} Consider $\Delta^2=[0,1,2]$. If $p_1$ is a path from $a$ to $b$ and $p_2$ is a path from $b$ to $c$, then let $f\colon\Lambda_1^2\to X$ take $[0,1]$ to $p_1$ and $[1,2]$ to $p_2$. See Figure \ref{F: fig20a}. The Kan condition lets us extend $f$ to $\bar f\colon\Delta^2 \to X$, and $\bar f[0,2]$ gives us a path from $a$ to $c$.  

\begin{figure}[!htp]
\begin{center}
\scalebox{.6}{\includegraphics{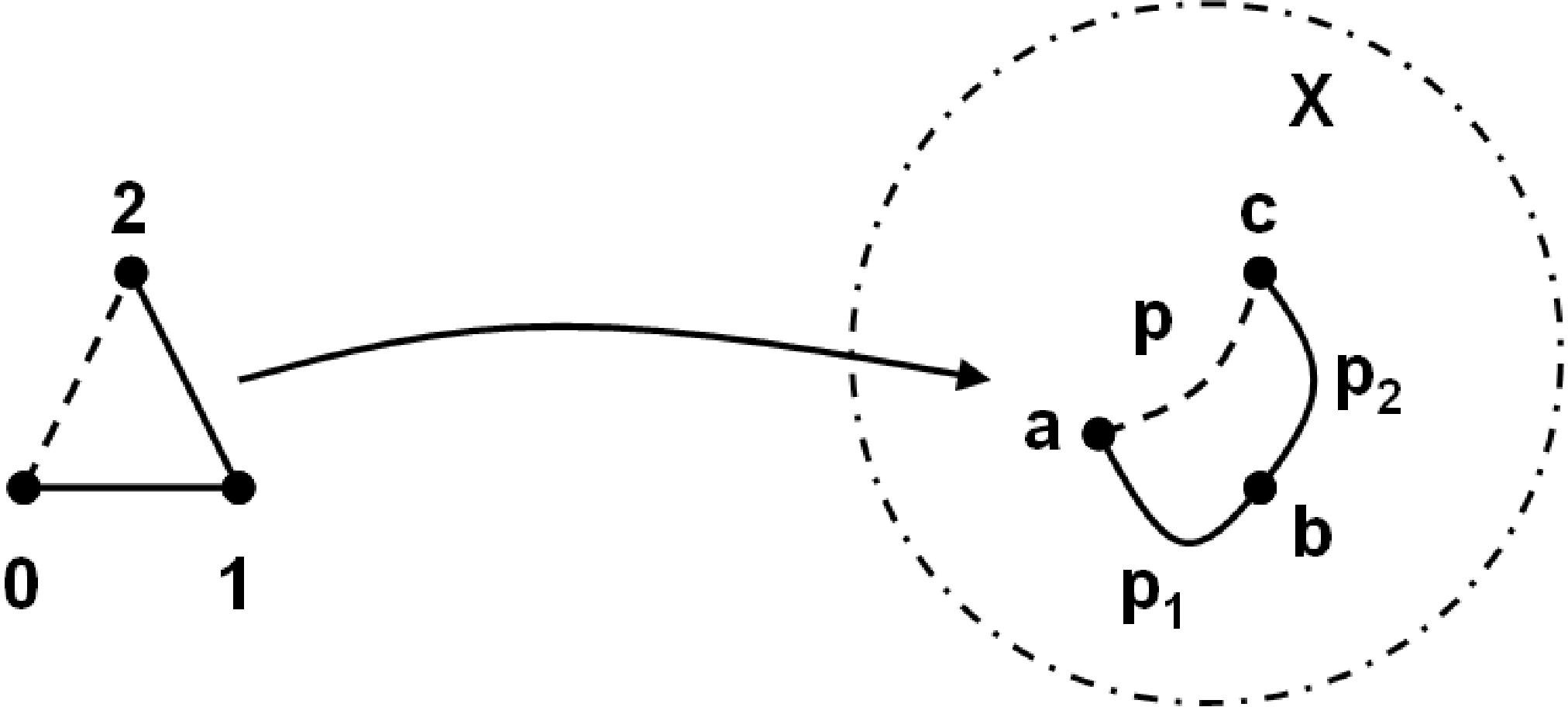}}
\end{center}
\caption{The transitivity relation on path connectedness via the Kan condition}\label{F: fig20a}
\end{figure}

\subparagraph{Symmetry.} This is only slightly more tricky than the transitivity condition.  See Figure \ref{F: fig20b}. Let $p$ be a path in $X$ from $a$ to $b$. We need a path the other way. Think of $p$ as the $[0,1]$ side of $\Delta^2$. Let the $[0,2]$ side of $\Delta^2$ represent $s_0[a]$, which must exist since $X$ is a simplicial set. Notice that $d_0s_0[a]=d_1s_0[a]=[a]$. At this point, we can label the three vertices $[0,1,2]$ of $\Delta^2$ as $[a,b,a]$, and we have a simplicial map on $\Lambda^2_0$ taking $[0,1]$ to $p$ and $[0,2]$ to $s_0[a]$. The Kan condition tells us that this map can be extended to all of $\Delta^2$ and $[1,2]$ gets taken to a path from $b$ to $a$. 
\begin{figure}[!htp]
\begin{center}
\scalebox{.6}{\includegraphics{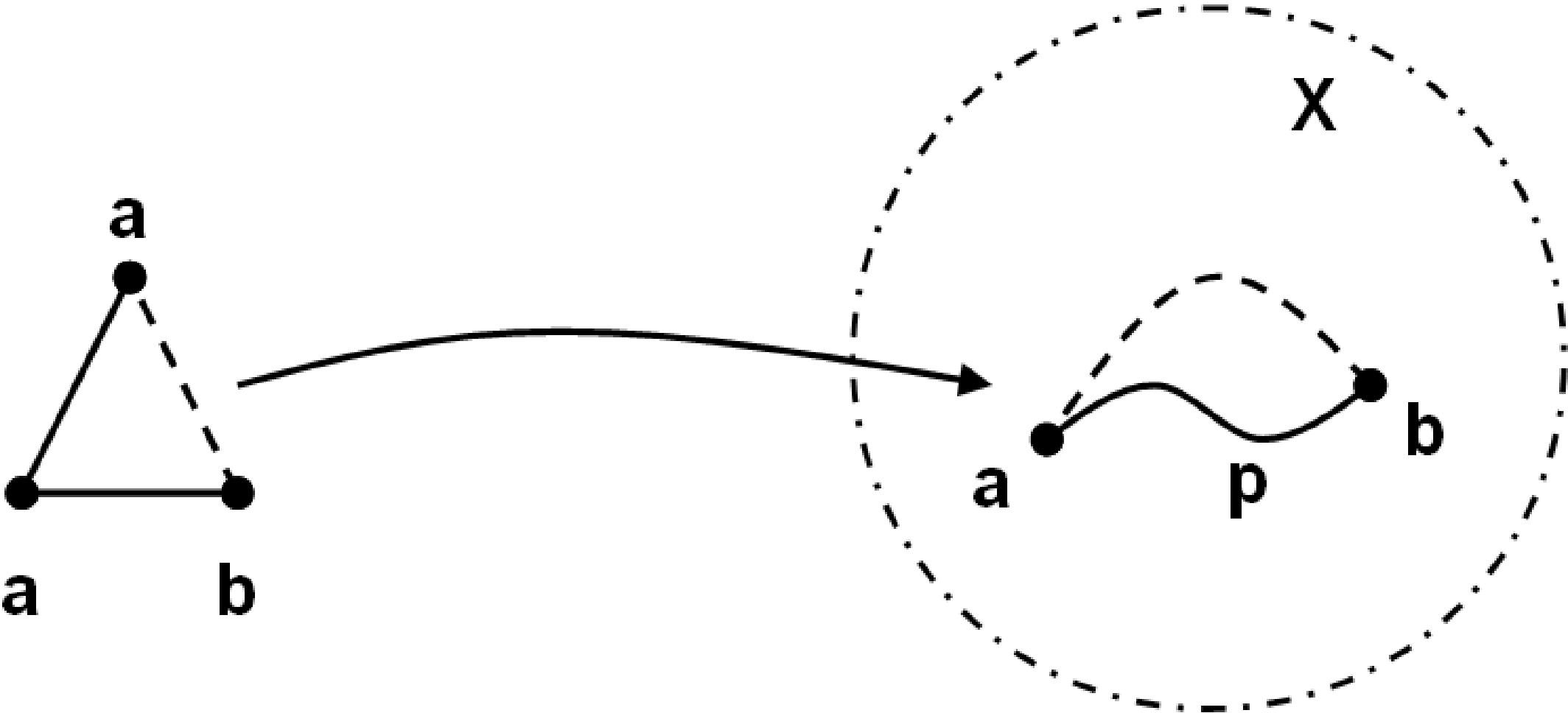}}
\end{center}
\caption{The symmetry relation on path connectedness}\label{F: fig20b}
\end{figure}
\end{proof}

Notice how important the Kan condition is here. 

Since we have demonstrated that being in the same path component is an equivalence relation, we have equivalence classes.

\begin{definition}
We denote the set of path components of $X$ (i.e. the equivalence classes of vertices of $X$ under the relation of being in the same path component) by  $\pi_0(X)$.
\end{definition}

So far, this is comfortingly familiar.

\subsection{Homotopies of maps}\label{S: homotopy}

There are at least two classical versions of the definition of simplicial homotopy, and at least two more modern versions for which we refer the reader to \cite{GoeJar}. Of the two classical versions, one has the expected form for a homotopy, $H\colon X\times I\to Y$. The other is more closely related to the homotopies we see in chain complexes $\widehat H\colon X_n\to Y_{n+1}$. We will look at both of these and see how they are related.

Perhaps the most natural definition of simplicial homotopy  looks something like  this:

\begin{definition}[Simplicial homotopy 1] 
Two simplicial maps $f,g\colon X\to Y$ are \emph{homotopic} if there is a simplicial map $H\colon X\times I\to Y$ such that $H|_{X\times 0}=g$ and $H|_{X\times 1}=f$ (i.e., if $g=H\circ i_0$ and $f=H\circ i_1$, where $i_0,i_1$ are the evident simplicial inclusion maps $i_0\colon X\times [0]\into X\times I$ and $i_1\colon X\times [1]\into X\times I$). 
\end{definition}

Unfortunately, here is the definition of simplicial homotopy one finds quite often in the literature:

\begin{definition}[Simplicial homotopy 2] \label{D: comb hom}
Two simplicial maps $f,g\colon X\to Y$ are \emph{homotopic} if for each $p$ there exist functions $h_i=h_i^p\colon X_p\to Y_{p+1}$ for each $i$, $0\leq i\leq p$, such that 
\begin{enumerate}
\item \begin{align*}d_0h_0&=f \\d_{p+1}h_p&=g \end{align*}
\item \begin{align*}
d_ih_j&=h_{j-1}d_i \qquad\text{if $i<j$}\\
d_{j+1}h_{j+1}&=d_{j+1}h_j\\
d_ih_j&=h_{j}d_{i-1} \qquad\text{if $i>j+1$}
\end{align*}
\item \begin{align*}
s_ih_j&=h_{j+1}s_i \qquad\text{if $i\leq j$}\\
s_ih_j&=h_{j}s_{i-1} \qquad\text{if $i>j$}.
\end{align*}
\end{enumerate}
\end{definition}

It will take some doing to see how these two definitions are related. This was one of the initial motivations for writing this exposition!

As usual, we will consider the universal example, $X=\Delta^p$, since once we understand how a homotopy works on a single simplex, we will also understand what happens along its faces and degeneracies, and everything else is determined by how the simplices are glued together. 

The key here is to recall Example \ref{E: prism} of Section \ref{S: product}, in which we showed how the prism $|\Delta^p\times I|$ is decomposed into simplices. In particular, it consists of $p+1$ nondegenerate $(p+1)$-simplices that we labeled $P_k\in (\Delta^p\times I)_{p+1}$, $0\leq k\leq p$. Suppose now that we have a homotopy $H\colon \Delta^p\times I\to Y$ from $f$ to $g$. Everything is determined by what $H$ does to the $P_k$, since every other nondegenerate simplex in $\Delta^p\times I$ is a face of one of these simplices. All other simplices in $\Delta^p\times I$ are degenerate, and so their images are determined by the images of the simplices of which they are degeneracies. 

How does this relate to the combinatorial Definition \ref{D: comb hom}? Let us denote the unique nondegenerate $p$-simplex of $\Delta_p$ by $E_p$. In this version, there are $p+1$ functions $h_i\colon E_p\to Y_{p+1}$. Each of the $p+1$ functions $h_i$ assigns to  $E_p$ a $(p+1)$-simplex of $Y$. Collectively, these give us the image of the prism over $E_p$ in $Y$. 

To see this, we use the notation $P_k=[0,\ldots,k,k',\ldots,p']$, $0\leq k\leq p$, for the $(p+1)$-simplices of the prism $\Delta^p\times I$ (see Example \ref{E: prism}). Given $H\colon \Delta^p\times I\to Y$, let  $h_k(E_p)$  correspond to the image $H(P_k)$ in $Y$. Now let's look at the conditions in  Definition \ref{D: comb hom} and see what they mean.

Starting with the first conditions, $d_0h_0(E_p)=d_0H(P_0)=H(d_0P_0)=H(d_0[0,0',\dots,p'])=H([0',\ldots, p'])=H\circ i_1(E_p)=f(E_p)$, using the first definition of homotopy. Similarly, $d_{p+1}h_p(E_p)=H([0,\ldots, p])=H\circ i_0(E_p)=g(E_p)$. So these conditions assure that the ends of the prism really are controlled by the maps $f$ and $g$.

The first and third equations of the second set of conditions mirror the observations made in Example \ref{E: prism} that most of the boundaries of the $(p+1)$-simplices of the prism $\Delta^p\times I$ are themselves  simplices of the prisms built on the boundary faces of $\Delta^p$. So these equations ensure that these faces of the $h_i(\Delta^p)$ are compatible with the actions of the homotopy maps $h_i^{j}$ of lower dimensions $j<p$ on the faces of $\Delta^p$. The second equation is the condition that the neighboring simplices $P_k$ and $P_{k+1}$ share a boundary. We invite the reader to glean these combinatorial details from the calculations in Example \ref{E: prism}.

The third set of equations can also be obtained in a fairly straightforward manner by working with the $P_k$. For example, we observe that for $i\leq j$, $s_iP_j=[0,\ldots,i,i,\ldots,j,j',\ldots,p']$, which is also the $(j+1)$st prism simplex on the degenerate simplex $[0,\ldots,i,i,\ldots,j-1,j,j+1,\ldots,p]$. In other words, the $i$th degeneracy of the $j$th prism $(p+1)$-simplex over $\Delta^p$ is the $(j+1)$st prism simplex over the $i$th degeneracy of $\Delta^p$. 
The geometric idea of these equations is a bit less obvious than in the preceding paragraphs, but really this is just the condition that the way the homotopy acts on degenerate simplices is determined by how it acts on the simplices of which they are degeneracies. 

Having described how the combinatorial conditions of the second definition correspond to the more geometric ideas of the first definition, we now leave it to the interested reader to verify the complete equivalence of the two definitions, in particular to verify that the data given by all the $h_i^j$ is enough to reconstruct $H$.

\smallskip

We would like  homotopy to be an equivalence relation, but this will not hold in general. For example, in our discussion of path connectedness, which we see in the current language corresponds directly to homotopies of maps $\Delta^0\to X$, we saw that path connectedness is not always an equivalence relation. However, the discussion of path connectedness might lead one to suspect that we will be safe in the world of Kan complexes, and this is so.

\begin{theorem}
Homotopy of maps $X\to Y$ is an equivalence relation if $Y$ is a Kan complex. If $f$ and $g$ are homotopic, we denote that by $f\sim g$. 
\end{theorem}

We invite the reader to prove this by extending the argument given for path connectedness. An indirect proof involving the  ``function complexes'' $\Homs(X,Y)$ can be found in \cite{MAY67}; see in particular the discussion on page 17 and Corollary 6.11 of \cite{MAY67}. 

It is also fairly straightforward to verify other expected elementary facts about homotopy; for instance if $f\sim f'$, then $fg\sim f'g$ and $gf\sim gf'$. Also, homotopic maps induce the same homomorphisms on homology groups (see Section \ref{S: simplicial sets} - this follows as for the usual proof in singular homology theory by using the triangulation of the homotopy prism; see, e.g. \cite{MK}). See \cite[Section I.5]{MAY67} for proofs of these facts.

\begin{remark}
Notice that homotopies $H\colon X\times I\to Y$ correspond to elements of the set $\Homs(X,Y)_1$ as defined in Section \ref{S: Hom}, just as elements of $\Homs(X,Y)_0$ correspond to simplicial maps. If $Y$ is a Kan complex, $\Homs(X,Y)$ will be a Kan complex as well (see \cite[Theorem 1.6.9]{MAY67}), so 
in this case it makes sense to observe that two simplicial maps $f,g\in \Homs(X,Y)_0$ will be  homotopic if and only if they lie in the same path component of $\Homs(X,Y)$. In other words, $f,g:X\to Y$ are homotopic if and only if they represent the same element of $\pi_0(\Homs(X,Y))$.
\end{remark}

\subsection{Relative homotopy}\label{S: relative homotopy}

The notions of subcomplexes and relative homotopy offer no surprises, but we record the definitions for clarity.

\begin{definition}
If $X$ is a simplicial set, then $A$ is a simplicial subset of $X$, denoted $A<X$, if $A$ itself is a simplicial set such that $A_n\subset X_n$ for all $n$ and the face and degeneracy maps of $A$ agree with those from $X$. A pair of simplicial sets is often denoted by $(X,A)$. 
Simplicial maps of pairs $(X,A)\to (Y,B)$ are simplicial maps $X\to Y$  such that the image of $A$ is contained in $B$.
\end{definition}

\begin{definition}
If $(X,A)$ are a simplicial set and simplicial subset and both $X$ and $A$ satisfy the Kan condition, then $(X,A)$ is called a Kan pair.
\end{definition}

\begin{example}
An important example of a simplicial subset of a simplicial set $X$ is a \emph{basepoint} for $X$, consisting of an element of $X_0$ and all of its degeneracies. We will denote basepoints by $*$. Notice that $*$ is isomorphic as a simplicial set to $\Delta^0$ and can be considered as an image $\Delta^0\to X$ of a simplicial map. Since $\Delta^0$ is a Kan complex, $(X,*)$ will be a Kan pair if $X$ is Kan. 
\end{example}

\begin{example}
Note that  a subcomplex of a Kan complex need not be Kan. For instance, we know from Example \ref{E: not Kan} that the simplex $\Delta^1$ is not a Kan complex. We also know that the singular set $\ms S(|\Delta^1|)$ on the space $|\Delta^1|$ is a Kan complex, by Example \ref{E: simp Kan}. But the former is a subcomplex of the latter, realized by the singular simplices that represent $|\Delta^1|$ as a simplicial complex.  Namely, $\Delta^1$ corresponds to the subcomplex of $\ms S(|\Delta^1|)$ generated by the singular $0$-simplices $\sigma_0\colon |\Delta^0|\to [0]$ and $\sigma_1\colon |\Delta^0|\to [1]$, by the singular $1$-simplex $\text{id}\colon |\Delta^1|\to|\Delta^1|$, and by their degeneracies.
\end{example}

\begin{definition}
Given a simplicial pair $(X,A)$, a homotopy $H\colon X\times I\to Y$ is a homotopy \emph{rel $A$} if the restriction of $H$ to $A\times I$ can be factored as $H|_{A\times I}=g\pi_1\colon A\times I\to Y$, where $g$ is a simplicial map $A\to Y$ and $\pi_1$ is the projection $A\times I\to A$ (i.e., if the homotopy $H$ is constant on the simplicial subset $A$). If $Y$ is Kan, then homotopy rel $A$ is an equivalence relation. 
\end{definition}

While considering simplicial pairs, there is another crucial theorem we should mention: the homotopy extension theorem for simplicial maps to Kan complexes.

\begin{theorem}[Homotopy extension theorem]\label{T: HET}
Let $(X,A)$ be a pair of simplicial sets and $Y$  a Kan complex. Suppose there is a simplicial map $f\colon X\to Y$ and a simplicial homotopy $H\colon  A\times I\to Y$ such that $H|_{A\times 0}=f|A$. Then there exists an extension $F\colon X\times I\to Y$ such that $F|_{A\times I}=H$ and $F|_{X\times 0}=f$. 
\end{theorem}

Unfortunately, the proofs I know would all take us too far afield, so we refer the reader to \cite[Chapter 1, Appendix A]{MOORE} for a combinatorial treatment or \cite[Section I.4]{GoeJar} for a more modern treatment.

\section{$\mathbf{\pi_n(X, *)}$}\label{S: pin}

In this section, we will discuss the homotopy groups of Kan complexes. This section is a bit more technical than the preceding ones, as we here need some theorems and not just definitions. This section should serve as good technical practice for the reader preparing to go on to read  further material on simplicial objects. 

Given a  Kan complex  with basepoint $(X,*)$, there are at least four ways to define $\pi_n(X,*)$:

\begin{enumerate} 

\item One can define these groups directly as homotopy classes of maps of (simplicial) spheres to $X$.

\item There is a purely combinatorial definition.

\item As in algebraic topology, one can first define appropriate iterated simplicial loop spaces $\Omega^n(X)$ and define $\pi_n(X)=\pi_0(\Omega^n(X))$. 

\item As a more topological alternative, one could try the topological homotopy groups of the realization of $X$, i.e. $\pi_n(|X|,|*|)$.
\end{enumerate}

We will focus on the relationship between the first two of these, referring the interested reader to \cite{MOORE} for the third approach. For hints at the relevance of the fourth approach, see Theorem \ref{T: equiv}, below, as well as the discussion in Section \ref{S: concluding} in general.

The definition of $\pi_n(X,*)$ in terms of spheres is straightforward once we decide what a sphere is. Example \ref{E: sphere} teaches us that there is more than one reasonable definition, or at least more than one simplicial set whose realization is a sphere. In fact, we will see that both versions treated in that example are acceptable. 

\begin{definition}[First definition of $\pi_n$]\label{D: pin1}
Given a Kan complex  with basepoint $(X,*)$, define $\pi_n(X,*)$, $n>0$, to be the set of homotopy equivalence classes of maps $(\bd \Delta^{n+1},*)\to (X,*)$. Here, we take for the basepoint of $\bd \Delta^{n+1}$ the simplicial subset of $\Delta^{n+1}$ generated by the vertex $[0]$, and all homotopies are relative to the basepoint.
\end{definition}

It might be a good exercise even at this point for the reader to prove that if a map $(\bd \Delta^{n+1},*)\to (X,*)$ extends to a map $(\Delta^{n+1},*)\to (X,*)$ then it is homotopic to the constant map sending all of $\bd \Delta^{n+1}$ to the basepoint $*$.

The requirement in the definition that $X$ be Kan is necessary for homotopy to be an equivalence relation. Of course, we want $\pi_n(X,*)$ to be a group, but this will have to wait a moment. Let's first work toward the more combinatorial definition. This takes a little bit of preliminary effort.

\begin{definition}\label{D: homotopy of simplices}
We say that two $n$-simplices $x,x'\in X_n$ are \emph{homotopic} if 
\begin{enumerate}
\item $d_ix=d_ix'$ for $0\leq i\leq n$, and 
\item there exists a simplex $y\in X_{n+1}$ such that 
\begin{enumerate}
\item $d_ny=x$, 
\item $d_{n+1}y=x'$, and
\item $d_iy=s_{n-1}d_ix=s_{n-1}d_ix'$, $0\leq i\leq n-1$. 
\end{enumerate}
\end{enumerate}
\end{definition}

\begin{figure}[!htp]
\begin{center}
\scalebox{.6}{\includegraphics{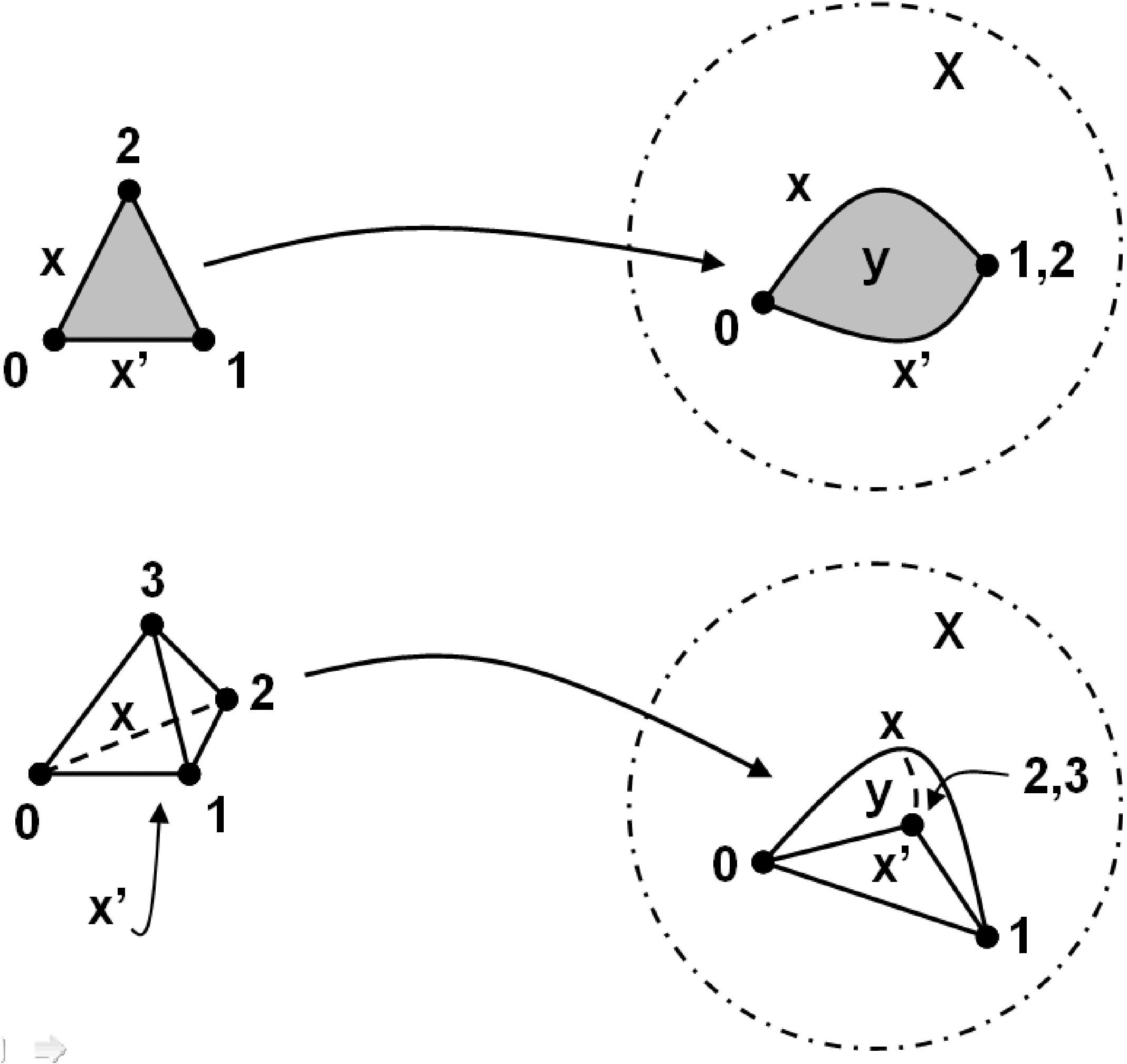}}
\end{center}
\caption{Above: a homotopy of $1$-simplices. Below: a homotopy of $2$-simplices. The picture in the bottom right depicts two $2$-simplices glued together along their boundaries.}\label{F: fig23}
\end{figure}

The idea  here  is that $x$ and $x'$ have the same boundaries and that $y$ provides the homotopy between them, rel boundary, by letting $x$ and $x'$ be two of the faces of $y$, while the rest of the faces of $y$ degenerate to the boundaries of $x$ and $x'$.  In other words, the simplex $y$ can be imagined as having $x$ on one face and $x'$ on another face so that $x$ and $x'$ share an edge in common and then their other corresponding edges are collapsed together as degeneracies; see Figure \ref{F: fig23}.

It can be shown directly that homotopy of simplices is an equivalence relation if $X$ is  a Kan complex. The argument is a generalization of the one showing that path connectedness is an equivalence relation. Again the idea is to arrange a simplex so that the pieces we know fall on certain faces of horns and the pieces we'd like to show exist  fall on the missing faces. Then these relations must exist due to the Kan extension condition.
We refer the interested reader to \cite[Section I.3]{MAY67}.

\begin{definition}[Second definition of $\pi_n$]\label{D: pin2}
Given a Kan complex  with basepoint $(X,*)$, we can also define $\pi_n(X,*)$, $n>0$, as the set of equivalence classes of $n$-simplices $x\in X_n$ with $d_ix\in *$ for all $i$, $0\leq i\leq n$, up to homotopy of simplices. 
\end{definition}

This version of the homotopy groups corresponds more closely to our second version of the sphere in Example \ref{E: sphere}. Recall that, as a simplicial set, this version of the sphere $S^n$ had only two nondegenerate simplices: one in dimension $n$ and one in dimension $0$. An $n$-simplex of $X$ all of whose faces live in $*$ can be thought of as the image of that simplicial version of $S^n$ in $X$. Thus this definition of $\pi_n(X,*)$ also makes some geometric sense. However, there are some obvious questions, such as:  Why do the first and second definitions of $\pi_n$ agree? And where is the group structure we expect?

To answer the first question, we need a series of lemmas:

\begin{lemma}\label{L: alt homotopy}
If $X$ is Kan and $d_ix=d_ix'$ for all $i$, we obtain the same equivalence relation as in Definition \ref{D: homotopy of simplices} if we instead require that $d_ry=x$, $d_{r+1}y=x'$ for some $0\leq r\leq n$,  and $d_iy=d_is_rx=d_is_rx'$ for  $i\neq r, r+1$.
\end{lemma}
\begin{proof}
We refer the reader to \cite{MAY67} for the full proof, which is contained within Lemma 5.5 there. The idea is to show that the case of the definition using $r, r+1$ is equivalent to the version with $r+1, r+2$ for each relevant $r$. This is done using an extension argument by which one creates an $(n+2)$-simplex which has the two desired homotopies on two of the sides. We illustrate a low-dimensional case in Figure \ref{F: fig25}: Suppose that $x, x'$ are $1$-simplices and that we have a $y$ with $d_0y=x$, $d_1y=x'$. We want to find a $z$ with $d_1z=x$, $ d_2z=x'$. We form the horn $\Lambda^2_0$, shown flattened on the right of Figure \ref{F: fig25}. We embed $y$ as $[0,1,3]$ (note that this maintains its orientation simplicially despite the oddities of the drawing). We let the other sides of the horn be appropriate degeneracies of $x'$. Notice that there is no trouble embedding this horn in $X$  extending $y\in X$. Now the Kan condition assures us that we can extend this embedding to all of $\Delta^3$, including the remaining face $[1,2,3]$. We can check that this last face can be taken as the desired $z$ (be careful to notice that $d_1[1,2,3]=[1,3]$ and $d_2[1,2,3]=[1,2]$). 

The idea in higher dimensions is precisely the same; the extra  faces of the horn that exist in higher dimensions contain other degeneracies of faces of $x$ - see \cite[Lemma 5.5]{MAY67}, \cite[Proposition 1.19]{Cu71}.
\begin{figure}[!htp]
\begin{center}
\scalebox{.65}{\includegraphics{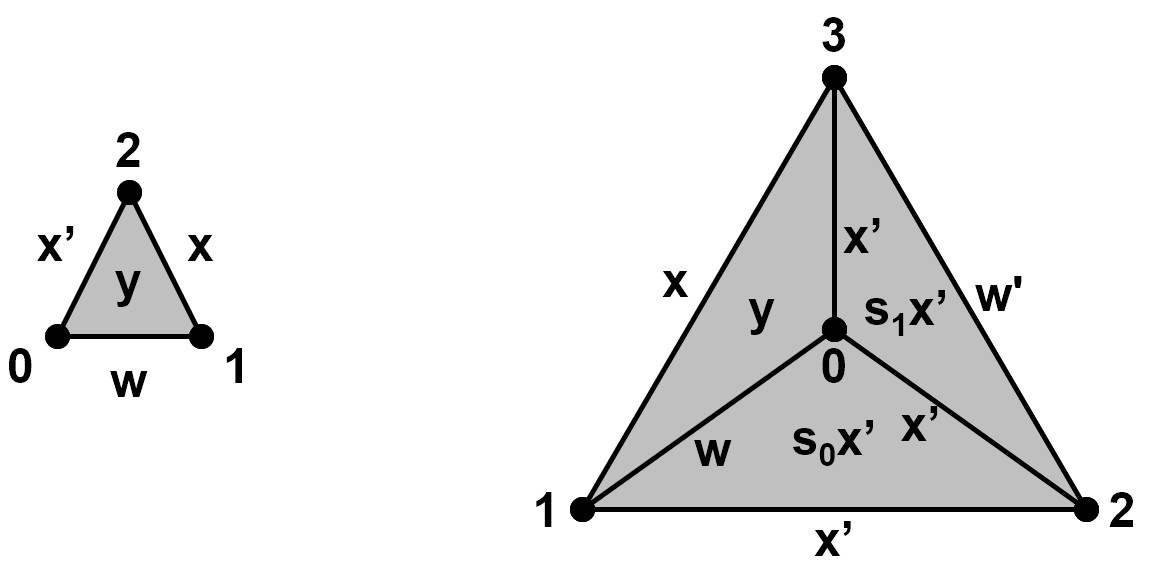}}
\end{center}
\caption{Shifting indices in the homotopy relation. Here $w$ stands for  $d_2s_0x=d_2s_0x'$, which is a degenerate $1$-simplex, both vertices being the first vertex of $x$, which is also the first vertex of $x'$. The edge labeled $w'$ is $d_0 s_1 x' = s_0 d_0 x' = s_0d_0x$.}\label{F: fig25}
\end{figure}
\end{proof}

\begin{lemma}\label{L: equiv. of homotopies}
If $X$ is Kan, two $n$-simplices $x,x'\in X$ with $d_ix, d_ix'\in *$ for all $i$, $0\leq i\leq n$, are homotopic in the sense of Definition \ref{D: homotopy of simplices} if and only if the  maps $f\colon \Delta^n \to X$ and $f'\colon \Delta^n\to X$ that represent $x$ and $x'$ are homotopic rel $\bd \Delta^n$ as maps. 
\end{lemma}
\begin{proof}
Of course to say that $f$ represents $x$ means that $f$ takes the nondegenerate $n$-simplex $E_n$ of $\Delta^n$ to $x\in X$. 

One direction of the argument is fairly straightforward. In order to show that $f$ and $f'$ are homotopic, it suffices to find a chain of $n+1$ simplices of dimension $n+1$, representing the images of nondegenerate simplices of the prism $\Delta^n\times I$, such that the ``top'' and ``bottom'' faces of the first and last simplex represents $x$ and $x'$. But if we know that $x$ and $x'$ are homotopic as simplices, we know there is one $(n+1)$-simplex $y$ connecting them with, say, $d_ny=x$, $d_{n+1}y=x'$, and $d_iy\in *$ for all other $i$. So now we just let $y$ be the $(n+1)$st simplex $h_n(\Delta^n)$, and we let $h_i(\Delta^n)=s_ix$ for $0\leq i\leq n$. In other words, we let the last nondegenerate simplex in $\Delta^n\times I$ do all the work of the homotopy, and we just collapse all the rest into the face representing $x$. See Figure \ref{F: fig24}. 

\begin{figure}[!htp]
\begin{center}
\scalebox{.6}{\includegraphics{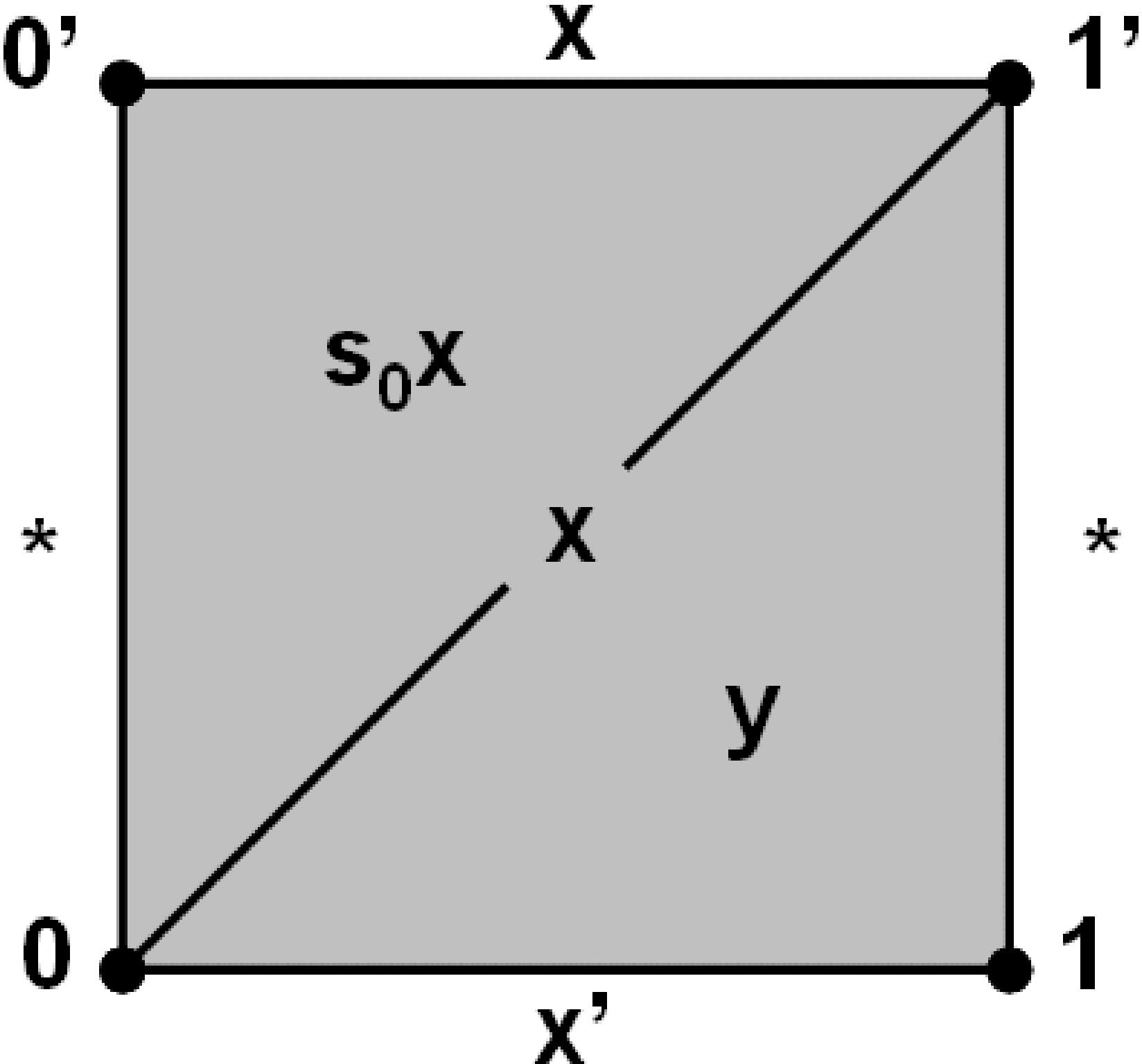}}
\end{center}
\caption{We label the vertices with the ``prism notation'' of Example \ref{E: prism}. The bottom simplex $y$ is a homotopy of the simplices $x$ and $x'$. Adjoining the degenerate simplex $s_0x$ shows how to obtain a model prism for  the homotopy from $x$ to $x'$.}\label{F: fig24}
\end{figure}

In the other direction, suppose we have an actual homotopy rel $*$ from $x$ to $x'$ thought of as inclusion maps. By definition, this gives us a prism $\Delta^n\times \Delta^1\in X$ whose top is $x$ and whose bottom is $x'$. We know from the discussion in Example \ref{E: prism} that each of the nondegenerate $(n+1)$-simplices of the prism has two $n$-faces that are not in $\bd \Delta^n\times \Delta^1$, and the rest are in $\bd \Delta^n\times \Delta^1$, all of which go to $*$ in $X$. Furthermore, it is not hard to check that the two  $n$-faces not in $\bd \Delta^n\times \Delta^1$ are consecutive faces. In particular, using the notation of Example \ref{E: prism}, these faces are $d_kP_k$ and $d_{k+1}P_k$. Thus by Lemma \ref{L: alt homotopy}, each $P_k$ is a homotopy between these two faces. Since the top and bottom faces of the prism are $x$ and $x'$, we obtain a simplicial homotopy between $x$ and $x'$ using the transitivity of simplicial homotopy.
\end{proof}

Thus, to show that our two definitions of $\pi_n(X,*)$ agree, it is only necessary to prove the following lemma, which is familiar in the context of algebraic topology. The proof is somewhat long, but we provide most of the details, as it is difficult to find a direct proof in the standard expositions.

\begin{lemma}\label{L: lemma}
If $X$ is a Kan complex, there is a bijection between homotopy classes of maps $f\colon (\bd \Delta^{n+1},*)\to (X,*)$ and homotopy classes of maps $g\colon  (\Delta^n, \bd \Delta^n)\to (X,*)$.
\end{lemma}
\begin{proof}
Given $g\colon (\Delta^n, \bd \Delta^n)\to (X,*)$, it is easy to construct an associated $f$ by identifying $\Delta^n$ with $d_0\Delta^{n+1}$. Then we let $f\colon  (\bd \Delta^{n+1},*)\to (X,*)$ be defined so that $f$ is given by $g$ on $d_0\Delta^{n+1}$ and by the unique map to $*$ on each $d_i\Delta^{n+1}$, $i>0$. It is also straightforward to see that any homotopy of $g$ rel $\bd \Delta^n$ determines a homotopy of $f$ rel $*$. 

Conversely, suppose we are given $f\colon (\bd \Delta^{n+1},[0])\to (X,*)$. We show that $f$ is homotopic to a function $\td f$ that takes $\Lambda_0^{n+1}$ to $*$. Then we can let $g$ be $\td f|_{d_0\Delta^{n+1}}$.

We first observe, as noted in the proof of Lemma \ref{L: equiv. of homotopies}, that to construct a homotopy between two  $k$-simplices $x$ and $x'$ in $X$, it suffices to find a simplex $y$ in $X$ with $d_ky=x$, $d_{k+1}y=x'$ since this can be considered one of the  blocks of a prism, and the rest of the prism can be filled up with degeneracies of $x$ or $x'$. 

Keeping this in mind, we proceed by induction with the following induction step: Suppose $f_{k-1}\colon \bd \Delta^{n+1}\to X$ is such that $f([0])\in *$ and $f(z)\in *$ for all simplices $z\in\bd \Delta^{n+1}$ of dimension $\leq k-1$ such that $[0]$ is a vertex of $z$, then there is a homotopy from $f_{k-1}$ to an $f_k$ that takes all simplices up to dimension $k$ having $[0]$ as a vertex to $*$. Furthermore, the homotopy can be performed rel the faces of dimension $\leq k-1$ having $[0]$ as a vertex.

Clearly we can take $f_0=f$. So suppose we have constructed $f_{k-1}$ for $k\geq 1$. We need only find the desired homotopy on the $k$-simplices of $\Delta^{n+1}$ that have $[0]$ as a vertex, and then we can apply the homotopy extension theorem, Theorem \ref{T: HET}. So let $z$ be a $k$-simplex of $\Delta^{n+1}$ with $0$ as a vertex. We know that $f_{k-1}(d_iz)\in *$ for $i\neq 0$. Now, consider the horn $\Lambda^{k+1}_0$, and note that we can map $\Lambda^{k+1}_0$ into $X$ such that the $k$-face corresponding to $d_{k+1}\Delta^{k+1}$ is $f_{k-1}(z)$ and such that all other $k$-faces are taken into $*$. Notice that this is possible precisely because $f_{k-1}(d_iz)\in *$ for $i\neq 0$. Now since $X$ is a Kan complex, we can extend this horn to a $(k+1)$-simplex $y$ in $X$ such that $d_{k+1}y=f_{k-1}(z)$ and $d_ky\in *$.  As noted, this is enough to construct a homotopy on $z$ from $f_{k-1}(z)$ to the unique map of $z$ into $*$. In addition, this is a homotopy rel those faces of $z$ that have $[0]$ as a simplex.
Notice also that it is possible to find such homotopies for all such $z$ independently and compatibly. In this way, we get a homotopy on the $k$-simplices of $\Delta^{n+1}$ having $[0]$ as a vertex from $f_{k-1}$ to the map to $*$.  Extending this homotopy by the homotopy extension theorem yields the desired homotopy to $f_k$. 

Continuing inductively, we obtain a map $f_{n+1}\colon \bd \Delta^{n+1}\to X$ homotopic to $f$ such that $\Lambda^{n+1}_0$ is taken to $*$. Now we can define $g$ to be the restriction of $f_{n+1}$ to $d_0\Delta^{n+1}$.

If $f,f'\colon (\bd \Delta^{n+1},[0])\to (X,*)$ are homotopic rel $[0]$, then we can show that the resulting $g$ and $g'$ are homotopic  by building a homotopy from the homotopy $H\colon \bd \Delta^{n+1}\times I\to X$ from $f$ to $f'$ to a homotopy $H_{n+1}\colon \bd \Delta^{n+1}\times I\to X$ such that  $H_{n+1}(\Lambda^{n+1}_0\times I)\in *$ and that extends the homotopies built over $f$ and $f'$ as in the preceding paragraphs. Then $H_{n+1}|_{d_0\Delta^{n+1}\times I}$ will be a homotopy from $g$ to $g'$. We leave the details to the reader. 
\end{proof}

Lemmas \ref{L: equiv. of homotopies} and \ref{L: lemma} together prove the following.

\begin{proposition}
If $X$ is a Kan complex, the definitions of $\pi_n(X,*)$ in Definitions \ref{D: pin1} and \ref{D: pin2} agree.
\end{proposition}

\paragraph{The group structure.}

One benefit of the version of $\pi_n(X,*)$ given in Definition \ref{D: pin2}, compared to the perhaps more geometrically transparent Definition \ref{D: pin1}, is the ease of  proving that $\pi_n(X,*)$ is a group and of describing the group operation.

\begin{figure}[!htp]
\begin{center}
\scalebox{.6}{\includegraphics{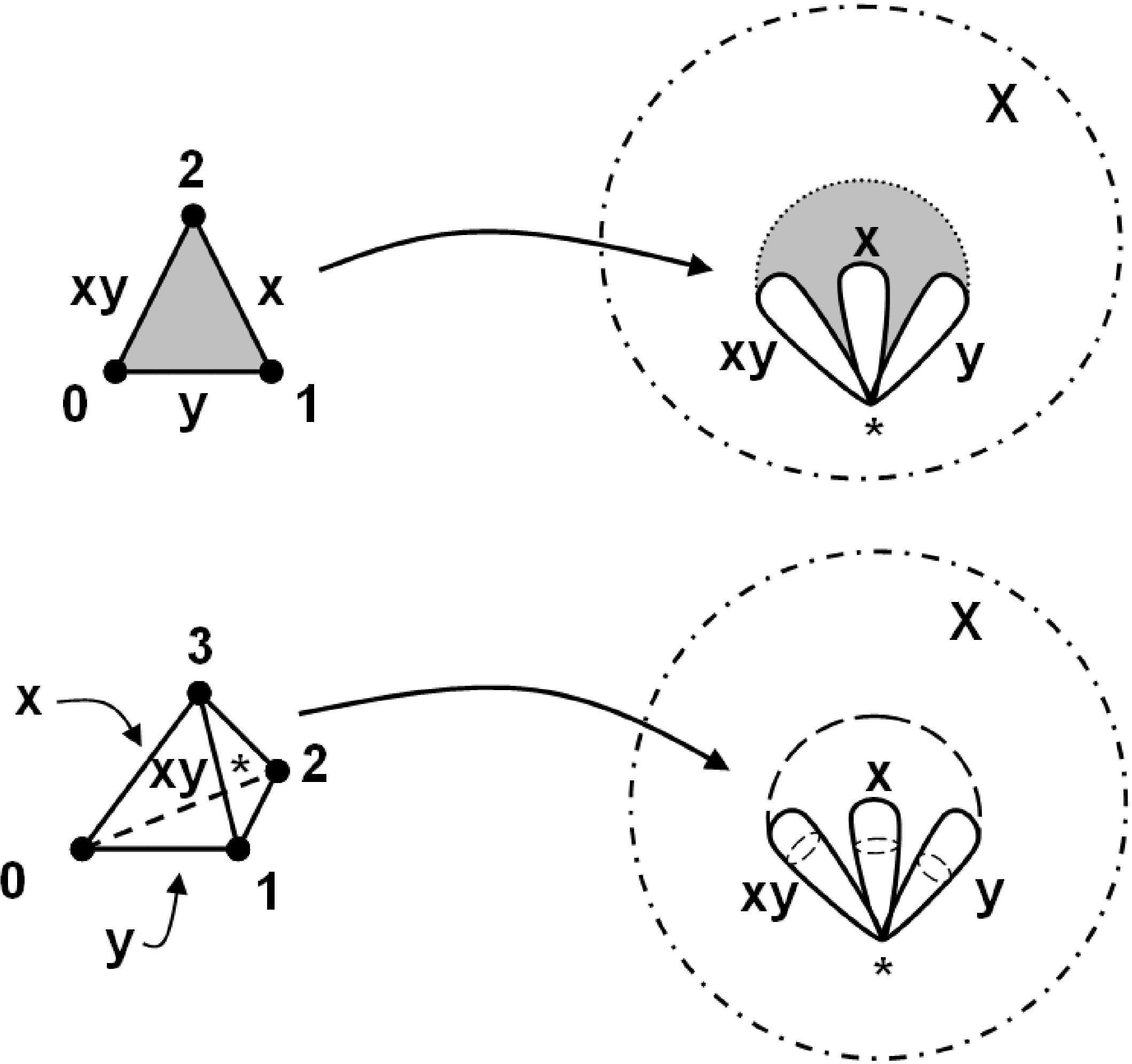}}
\end{center}
\caption{The product of $x$ and $y$ in $\pi_1(X,*)$ (above) or $\pi_2(X,*)$ (below).}\label{F: fig26}
\end{figure}

\begin{definition}\label{D: pi product}
Let $x,y$ be two $n$-simplices, $n\geq 1$, in the Kan complex $X$ such that $d_ix=d_iy\in *$ for all $i$. Let $\Lambda^{n+1}_n$ be the horn of $\Delta^{n+1}$ in $X$ such that the face corresponding to $d_{n+1}\Delta^{n+1}$ is $y$, the face corresponding to $d_{n-1}\Delta^{n+1}$ is $x$, and the faces corresponding to all other sides of the horn are in $*$. Let $z$ be an extension of the horn to $\Delta^{n+1}$ as guaranteed by the Kan condition. Then define $xy$ as the homotopy class of $d_nz$ in $\pi_n(X,*)$. See Figure \ref{F: fig26}.

It will be useful to say that $z$ \emph{realizes} the product $xy$. 
\end{definition}

It can be shown that the definition is independent of the choices made:

\begin{proposition}
The product of Definition \ref{D: pi product} yields a well-defined function $\pi_n(X,*)\times \pi_n(X,*)\to \pi_n(X,*)$.
\end{proposition}
\begin{proof}
The proof is by various applications of the Kan extension condition. See \cite[Lemma 4.2]{MAY67}. This would also be a good exercise for the reader.
\end{proof}

The idea of the product on the simplicial $\pi_n(X,*)$ is not far from that for the product in the topological homotopy groups. First, suppose one has a map of the $(n+1)$-ball $D^{n+1}$ to a topological space $X$ such that the equator of the boundary sphere $S^n$ is mapped to the basepoint of $X$. Then the restrictions of the map to the  upper and lower hemispheres of $S^{n}$ determine elements of $\pi_n(X,*)$, and the map of all of $D^{n+1}$ determines a homotopy between them. Secondly, recall that, roughly speaking, the product of two elements $x,y$ in the topological $\pi_n(X,*)$ can be represented by a map of a sphere that agrees with $x$ and $y$ on two disjoint disks in $S^n$ and takes the rest of $S^n$ to the basepoint.

Definition \ref{D: pi product} puts these ideas together. 
In the simplicial world, we can think of $d_n\Delta^{n+1}$ as being one hemisphere of $\bd \Delta^{n+1}$ and the rest of $\bd \Delta^{n+1}$ as the other hemisphere. Then in  Definition \ref{D: pi product}, the $(n+1)$-simplex $z$ can be thought of as providing a homotopy between $d_nz$ and what is happening on the rest of $\bd z$ (notice that, indeed, $\bd d_nz\in *$). But the rest of $\bd z$ contains $x$ and $y$ on two separate faces and everything else goes to $*$, just as for the topological product.

Of course we expect $\pi_n(X,*)$ to be a group if $n>0$.

\begin{theorem}
With the product of Definition \ref{D: pi product}, $\pi_n(X,*)$ is a group. 
\end{theorem}
\begin{proof}
The constructions corresponding to the necessary axioms are pictured in Figure \ref{F: fig27}.

\begin{figure}[!htp]
\begin{center}
\scalebox{.6}{\includegraphics{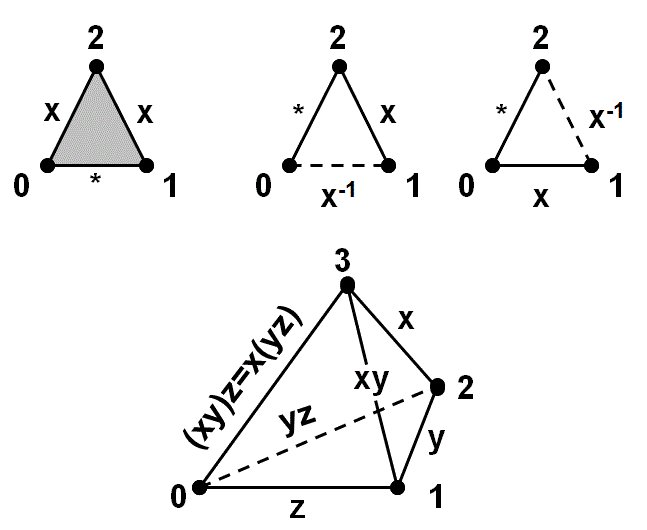}}
\end{center}
\caption{Above left: the identity $x*=x$. Above right: construction of the right and left inverses of $x$. Below: Associativity $(xy)z=x(yz)$.}\label{F: fig27}
\end{figure}

The constant map $\Delta^n\to *$ (which we will also denote by $*$) is the identity element. Indeed, given $x\in X$ representing an element of $\pi_n(X,*)$, the $(n+1)$-simplex $s_{n}x$ will have $d_{n+1}s_nx=d_ns_nx=x$, while for $i<n$, $d_is_nx=s_{n-1}d_ix\in *$. This realizes $x=*x$. Similarly, consideration of $s_{n-1}x$  gives $x=x*$.

It is also easy to construct inverses:  given $x\in X$ representing an element of $\pi_n(X,*)$, there is no problem mapping the horn $\Lambda^{n+1}_{n+1}$ into $X$ such that the face corresponding to $d_{n-1}\Delta^{n+1}$ goes to $x$ and all other faces land in $*$. The Kan condition lets us extend this to a map of $\Delta^{n+1}$ into $X$ and then the face corresponding to $d_{n+1}\Delta^{n+1}$ is a right inverse to $x$. Similarly, we can find a left inverse using $\Lambda^{n+1}_{n-1}$ and putting $x$ on the face corresponding to $d_{n+1}\Delta^{n+1}$. 

Finally, we show that the group operation is associative, which takes a bit more work. 
Let $x,y,z$ be simplices in $X$ representing elements of $\pi_n(X,*)$. The idea is to form an $(n+2)$-horn with appropriate $(n+1)$-faces realizing $xy$, $yz$, and $(xy)z$ and then to use a Kan extension argument to show that the new face guaranteed by the extension realizes $x(yz)$ on the same $n$-face that already represents $(xy)z$.

In more detail, we choose $(n+1)$-simplices $w_{n-1}$ and $w_{n+2}$ that respectively realize the products $xy$ and $yz$, and we choose a simplex $w_{n+1}$ realizing the product $(xy)z$, where $xy$ is represented by $d_nw_{n-1}$. Now, we can find a horn $\Lambda^{n+2}_n$ in $X$ such that the faces corresponding to $d_i\Delta^{n+2}$ are the $w_i$ for $i=n-1,n+1,n+2$ and $*$ otherwise. To see that this data is consistent to form the horn, we  need to check  the appropriate faces, most of which are in $*$, to see that they correspond.  The only faces of $\Delta^{n+2}$ we don't need to check are those of the form $d_id_n\Delta^{n+2}$ since $d_n\Delta^{n+2}$ isn't in the horn. By the simplicial axioms, these also correspond to the faces  $d_{n-1}d_i\Delta^{n+2}$ for  $i<n$ and $d_nd_{i+1}\Delta^{n+2}$ for $i\geq n$. This leaves the following faces to check: We have $d_nd_{n-1}\Delta^{n+2}_n=d_nw_{n-1}=xy=d_{n-1}w_{n+1}=d_{n-1}d_{n+1}\Delta^{n+2}_n$ and $d_{n+1}d_{n-1}\Delta^{n+2}_n=d_{n+1}w_{n-1}=y=d_{n-1}w_{n+2}=d_{n-1}d_{n+2}\Delta^{n+2}_n$.  We also have  $d_{n+1}d_{n+1}\Delta^{n+2}_n=d_{n+1}w_{n+1}=z=d_{n+1}w_{n+2}=d_{n+1}d_{n+2}\Delta^{n+2}_n$.  All other sides in the proposed horn are in $*$, and so the data is consistent.

We can extend this horn to an $(n+2)$-simplex $u$ by the Kan condition. So now by definition of $w_{n+1}$,  $(xy)z=d_nw_{n+1}=d_nd_{n+1}u$, which, using the axioms, is also equal to $d_nd_nu$. But this also represents the product of $d_{n-1}d_nu=d_{n-1}d_{n-1}u=d_{n-1}w_{n-1}=x$ with $d_{n+1}d_nu=d_nd_{n+2}u=d_nw_{n+2}=yz$. So $d_nd_nu$ also represents the product $x(yz)$, proving associativity. 
\end{proof}

Also as expected, $\pi_n(X,*)$ is an abelian group for  $n\geq 2$, but this is a bit more difficult to prove. We refer the reader to \cite[Proposition 4.4]{MAY67}.

\paragraph{Relative homotopy groups.}

If $(X,A,*)$ is a Kan triple (meaning $A$ is a Kan subcomplex of the Kan complex $X$ and $*$ is a basepoint in $A$), there are also relative homotopy groups $\pi_n(X,A,*)$. Corresponding to our first definition of $\pi_n(X,*)$ and the topological notion of relative homotopy, we could define $\pi_n(X,A,*)$ to be relative homotopy classes of maps $(\Delta^n,\bd \Delta^n,[0])\to (X,A,*)$, where the homotopies are required to keep the image of $\bd \Delta^n\times I$ in $A$ and the image of $[0]\times I$ in $*$. For a  version of $\pi_n(X,A,*)$ corresponding to our second definition of $\pi_n(X,*)$, we first need a relative notion of homotopy of simplices:

\begin{definition}\label{D: rel homotopy of simplices}
If $A$ is a subcomplex of $X$, we say that two $n$-simplices $x,x'\in X_n$ are \emph{homotopic rel $A$} if $d_ix=d_ix'$ for $1\leq i\leq n$, $d_0x$ is homotopic to $d_0x'$ in $A$ via an $n$-simplex $y$, and there exists a simplex $w\in X_{n+1}$ such that $d_0w=y$, $d_nw=x$, $d_{n+1}w=x'$, and $d_iw=s_{n-1}d_ix=s_{n-1}d_ix'$, $1\leq i\leq n-1$.  
\end{definition}

This definition is very similar to that for homotopy of simplices except instead of requiring $d_0x=d_0x'$, we let $d_0x$ and $d_0x'$ be  two simplices that are themselves homotopic in $A$, and the homotopy between $x$ and $x'$, provided by $w$, contains within it the homotopy between $d_0x$ and $d_0x'$. 

Using this relative notion of homotopy, we can define $\pi_n(X,A,*)$.

\begin{definition}
Given a Kan triple $(X,A,*)$, we define $\pi_n(X,A,*)$, $n>0$, as the set of equivalence classes of $n$-simplices $x\in X$ with $d_0x\in A$ and $d_ix\in *$ for all $i$, $1\leq i\leq n$, up to relative homotopy of simplices. 
\end{definition}

$\pi_n(X,A,*)$ is also a group for $n\geq 2$ and an abelian group for $n\geq 3$. We will define the product; the proofs of well-definedness and that we have a group are analogous to those for $\pi_n(X,*)$.

\begin{figure}[!htp]
\begin{center}
\scalebox{.6}{\includegraphics{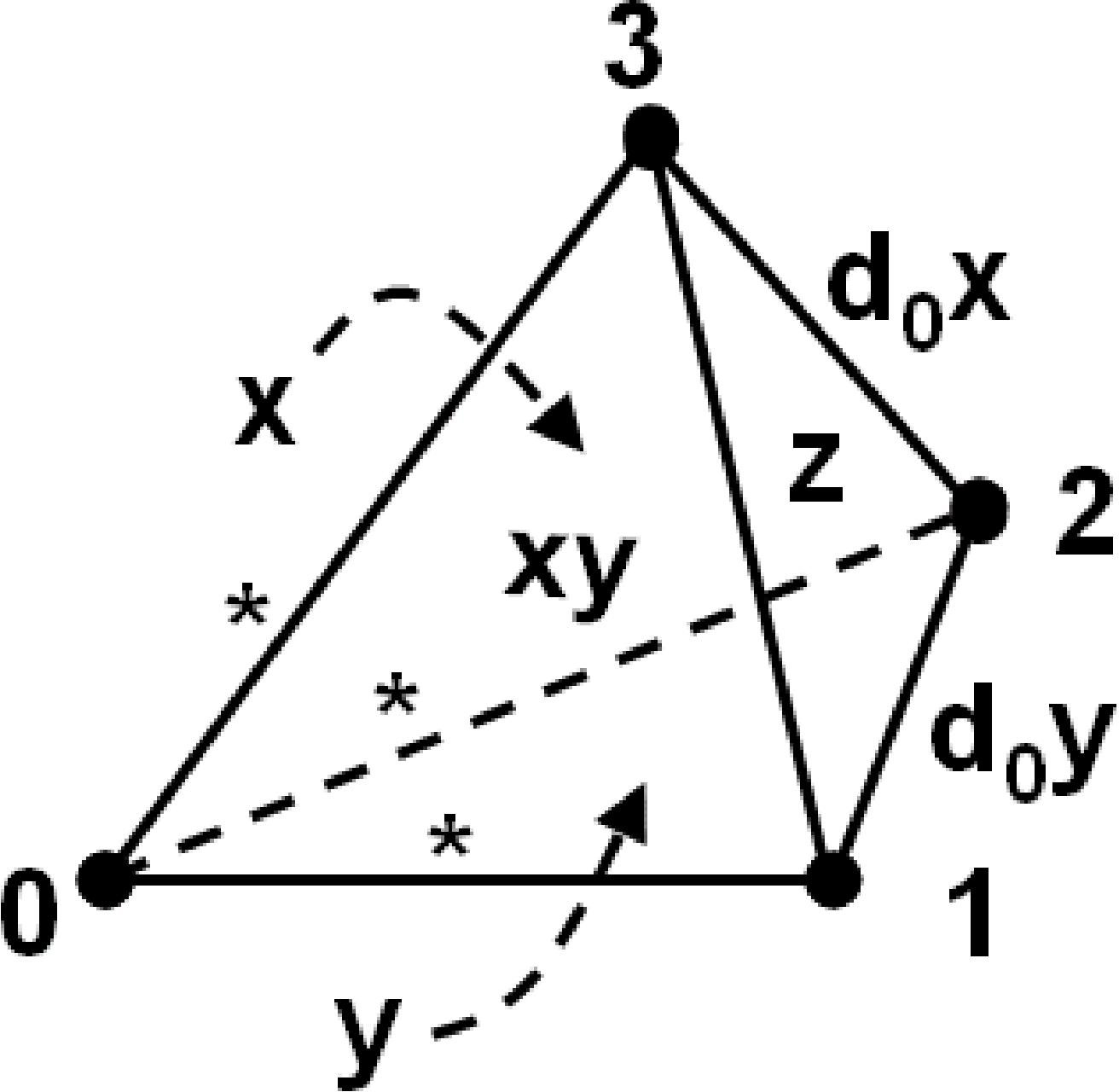}}
\end{center}
\caption{The product of two elements of $x,y\in \pi_2(X,A,*)$. The $1$-simplex with endpoints $1$ and $3$ is the product $(d_0x)(d_0y)$ in $\pi_1(A)$.}\label{F: fig28}
\end{figure}

\begin{definition}
Suppose $x,y$ represent elements of $\pi_n(X,A,*)$, $n\geq 2$. Let $z$ represent the product between $d_0x$ and $d_0y$ in $\pi_{n-1}(A,*)$; in other words, let $z\in A_{n+1}$ be such that $d_{n-2}z=d_0x$ and $d_{n}z=d_0y$ so that $d_{n-1}z=(d_0x)(d_0y)$. Now map the horn $\Lambda^{n+1}_n$ into $X$ such that the sides corresponding to $d_0\Delta^{n+1}$, $d_{n-1}\Delta^{n+1}$, and $d_{n+1}\Delta^{n+1}$ are $z$, $x$, and $y$, respectively, and all other faces go to $*$. One can check that this is consistent data. Then let $w$ be an extension of the horn, which exists because $X$ is Kan, and define $xy$ to be $d_nw$.  See Figure \ref{F: fig28}.
\end{definition}

An excellent exercise for the reader at this point would be to show that there is a long exact sequence 

$$\cdots \to\pi_n(A,*)\to\pi_n(X,*)\to\pi_n(X,A,*)\to \pi_{n-1}(A,*)\to\cdots .$$

\section{Concluding remarks}\label{S: concluding}

It is difficult to know where to end a survey of the type we have undertaken here. On the one hand, although we have included some material from its later chapters,
we have not even covered the entire first chapter of May's textbook \cite{MAY67}! On the other hand, our goal has never been to provide a completely rigorous or comprehensive treatise on simplicial theory, but to provide the reader with an introduction to some of the most important elementary ideas while maintaining a bridge to the geometric pictures that the combinatorics are based upon. We hope that we have prepared the interested student to move on to the more standard texts on simplicial objects with some picture (literally) of what's going on there. 

And what is going on there? Just about everything in topological homotopy theory and then some. Just a glance at the table of contents of \cite{MAY67} turns up many familiar concepts from homotopy theory: fibrations, fiber bundles, Postnikov systems, function spaces, Hurewicz theorems, Eilenberg-Mac Lane complexes, $k$-invariants, cup and cap products, the Serre spectral sequence, etc. This is not surprising in light of the following theorem; we refer the reader to Curtis \cite[Section 12]{Cu71}, or to \cite[Section I.11]{GoeJar} for a modern proof.

\begin{theorem}\label{T: equiv}
The homotopy category of Kan complexes, consisting of Kan complexes and homotopy classes of maps between them, is equivalent to the category of CW complexes and homotopy classes of continuous maps.
\end{theorem}

The functors that realize this equivalence are the realization functor of simplicial complexes and the singular set functor that assigns the singular set to a topological space. Thus this theorem is closely related to the adjunction theorem, Theorem \ref{T: adjoint}. So this tells us that everything we have been doing in the simplicial realm is a reflection of ordinary homotopy theory. Yet, despite the geometric point of view we have been emphasizing here, simplicial theory is purely combinatorial and algebraic, accessible by discrete tools that may not be evident in pure topology. Thus, using simplicial theory, one can hope to study topological homotopy theory via these combinatorial tools. Furthermore, we touched upon how the combinatorial simplicial methods can be transported to other contexts, such as simplicial groups. They can also be abstracted to broader categorical settings, leading to the theory of simplicial model categories. We hope to have introduced enough of the background also to enable the reader to pursue these more modern approaches, such as can be found in \cite{GoeJar}, with some understanding of their original motivation in concrete homotopy theory.

We leave the reader with some bibliographical notes on the sources we have used. 

Our primary sources were May's \emph{Simplicial Objects in Algebraic Topology} \cite{MAY67} and Moore's lecture notes \emph{Seminar on algebraic homotopy theory} \cite{MOORE}. May's book, first published in 1967, is the most comprehensive reference of its time, featuring a direct combinatorial approach. Moore's notes are from nearly a decade earlier, but they are perhaps a bit more accessible to the geometrically-minded reader; they take a different approach to homotopy groups, defining them as $\pi_0$ of simplicial loop spaces. Our primary modern source was \emph{Simplicial Homotopy Theory} \cite{GoeJar} by Goerss and Jardine. It starts off directly from the modern model category point of view, without much need for the combinatorial underpinnings (some knowledge of the combinatorial approach, however, will aid the reader). Despite the abstractness of the material, I found this book quite readable. The book \emph{Calculus of Fractions and Homotopy Theory} \cite{GabZis} by Gabriel and Zisman, though contemporary with May's book, is something of a bridge between the classical combinatorics and some of the more current axiomatic ideas. We should also mention in this paragraph the long survey \emph{Simplicial Homotopy Theory} \cite{Cu71} by Curtis, and since initially writing this exposition I have become aware of another introductory survey paper \emph{Introduction to Combinatorial Homotopy Theory} by Sergeraert \cite{Ser08}. 
As one might expect, each of these sources contains somewhat different material and sometimes different approaches to the same material, thus it is well worth consulting each of them depending on the reader's interests  in terms of both  material and style.

Besides these longer expositions,  introductory chapters on simplicial theory can be found within many other textbooks and surveys. In particular, I know of sections on simplicial theory 
in Selick's \emph{Introduction to Homotopy Theory} \cite{SELICK},
Smirnov's \emph{Simplicial and Operad Methods in Algebraic Topology} \cite{SMIRNOV}, Gelfand and Manin's \emph{Methods of Homological Algebra} \cite{GelMan2}, and Weibel's \emph{An Introduction to Homological Algebra} \cite{WEIB}. As one might expect, these last references are a good source for applications of simplicial theory to homological algebra.  There are also  review sections on simplicial sets in Bousfield and Kan's \emph{Homotopy Limits, Completions, and Localizations} \cite{BouKan} and in \emph{Mixed Hodge Structures} \cite{PetSte} by Peters and Steenbrink.
The breadth of topics covered by those titles alone should give the reader some impression of just how varied the applications of simplicial theory are.

\providecommand{\bysame}{\leavevmode\hbox to3em{\hrulefill}\thinspace}
\providecommand{\MR}{\relax\ifhmode\unskip\space\fi MR }
\providecommand{\MRhref}[2]{%
  \href{http://www.ams.org/mathscinet-getitem?mr=#1}{#2}
}
\providecommand{\href}[2]{#2}

\pagebreak
\begin{center}
{\Large
Errata to

An elementary illustrated introduction to simplicial sets 

Rocky Mountain Journal of Mathematics 42 (2012), 353-424 
}
\end{center}

\vskip.75in

The following corrections have been made to this version of the paper but remain as errata in the published version. Thanks to Jim Davis, Donghan Wang, Daniel de Carvalho, Abdullah Malik, and Marie Labeye for each pointing out errors. 

\begin{enumerate}
\item Section 2.1, ``the intersection of any two simplices of $X$ is a face of each them'' corrected to ``the intersection of any two simplices of $X$, if non-empty, is a face of each them.''

\item Example 3.7, grammatical error corrected

\item Examples 4.4 and 4.7, specified \emph{ordered} simplicial complex

\item In Example 5.4, corrected to have $d_0([0,1],[0,1])=(1,1)$ and $d_1([0,1],[0,1])=(0,0)$, not the other way around.

\item Penultimate paragraph of Section 5.5, some $S_i$ corrected to $P_i$. 

\item In the last paragraph of Sectin 5.5:

\begin{enumerate}

\item the expression 
$$d_iP_k=(s_{k-1}d_i E_p, s_{p-1}\cdots  s_{k+1}s_{k-1}\cdots s_i (d_is_i)s_{i-1}\cdots s_{0} e)=(s_{k-1}d_i E_p, s_{p-1}\cdots s_{k}s_{k-2} \cdots s_{0} e)$$
has been corrected to
$$d_iP_k=(s_{k-1}d_i E_p, s_{p-1}\cdots  s_{k}s_{k-2}\cdots s_i (d_is_i)s_{i-1}\cdots s_{0} e)=(s_{k-1}d_i E_p, s_{p-1}\cdots s_{k}s_{k-2} \cdots s_{0} e)$$

\item the expression $$d_iP_k=(s_{k}d_{i-1} E_p, s_{p-1}\cdots s_{k}s_{k-2} \cdots s_{0} e)$$
has been corrected to $d_iP_k=(s_{k}d_{i-1} E_p, s_{p-1}\cdots s_{k+1}s_{k-1} \cdots s_{0} e)$

\item the expressions
\begin{align*}
d_kP_k&=(d_ks_k E_p, s_{p-1}\cdots s_{k+1}s_{k-2} \cdots s_{0} e)=(E_p, s_{p-1}\cdots s_{k+1}s_{k-2} \cdots s_{0} e)\\
d_{k+1}P_k&=(d_{k+1}s_k E_p, s_{p-1}\cdots s_{k+2}s_{k-1} \cdots s_{0} e)=(E_p, s_{p-1}\cdots s_{k+2}s_{k-1} \cdots s_{0} e).
\end{align*}

have been corrected to

\begin{align*}
d_kP_k&=(d_ks_k E_p, s_{p-1}\cdots s_{k}s_{k-2} \cdots s_{0} e)=(E_p, s_{p-1}\cdots s_{k}s_{k-2} \cdots s_{0} e)\\
d_{k+1}P_k&=(d_{k+1}s_k E_p, s_{p-1}\cdots s_{k+1}s_{k-1} \cdots s_{0} e)=(E_p, s_{p-1}\cdots s_{k+1}s_{k-1} \cdots s_{0} e).
\end{align*}
\end{enumerate}

\item In Example 6.4, corrected to note that $BG$ is a simplicial group only if $G$ is abelian. 

\item Prior to Example 7.5, there was a $K$ that has been corrected to an $X$.

\item In Example 7.6, there has been a font correction from $S$ to $\ms S$

\item Statement of Lemma 9.5 corrected to note that the boundaries of the simplices should all be in the basepoint. 

\item Paragraph after Definition 9.2, ``she'' corrected to ``the''

\item Figure 28, the edge labeled $w'$ was previously incorrectly identified as $w$. The label has been corrected and the caption updated.

\item Last paragraph in the proof of Lemma 9.5, $\bd \Delta^n\times \Delta^0$ corrected to $\bd \Delta^n\times \Delta^1$

\item Last paragraph in the proof of Lemma 9.6, expressions $H_{k+1}$ corrected to $H_{n+1}$ 

\item Definition 9.13, reference added to Figure \ref{F: fig28}
\end{enumerate}

\end{document}